\documentclass[10pt]{article}
\usepackage{psfrag,amsmath}
\usepackage{amssymb}
\usepackage{mathrsfs}
\usepackage{amsfonts}
\usepackage{amsmath}
\usepackage{mathrsfs,color}
\usepackage{epstopdf}
\usepackage{graphicx}
\usepackage{mathrsfs,amscd,amssymb,amsthm,amsmath,bm,url}
\usepackage{tikz}
\usetikzlibrary{positioning,backgrounds}
\usetikzlibrary{fadings}
\usetikzlibrary{patterns}
\usetikzlibrary{calc}
\usetikzlibrary{shadings}
\pgfdeclarelayer{background}
\pgfdeclarelayer{foreground}
\pgfsetlayers{background,main,foreground}

\makeatletter

\renewcommand{\@seccntformat}[1]{{\csname the#1\endcsname}{\normalsize .}\hspace{.5em}}
\makeatother

\def \[{\begin{equation}}
\def \]{\end{equation}}

\def \diam{{\rm diam}}

\def \ex{{\rm ex}}
\newtheorem{thm}{Theorem}

\newtheorem{claim}[thm]{Claim}

\newtheorem{lem}[thm]{Lemma}

\newtheorem{pb}{Problem}

\newenvironment{wst}
{\setlength{\leftmargini}{1.5\parindent}
 \begin{itemize}
 \setlength{\itemsep}{-1.1mm}}
{\end{itemize}}
\begin{document}

\setlength{\baselineskip}{0.20in}
\begin{center}{\Large\bf {Extensions on spectral extrema of $C_5/C_6$-free graphs with given size}\footnote{Financially supported  by the National Natural Science Foundation of China (Grant Nos. 12171190, 11671164) and the excellent doctoral dissertation cultivation grant from Central China Normal University (Grant No. 2022YBZZ033)}}
\vspace{4mm}

{\large Wanting Sun$^a$,\ \ Shuchao Li$^{a,}$\footnote{Corresponding author. \\
\hspace*{5mm}{\it Email addresses}: wtsun2018@sina.com (W.T. Sun), \ lscmath@ccnu.edu.cn (S.C. Li), \ weiweimath@sina.com (W. Wei).},\ \ Wei Wei$^b$}\vspace{2mm}

$^a$Hubei Key Laboratory of Mathematical Science, and Faculty of Mathematics and Statistics,\\  Central China Normal
University, Wuhan 430079, PR China

$^b$School of Mathematics, Physics and Statistics, Shanghai University of Engineering Science,\\ Shanghai 201620, PR China
\end{center}

\noindent {\bf Abstract}:\ Let $\mathcal{F}$ denote a set of graphs. A graph $G$ is said to be $\mathcal{F}$-free if it does not contain any element of $\mathcal{F}$ as a subgraph. The Tur\'an number is the maximum possible number of edges in an $\mathcal{F}$-free graph with $n$ vertices. It is well known that classical Tur\'an type extremal problem aims to study the Tur\'an number of fixed graphs. In 2010, Nikiforov \cite{Nik2} proposed analogously a spectral Tur\'an type problem which asks to determine the maximum spectral radius of an $\mathcal{F}$-free graph with $n$ vertices. It attracts much attention and many such problems remained elusive open even after serious attempts, and so they are considered as one of the most intriguing problems in spectral extremal graph theory. It is interesting to consider another spectral Tur\'an type problem which asks to determine the maximum spectral radius of an $\mathcal{F}$-free graph with $m$ edges.  Denote by $\mathcal{G}(m,\mathcal{F})$ the set of $\mathcal{F}$-free graphs with $m$ edges having no isolated vertices. Each of the graphs among $\mathcal{G}(m,\mathcal{F})$ having the largest spectral radius is called a maximal graph. Let $\theta_{p,q,r}$ be a theta  graph  formed  by  connecting  two distinct vertices with three independent paths of length $p,q$ and $r,$ respectively (length refers to the number of edges). In this paper, we firstly determine the unique maximal graph among $\mathcal{G}(m,\theta_{1,2,3})$ and  $\mathcal{G}(m,\theta_{1,2,4}),$ respectively. Then we determine all the maximal graphs among $\mathcal{G}(m,C_5)$ (resp. $\mathcal{G}(m,C_6)$) excluding the book graph. These results extend some earlier results.

\vspace{2mm} \noindent{\it Keywords:}
Tur\'an-type extremal problem; $\mathcal{F}$-free graphs; Spectral radius
\vspace{2mm}


\section{\normalsize Introduction}

For a simple graph $G=(V(G),E(G)),$ we use $n:=|V(G)|$ and $m:=|E(G)|$ to denote the order and the size of $G,$ respectively. Since isolated vertices do not have an effect on the spectral radius, throughout this paper we consider graphs without isolated vertices. Let $\mathcal{F}$  be a set of graphs,  we say that $G$ is \textit{$\mathcal{F}$-free} if it does not contain any element in $\mathcal{F}$ as a subgraph. When the forbidden set $\mathcal{F}$ is a singleton, say $\{F\}$, then we write $F$-free for $\mathcal{F}$-free.

In 2013, F\"{u}redi and Simonovits \cite{survey} posed the following problem:
\begin{pb}[F\"{u}redi-Simonovits type problem]\label{pb1}
Assume $\mathbb{U}$ is a family of graphs and $G$ is in $\mathbb{U}.$ For a specific pair of parameters $(\tau,\upsilon)$ on $G$, our aim is to maximize the second parameter $\upsilon$ under the condition that $G$ is $\mathcal{F}$-free and its first parameter $\tau$ is given.
\end{pb}
If the pair of parameters above are the order and size of a graph, i.e., $(\tau,\upsilon)=(n,m),$ then the F\"{u}redi-Simonovits type problem is just the classical Tur\'{a}n type problem: determine the maximum number of edges, $\ex(n,\mathcal{F})$, of an $n$-vertex $\mathcal{F}$-free graph. The value $\ex(n,\mathcal{F})$ is called the \textit{Tur\'an number}. The research for the Tur\'an number attracts much attention, and it has become to be one of the most attractive fundamental problems in extremal graph theory (see \cite{survey,V2011} for surveys).

Let $A(G)$ be the adjacency matrix of a graph $G$. The largest modulus of all eigenvalues of $A(G)$ is the \textit{spectral radius} of $G$ and denoted by $\lambda(G)$. As usual, let $P_n,\,C_n,\,K_n,\,K_{a,n-a}$ and $K_{1,n-1}$ be the path, the cycle, the complete graph, the complete bipartite graph and the star on $n$ vertices, respectively.

In Problem~\ref{pb1}, if one lets $(\tau,\upsilon)=(n,\lambda(G))$, i.e., the pair of parameters are the order and the spectral radius on $\mathbb{U},$ then it becomes to be the spectral Tur\'{a}n type problem (also known as Brualdi-Solheid-Tur\'{a}n type problem, see \cite{Nik2}): what is the maximal spectral radius of an $\mathcal{F}$-free graph with order $n ?$ Over the past decade, much attention has been paid to the Brualdi-Solheid-Tur\'{a}n type problem. For more details, one may consult the references, such as for $\mathcal{F}= \{K_r\}$ \cite{Nik3,Wilf}, $\mathcal{F}= \{K_{s,t}\}$ \cite{Bab,Nik3,Nik4}, $\mathcal{F}= \{P_k\}$ \cite{Nik2}, $\mathcal{F}= \{C_4\}$ \cite{Nik6,ZW} and $\mathcal{F}= \{C_6\}$ \cite{zhai}.

In Problem~\ref{pb1}, if one lets $(\tau,\upsilon)=(m,\lambda(G))$, i.e., the pair of parameters are the size and the spectral radius on $\mathbb{U},$ then it becomes to be another spectral Tur\'{a}n type problem (also known as Brualdi-Hoffman-Tur\'an type problem, see \cite{Brua}): what is the maximal spectral radius of an $\mathcal{F}$-free graph of given size $m$? Each of the corresponding extremal graphs is called the \textit{maximal graph}. Nosal \cite{2} showed that every triangle-free graph $G$ with $m$ edges satisfies $\lambda(G)\leqslant \sqrt{m}.$ Very recently, Lin, Ning and Wu \cite{lin} slightly improved the bound to $\lambda(G)\leqslant \sqrt{m-1}$ when $G$ is non-bipartite and triangle-free of size $m$. Nikiforov \cite{Nik,Nik6} identified the unique graph with fixed size which attains the maximum spectral radius among $K_{r+1}$-free graphs. Very recently, Zhai, Lin and Shu \cite{ZLS} studied the Brualdi-Hoffman-Tur\'an type problem on $C_5$-free/$C_6$-free graphs with given size. For more results on this topic, we refer the reader to see \cite{Cheng,Nik5,ZS}. It is consequently the aim of this article to make some contribution on Brualdi-Hoffman-Tur\'an type problem.

A \textit{generalized theta graph}, say $\theta_{l_1, l_2,\ldots,l_t},$ is the graph obtained by connecting two vertices with $t$ internally disjoint paths of lengths $l_1, l_2,\ldots,l_t,$ where $l_1\leqslant l_2\leqslant \cdots\leqslant l_t$ and $l_2\geqslant 2.$ With the restriction $t=3,$ we obtain a \textit{theta graph} $\theta_{p, q, r}.$ Clearly, $\theta_{l_1,l_2}\cong C_{l_1+l_2}.$ And so, the research on F\"{u}redi-Simonovits type problems with $\theta_{l_1, l_2,\ldots,l_t}$ as the forbidden subgraph is a natural generalization of those on cycle-free graphs. For more advances along this line, we refer the reader to \cite{Bukh,Fa,Vers}.

Let $\mathcal{G}(m,\mathcal{F})$ denote the set of $\mathcal{F}$-free graphs with $m$ edges having no isolated vertices. {If $\mathcal{F}=\{F\}$, then we write $\mathcal{G}(m,F)$ for $\mathcal{G}(m,\mathcal{F})$.} The \textit{join} of simple graphs $G$ and $H,$ written $G\vee H,$ is the graph obtained from the disjoint union $G\cup H$ by adding the edges $\{xy:x\in V(G),y\in V(H)\}.$ Let $S_{n,k}:=K_k\vee (n-k)K_1$ be the \textit{complete split graph}.

In this paper, we first characterize the unique graph with the maximum spectral radius among $\mathcal{G}(m,\theta_{1,2,3})$ and $\mathcal{G}(m,\theta_{1,2,4}),$ respectively.
\begin{thm}\label{thm1.04}
Let $G$ be a graph in $\mathcal{G}(m,\theta_{1,2,r}).$ Then the following holds.
 \begin{wst}
 \item[{\rm (i)}] If $r=3$ and $m\geqslant 8,$ then $\lambda(G)\leqslant \frac{1+\sqrt{4m-3}}{2}$ and equality holds if and only if $G\cong S_{\frac{m+3}{2},2};$
 \item[{\rm (ii)}] If $r=4$ and $m\geqslant 22,$ then $\lambda(G)\leqslant \frac{1+\sqrt{4m-3}}{2}$ and equality holds if and only if $G\cong S_{\frac{m+3}{2},2}.$
 \end{wst}
\end{thm}

{Notice that $\theta_{1,2,r}$ can be viewed as a graph obtained from a cycle $C_{r+2}$ by adding one edge between two vertices with distance two. Hence, if a graph $G$ is $C_{r+1}$-free, then it is $\theta_{1,2,r}$-free. Notice that $S_{\frac{m+3}{2},2}$ is $\{C_5,C_6\}$-free. Then by Theorem \ref{thm1.04}, we have the following result, which was obtained by Zhai, Lin and Shu \cite{ZLS}}.
\begin{thm}[\cite{ZLS}]\label{thm004}
Let $G$ be a graph of size $m.$ If $G\in \mathcal{G}(m,C_5)$ with $m\geqslant 8$, or $G\in \mathcal{G}(m,C_6)$ with $m\geqslant 22,$ then $\lambda(G)\leqslant \frac{1+\sqrt{4m-3}}{2}.$ Equality holds if and only if $G\cong S_{\frac{m+3}{2},2}.$
\end{thm}

We shall notice that the size of the unique maximal graph in Theorem \ref{thm004} is odd. In what follows we use a unified approach (regardless of the parity of $m$) to determine all the graphs with maximum spectral radius among $\mathfrak{G}(m,C_5):=\mathcal{G}(m,C_5)\setminus \{S_{\frac{m+3}{2},2}\}$ or $\mathfrak{G}(m,C_6):=\mathcal{G}(m,C_6)\setminus \{S_{\frac{m+3}{2},2}\}$ for $m\geqslant 22.$

Let $S_{n,k}^t$ be the graph obtained from $S_{{n-t},k}$ by attaching $t$ pendant vertices to the maximum degree vertex of $S_{{n-t},k}$, and  $S_n^k$ be the graph obtained from $K_{1,n-1}$ by adding $k$ disjoint edges within its independent set. Let $\rho_1(m)$ be the largest zero of $\psi_1(x)$, where
\begin{equation*}
    \psi_1(x)=\left\{
    \begin{array}{ll}
        x^4-m x^2-(m-2)x+\frac{m}{2}-1,& \textrm{if $m$ is even;}\\[5pt]
        x^4-m x^2-(m-3)x+m-3,& \textrm{if $m$ is odd.}
    \end{array}
    \right.
\end{equation*}
and let $\rho_2(m)$ be the largest zero of $\psi_2(x),$ where
 \begin{equation*}
    \psi_2(x)=
    \left\{
    \begin{array}{ll}
        x^3-2 x^2-(m-3)x+m-6,& \textrm{if $m$ is even and $22\leqslant m\leqslant 72;$}\\[5pt]
        x^4-m x^2-(m-2)x+\frac{m}{2}-1,& \textrm{if $m$ is even and $m\geqslant 74;$}\\[5pt]
        x^5-x^4-(m-1)x^3-2x^2+\frac{3m-17}{2}x-\frac{m-7}{2},& \textrm{if $m$ is odd and $23\leqslant m\leqslant 71;$}\\[5pt]
        x^4-m x^2-(m-3)x+m-3,& \textrm{if $m$ is odd and $m\geqslant 73.$}
    \end{array}
    \right.
\end{equation*}

The subsequent two results identify all the graphs having the largest spectral radii among $\mathfrak{G}(m,C_5)$ and $\mathfrak{G}(m,C_6)$, respectively.
\begin{thm}\label{thm1.3}
 Let ${G}$ be in $\mathfrak{G}(m,C_5)$ with $m\geqslant 22.$ Then $\lambda(G)\leqslant \rho_1(m).$
Equality holds if and only if $G\cong S_{\frac{m+4}{2},2}^1$ if $m$ is even and $G\cong S_{\frac{m+5}{2},2}^2$ if $m$ is odd.
\end{thm}
\begin{thm}\label{thm1.03}
 Let ${G}$ be in $\mathfrak{G}(m,C_6)$ with $m\geqslant 22.$ Then $\lambda(G)\leqslant \rho_2(m).$
Equality holds if and only if
 \begin{equation*}
    G\cong
    \left\{
    \begin{array}{ll}
    K_1\vee S_{\frac{m}{2}}^1,& \textrm{if $m$ is even and $22\leqslant m\leqslant 72;$}\\[5pt]
    S_{\frac{m+4}{2},2}^1,& \textrm{if $m$ is even and $m\geqslant 74;$}\\[5pt]
    K_1\vee (S_{\frac{m-1}{2}}^1\cup K_1),& \textrm{if $m$ is odd and $23\leqslant m\leqslant 71;$}\\[5pt]
    S_{\frac{m+5}{2},2}^2,& \textrm{if $m$ is odd and $m\geqslant 73.$}
    \end{array}
    \right.
\end{equation*}
\end{thm}

In Theorems \ref{thm1.3} and \ref{thm1.03}, if we let $m\geqslant 22$ be an even integer, then the corresponding graphs with maximum spectral radius among $\mathcal{G}(m,C_5)$ (resp. $\mathcal{G}(m,C_6)$) coincide with \cite[Theorem 2.1]{Lou}.

{The remainder of the paper is organized as follows:} In Section 2, we recall some important known results. In Section 3 we give the proof of Theorem~\ref{thm1.04}. In Section~4 we give the proofs of Theorems~\ref{thm1.3} and \ref{thm1.03}. In the last section, we give some brief comments on our contribution and propose some further research questions.
\section{\normalsize Preliminaries}
In this section, we describe some preliminary results, which plays an important role in the subsequent sections.  Unless otherwise stated, we follow the traditional notation and terminology; see \cite{0001,0009}. In order to formulate these results, we need some additional notation.

For a graph $G$ with a vertex subset $S\subseteq V(G)$, denote by $G[S]$ the subgraph of $G$ induced by $S.$ Let $e(G):=|E(G)|$ be the size of $G.$ For two vertex subsets $S$ and $T$ of $G$ (where $S\cap T$ may not be empty), let $e(S, T)$ denote the number of edges with one endpoint in $S$ and the other in $T$, and $e(S, S)$ is simplified by $e(S).$ For two distinct vertices $u,v\in V(G)$ with $uv\not\in E(G),$ we write $G+uv:=(V(G),E(G)\cup \{uv\}).$

For a vertex $v\in V(G)$, let $N_G(v)$ be the neighborhood of $v$ in $G$, and $N_G[v]:=N_G(v)\cup \{v\}$ be the closed neighborhood of $v$ in $G.$ Denote by {$d_G(v):=|N_G(v)|$} the \textit{degree} of $v$ in $G$. Here, as elsewhere, we drop the index referring to the underlying graph if the reference is clear. Let $N_S(v):=N_G(v)\cap S$ and $d_S(v):=|N_S(v)|.$ {For each nonnegative integer $i,$ put} $N_i(v):=\{u: u\in N(v),\,u\ \text{lies in the component}\ K_{i+1}\ \text{of}\ G[N(v)]\}$ and $N^2(v):=(\cup_{u\in N(v)}N(u))\setminus N[v].$ {In addition, for a subgraph $H$ of $G$ and a vertex $u\in V(G)\setminus V(H),$ denote $N_H(u):=N_G(u)\cap V(H)$ and $d_H(u):=|N_H(u)|.$} 

Let $G$ be a connected graph. Then $A(G)$ is irreducible and nonnegative. From the Perron-Frobenius Theorem, we know that the largest eigenvalue of $A(G)$ is the spectral radius of $G$ and there exists a unique positive unit eigenvector ${\bf x}$ of $A(G)$ corresponding to $\lambda(G),$ which is called the \textit{Perron vector} of $G.$ It will be convenient to associate a labeling of vertices of $G$ (with respect to ${\bf x}$) in which $x_r$ is a label of the vertex $r$. If $x_{{u}} = \max\{x_v: v\in V({G})\},$ then we call $u$ an \textit{extremal vertex} in ${G}.$ We say that two distinct vertices $u$ and $v$ are \textit{equivalent} in $G,$ if there exists an automorphism $\phi:G\rightarrow G$ such that $\phi(u)=v.$ It is well known that $x_u=x_v$ if $u$ and $v$ are equivalent vertices in $G$ (see \cite[Proposition 16]{0007}). This result will be frequently used in our proofs. 

\begin{lem}[\cite{Nik,Nik6}]\label{thm4}
If $G\in \mathcal{G}(m,K_{r+1})$ for some integer $r\geqslant 2,$ then $\lambda(G)\leqslant \sqrt{2m(1-1/r)}.$ Equality holds if and only if $G$ is a complete bipartite graph for $r=2$, and $G$ is a complete regular $r$-partite graph for $r\geqslant 3.$
\end{lem}
\begin{lem}[\cite{Nik6}]\label{lem1.1}
Let $G$ and $G'$ be two connected graphs on the same vertex set. Assume that $N_G(u)\subsetneqq N_{G'}(u)$ for some vertex $u$. If the Perron vector ${\bf x}$ of $G$ satisfies ${\bf x}^TA(G'){\bf x}\geqslant {\bf x}^TA(G){\bf x},$ then $\lambda(G')>\lambda(G).$ In particular, assume that $u$ and $v$ are two distinct vertices of $G$ with $x_u\geqslant x_v$ and $\{v_i:1\leqslant i\leqslant s\}\subseteq N_G(v)\setminus N_{G}[u].$ If $G''=G-\{vv_i:1\leqslant i\leqslant s\}+\{uv_i:1\leqslant i\leqslant s\},$ then $\lambda(G'')>\lambda(G).$
\end{lem}
\begin{lem}[\cite{Sun}]\label{lem4.0}
Let $G$ be a graph and let $v$ be in $V(G)$ with $d_G(v)\geqslant 1.$ Then $\lambda(G)\leqslant \sqrt{\lambda^2(G-v)+2d_G(v)-1}.$ Equality holds if and only if either $G\cong K_n$ or $G\cong K_{1,n-1}$ with $d_G(v)=1.$
\end{lem}
\begin{lem}[\cite{Pap}]\label{lem:4.0}
Let $G$ be a connected graph with order $n$, and let ${\bf x}=(x_{1},\ldots,x_{n})^T$ be the Perron vector of $G.$ Then $x_{i}\leqslant \frac{1}{\sqrt{2}}$ for $1\leqslant i\leqslant n.$
\end{lem}
Let $G$ be a graph and $H_{t,s}$ be a bipartite graph with bipartition {$T\cup S$} satisfying $|T|=t$ and $|S|=s.$ Define $H_{t,s}\circ G$ to be the graph obtained by joining a maximum degree vertex of $G$ and {each one in the vertex subset $T$ of $H_{t,s}$ with an edge.} For a nonnegative integer $k,$ let $R_k$ be the graph obtained from $k$ copies of $K_4$ by sharing a common vertex.
\begin{lem}[\cite{ZLS}]\label{lem5.2}
  If $k\geqslant 1$ and $m=6k+t\geqslant 8,$ then $\lambda(H_{t,0}\circ R_k)<\frac{1+\sqrt{4m-3}}{2}.$
\end{lem}
Let $H$ be a real matrix, whose columns and rows are indexed by $[n]:=\{1,2,\ldots,n\}.$ Assume that $\pi:=V_1\cup V_2\cup \ldots\cup V_t$ is a partition of $[n]$. Then $H$ can be partitioned with respect to $\pi$ as
\begin{equation*}
    H=\left(
        \begin{array}{ccc}
          H_{11} & \cdots & H_{1t} \\
          \vdots & \ddots & \vdots \\
          H_{t1} & \cdots & H_{tt} \\
        \end{array}
      \right),
\end{equation*}
where $H_{ij}$ denotes the submatrix of $H,$ indexed by the rows and columns of $V_i$ and $V_j$ respectively. Let $\pi_{ij}$ be the average row sum of $H_{ij}$ for $1\leqslant i,j \leqslant t.$ As usual the matrix $H_{\pi}=(\pi_{ij})$ is called the \textit{quotient matrix} of $H$. Moreover, if the row sum of $H_{ij}$ is constant for $1\leqslant i,j \leqslant t,$ then we call $\pi$ an \textit{equitable partition} of $H$.
\begin{lem}[\cite{YLH}]\label{lem3.02}
Let $H$ be a real nonnegative matrix with an equitable partition $\pi$. Then the largest eigenvalues of $H$ and $H_{\pi}$ are equal.
\end{lem}
A \textit{cut vertex} of a graph is a vertex whose deletion increases the number of components. A graph is called $2$-\textit{connected}, if it is a connected graph without cut vertices. A \textit{block} is a maximal $2$-connected subgraph of a graph. An \textit{end-block} is a block containing at most one cut vertex.
\begin{lem}[\cite{ZLS}]\label{lem5.1}
Let $G$ be a graph in $\mathcal{G}(m,F)$ such that $\lambda(G)$ is as large as possible, {where $F$ is $2$-connected. The following claims hold.
\begin{wst}
\item[{\rm (i)}] $G$ is connected;
\item[{\rm (ii)}] If $u$ is the extremal vertex of $G,$ then there is no cut vertex in $V(G)\setminus \{u\},$ and hence $d(v)\geqslant 2$ for all $v\in V(G)\setminus N[u];$
\item[{\rm (iii)}] If $F$ is $C_4$-free, then, for any two non-adjacent vertices $v_1$ and $v_2$ of degree two in $G$, one has $N_G(v_1)=N_G(v_2)$.
\end{wst}}
\end{lem}
\begin{lem}\label{lem2.1}
Let $G$ be a graph having the maximum spectral radius among $\mathfrak{G}(m,F),$ where {$m\geqslant 3$ and} $F\in \{C_5,C_6\}.$ Then $G$ is connected. Furthermore, if $u$ is the extremal vertex of $G,$ then there exists no cut vertex in $V(G)\setminus\{u\},$ and so $d(v)\geqslant 2$ for all $v\in V(G)\setminus N[u].$
\end{lem}
\begin{proof}
{Here, we only give the proof for $F=C_5$, the rest case can be proved by a similar discussion, whose procedure is omitted here.}

At first, suppose to the contrary that $G$ is disconnected. Let $G_1$ and $G_2$ be two connected components of $G$ and ${\lambda}(G)=\lambda(G_1).$ Choose $v_1v_2\in E(G_2)$ and assume $v'$ is a vertex in $G_1$. {Let $G'$ be a graph obtained from $G-v_1v_2+v_1v'$ by deleting all of its isolated vertices.} Notice that $C_5$ is $2$-connected and $v_1v'$ is a cut edge of $G'$. Hence $G'$ is $C_5$-free and $G'\not\cong S_{\frac{m+3}{2},2}$. Clearly, $G_1$ is a proper subgraph of $G'$. Therefore, $\lambda(G')>\lambda(G_1)={\lambda}(G),$ which contradicts the choice of $G$. Hence, $G$ is connected.

Suppose that there is at least one cut vertex of $G$ in $V(G)\setminus\{u\}.$ Let $B$ be an end-block of $G$ with $u\not\in V(B)$ and let $v\in V(B)$ be a cut vertex of $G.$  Let $G''=G-\{wv:w\in V(B)\cap N(v)\}+\{wu:w\in V(B)\cap N(v)\}.$ Obviously, $u$ is a cut vertex of $G''$. Then $G''\not\cong S_{\frac{m+3}{2},2}.$ It is routine to check that $G''$ is $C_5$-free. That is, $G''\in \mathfrak{G}(m,C_5).$ In view of Lemma~\ref{lem1.1}, we have $\lambda(G'')>{\lambda}(G),$ a contradiction.

This completes the proof.
\end{proof}
Let $G$ be a connected graph and let ${\bf x}$ be its Perron vector. Then $A(G){{\bf x}}={\lambda}(G){{\bf x}}.$ Assume that $u$ is an arbitrary vertex of $G.$  Hence
\[\label{eq:1.1}
  {\lambda}(G)x_{u}=(A(G){\bf x})_{u}=\sum_{v\in N_0(u)}x_v+\sum_{v\in N(u)\setminus N_0(u)}x_v.
\]
In addition, ${\bf x}$ is also a characteristic eigenvector of $A^2(G)$ corresponding to ${{\lambda}^2(G)},$ i.e., ${A^2(G)}{\bf x}={{\lambda}^2(G)}{\bf x}.$ {For convenience, let $A^2(G)=(a_{uv}^{(2)}),$ where $a_{vu}^{(2)}$ denotes the number of walks of length 2 from $v$ to $u$ in $G$.}  It follows that
\[\label{eq:1.2}
  {{\lambda}^2(G)}x_{u}={\sum_{v\in V(G)}a_{uv}^{(2)}x_v}=d(u){x_u}+\sum_{v\in N(u)\setminus N_0(u)}d_{N(u)}(v)x_v+\sum_{w\in N^2(u)}d_{N(u)}(w)x_w.
\]
One may see that \eqref{eq:1.1} and \eqref{eq:1.2} will
be frequently used in the proofs of our main results.
\section{\normalsize Proof of Theorem \ref{thm1.04}}
In this section, we present the proof of Theorem \ref{thm1.04}, which determines the unique graph among $\mathcal{G}(m,\theta_{1,2,3})$ (resp. $\mathcal{G}(m,\theta_{1,2,4})$) having the largest spectral radius.

Let $\hat{G}$ be a graph in $\mathcal{G}(m,\theta_{1,2,3})$ (resp. $\mathcal{G}(m,\theta_{1,2,4})$) with the maximum spectral radius. In view of Lemma~\ref{lem5.1}(i), we know that $\hat{G}$ is connected. Assume that ${\bf \hat{x}}$ is the Perron vector of $\hat{G}$ and let $\hat{u}$ be the extremal vertex of $\hat{G}$. Notice that $\lambda(S_{\frac{m+3}{2},2})=\frac{1+\sqrt{4m-3}}{2}$ (see, \cite[Lemma 4.3]{ZLS}). On the other hand, since $S_{\frac{m+3}{2},2}$ is {$\theta_{1,2,3}$-free (resp. $\theta_{1,2,4}$-free),} one has
$$
  \hat{\lambda}:=\lambda(\hat{G})\geqslant \lambda(S_{\frac{m+3}{2},2})=\frac{1+\sqrt{4m-3}}{2}.
$$

Let $W=V(\hat{G})\setminus N[\hat{u}].$ Recall that $N_i(\hat{u})=\{u\in N(\hat{u}):u\ \text{lies in the component}\ {K_{i+1}}\ \text{of}\ \hat{G}[N(\hat{u})]\}.$ For $i\in\{0,1\},$ define $N_i^2(\hat{u})=\{w\in W:N(w)\cap N_i(\hat{u})\neq \emptyset\}.$ Denote by $N_+(\hat{u})=N(\hat{u})\setminus N_0(\hat{u})$ and $N_+^2(\hat{u})=\{w\in W:N(w)\cap N_+(\hat{u})\neq \emptyset\}.$

\subsection{\normalsize Proof of Theorem \ref{thm1.04}(i)}
{In this subsection, we prove Theorem \ref{thm1.04}(i). Now, $\hat{G}$ is a graph in $\mathcal{G}(m,\theta_{1,2,3})$ with the maximum spectral radius.} Notice that $\hat{G}$ is $\theta_{1,2,3}$-free, one obtains that $\hat{G}[N(\hat{u})]$ contains no path of length $3.$ Hence each component of $\hat{G}[N(\hat{u})]$ is either a triangle or a star $K_{1,r}$ for some $r\geqslant 0,$ where $K_{1,0}$ is a singleton component. Furthermore, $N_0^2(\hat{u})\cap N_1^2(\hat{u})=\emptyset,$ $|N(w)\cap N(\hat{u})|\leqslant 2$ for all $w\in N_1^2(\hat{u})$ and $|N(w)\cap N(\hat{u})|=1$ for all $w\in N_+^2(\hat{u})\setminus N_1^2(\hat{u}).$

Note that $\hat{\lambda}\geqslant \frac{1+\sqrt{4m-3}}{2}>3$ if $m\geqslant 8.$ Then ${\hat{\lambda}}^2-\hat{\lambda}\geqslant m-1.$ In view of \eqref{eq:1.1} and \eqref{eq:1.2}, one has
\begin{align}\label{eq:5.1}
  (m-1)x_{\hat{u}}&\leqslant ({\hat{\lambda}}^2-\hat{\lambda})x_{\hat{u}}= |N(\hat{u})|x_{\hat{u}}+\sum_{u\in N_+(\hat{u})}(d_{N(\hat{u})}(u)-1)x_u-\sum_{u\in N_0(\hat{u})}x_u+\sum_{w\in W}d_{N(\hat{u})}(w)x_w\\\notag
  &\leqslant \left(|N(\hat{u})|+2e(N_+(\hat{u}))-|N_+(\hat{u})|-\sum_{u\in N_0(\hat{u})}\frac{x_u}{x_{\hat{u}}}+e(W,N(\hat{u}))\right) x_{\hat{u}}.
\end{align}
It follows that
\[\label{eq:5.3}
  e(W)\leqslant e(N_+(\hat{u}))-|N_+(\hat{u})|-\sum_{u\in N_0(\hat{u})}\frac{x_u}{x_{\hat{u}}}+1
\]
and equality holds only if ${\hat{\lambda}}^2-\hat{\lambda}=m-1,$ and $x_u=x_{\hat{u}}$ for all $u\in \{w\in N_+(\hat{u}):d_{N(\hat{u})}(w)\geqslant 2\}\cup \{w\in W:d_{N(\hat{u})}(w)\geqslant 1\}.$ 
Recall that $e(N_+(\hat{u}))\leqslant |N_+(\hat{u})|.$ Hence, $e(W)\leqslant 1.$

In order to complete the proof of Theorem \ref{thm1.04}(i), we need the following three lemmas.
\begin{lem}\label{cl5.1}
  $e(W)=0.$
\end{lem}
\begin{proof}
Suppose that $e(W)=1,$ and let $w_1w_2$ be the unique edge in $\hat{G}[W].$ It follows from \eqref{eq:5.3} that $N_0(\hat{u})=\emptyset$ and $e(N_+(\hat{u}))-|N_+(\hat{u})|=0,$ which implies that each component of $\hat{G}[N(\hat{u})]$ is $K_3.$ {Note that $\hat{G}$ is $\theta_{1,2,3}$-free. Then $d_{N(\hat{u})}(w)=1$ for each $w\in N^2(u).$} Together with Lemma~\ref{lem5.1}(ii), one has $W=N^2(u)=\{w_1,w_2\}.$ Hence $x_{w_1}=x_{w_2}$ and $\hat{\lambda}x_{w_1}\leqslant x_{w_2}+x_{\hat{u}},$ i.e., $x_{w_1}\leqslant \frac{x_{\hat{u}}}{\hat{\lambda}-1}.$ By \eqref{eq:5.1}, we have
$$
  (m-1)x_{\hat{u}}\leqslant (|N(\hat{u})|+e(N(\hat{u}))) x_{\hat{u}}+x_{w_1}+x_{w_2}\leqslant \left(|N(\hat{u})|+e(N(\hat{u}))+\frac{2}{\hat{\lambda}-1}\right) x_{\hat{u}}.
$$
Therefore, $2\leqslant \frac{2}{\hat{\lambda}-1}$ and so $\hat{\lambda}\leqslant 2,$ a contradiction. {Hence $e(W)=0,$} as desired.
\end{proof}
\begin{lem}\label{cl5.2}
$W=N_0^2(\hat{u}).$
\end{lem}
\begin{proof}
Suppose to the contrary that $W\setminus N_0^2(\hat{u})\neq \emptyset.$ {Let $w$ be in $W\setminus N_0^2(\hat{u}).$} Recall that $N_0^2(\hat{u})\cap N_1^2(\hat{u})=\emptyset$ and $|N(w)\cap N(\hat{u})|=1$ for all $w\in N_+^2(\hat{u})\setminus N_1^2(\hat{u}).$ 
Together with Lemma~\ref{lem5.1}(ii) and Lemma~\ref{cl5.1}, we know that {$N(w)=\{u_1,u_2\},$ where $u_1u_2$ is a $K_2$ component of $\hat{G}[N(\hat{u})].$ Furthermore, $N(w')=\{u_1,u_2\}$ for each vertex $w'\in N_W(u_1)\cup N_W(u_2).$} Let $G=\hat{G}-\{w'u_1:w'\in N_{W}(u_1)\}+\{w'\hat{u}:w'\in N_{W}(u_1)\}.$ It is easy to see that $G\in \mathcal{G}(m,\theta_{1,2,3})$ and $\lambda(G)>\hat{\lambda}$ (based on  Lemma~\ref{lem1.1}), which contradicts the choice of $\hat{G}.$

This completes the proof.
\end{proof}
\begin{lem}\label{cl5.3}
$\hat{G}\cong S_{\frac{m+3}{2},2}$ or $\hat{G}\cong H_{t,s}\circ R_k$ for some bipartite graph $H_{t,s}$ and some integer $k\geqslant 1.$
\end{lem}
\begin{proof}
By \eqref{eq:5.3} and Lemma \ref{cl5.1}, one has $e(N_+(\hat{u}))\geqslant |N_+(\hat{u})|+\sum_{u\in N_0(\hat{u})}\frac{x_u}{x_{\hat{u}}}-1.$ If $N_0(\hat{u})\neq \emptyset,$ then $\sum_{u\in N_0(\hat{u})}\frac{x_u}{x_{\hat{u}}}>0$ and so $e(N_+(\hat{u}))\geqslant |N_+(\hat{u})|.$ It follows that $e(N_+(\hat{u}))=|N_+(\hat{u})|$ and $\hat{G}[N_+(\hat{u})]\cong kK_3.$ If $k=0,$ then $\hat{G}$ is a bipartite graph. Based on Lemma \ref{thm4}, one has $\hat{\lambda}\leqslant \sqrt{m}<\frac{1+\sqrt{4m-3}}{2}$ if $m\geqslant 8,$ a contradiction. {Together with Lemmas \ref{cl5.1} and \ref{cl5.2}, we obtain that} ${\hat{G}}\cong H_{t,s}\circ R_k$ for some integer $k\geqslant 1$ and some bipartite graph $H_{t,s},$ as desired.

Next, we consider the case that $N_0(\hat{u})= \emptyset.$ Thus, $N_0^2(\hat{u})=W=\emptyset$ (based on Lemma \ref{cl5.2}) and $e(N_+(\hat{u}))- |N_+(\hat{u})|\geqslant -1.$ Assume that $\hat{G}[N_+(\hat{u})]$ contains exactly $c$ non-trivial tree components. Notice that $c\leqslant 1.$ If $c=0,$ then $\hat{G}\cong R_k$ for $k=\frac{m}{6}.$ It is routine to check that $\hat{\lambda}<\frac{1+\sqrt{4m-3}}{2},$ a contradiction. Hence $c=1$ and $\hat{G}\cong K_1\vee (K_{1,r}\cup kK_3)$ for some $k\geqslant 0$ and $r\geqslant 1.$ If $k\geqslant 1,$ then we assume that $H$ is a $K_4$ block of $\hat{G}$ with $V(H)=\{u_1,u_2,u_3,\hat{u}\}.$ Obviously, $x_{u_1}=x_{u_2}=x_{u_3}$ and $\hat{\lambda}x_{u_1}=2x_{u_1}+x_{\hat{u}}.$ Therefore, $x_{u_1}=\frac{x_{\hat{u}}}{\hat{\lambda}-2}<x_{\hat{u}}$ (since $\hat{\lambda}>3).$ Thus, the inequality in \eqref{eq:5.3} is strict and so $e(W)<0,$ a contradiction. Hence $k=0$ and $\hat{G}\cong S_{\frac{m+3}{2},2}.$

This completes the proof.
\end{proof}

\begin{proof}[\bf Proof of Theorem \ref{thm1.04}(i)]
{In view of Lemma \ref{cl5.3}, in order to prove Theorem \ref{thm1.04}(i), it suffices to show that $\hat{G}\not\cong H_{t,s}\circ R_k$ for any bipartite graph $H_{t,s}$ and any integer $k\geqslant 1.$}

Suppose that $\hat{G}\cong H_{t,s}\circ R_k,$ where $H_{t,s}$ is a bipartite graph and $k$ is a positive integer. Let $H$ be a $K_4$ block of $\hat{G}$ with  $V(H)=\{u_1,u_2,u_3,u\}$ and $d_{\hat{G}}(u)\geqslant 4.$ Now, we give the following claim to determine the unique extremal vertex of $\hat{G}.$
\begin{claim}\label{cla3.1}
  $x_v<x_u$ for all $v\neq u,$ i.e., $u=\hat{u}.$
\end{claim}
\begin{proof}[\bf Proof of Claim \ref{cla3.1}]
Note that $\hat{\lambda}>{3}.$ Hence, $x_{u_1}=x_{u_2}=x_{u_3}=\frac{x_{u}}{\hat{\lambda}-2}<x_{u}$. So, if there exists a vertex $v\,(\neq u)$ such that $x_v\geqslant x_u,$ then $v\in V(H_{t,s}).$ Let $G_1=\hat{G}-\{uu_i:1\leqslant i\leqslant 3\}+\{vu_i:1\leqslant i\leqslant 3\}.$ It is easy to see that $G_1\in \mathcal{G}(m,\theta_{1,2,3}).$ In addition, it follows from Lemma~\ref{lem1.1} that $\lambda(G_1)>\hat{\lambda},$ which contradicts the choice of $\hat{G}.$ Hence $x_v<x_u$ for all $v\neq u,$ i.e., $u=\hat{u}.$
\end{proof}

For convenience, assume that $T\cup S$ is a bipartition of $H_{t,s}$ with $T=\{v_1,\ldots,v_t\}$ satisfying $x_{v_1}\geqslant \cdots\geqslant x_{v_t}$ and $S=\{w_1,\ldots,w_s\}$ satisfying $x_{w_1}\geqslant \cdots\geqslant x_{w_s}.$ By Lemma~\ref{lem1.1}, one obtains $v_iw_j\in E(\hat{G})$ only if $v_{p}w_q\in E(\hat{G})$ for each $p\leqslant i$ and $q\leqslant j$.

By Lemma \ref{lem5.2}, we know that $s\geqslant 1.$ If $s=1,$ then $v_1w_1\in E(\hat{G}).$ Define $G_2=\hat{G}-v_1w_1+\hat{u}w_1.$ Then $G_2\in \mathcal{G}(m,\theta_{1,2,3})$ and $\lambda(G_2)>\hat{\lambda}$ (by Lemma~\ref{lem1.1}), a contradiction. Hence $s\geqslant 2.$ In view of Lemma~\ref{lem5.1}(ii), one has $d(w_i)\geqslant 2$ for each $w_i\in S$ and so $e(T,S)\geqslant 2s.$ Furthermore, $N_S(v_1)=N_S(v_2)=S$ and $x_{v_1}=x_{v_2}.$

If $e(T,S)=2s,$ then $N(w_1)=\{v_1,v_2\}.$ Therefore, $\hat{\lambda}x_{w_1}=2x_{v_1},$ that is, $x_{w_1}=\frac{2x_{v_1}}{\hat{\lambda}}<x_{v_1}.$ If $e(T,S)=2s+1,$ then $N(w_1)=\{v_1,v_2,v_3\}.$ Hence $\hat{\lambda}x_{w_1}=2x_{v_1}+x_{v_3},$ i.e., $x_{w_1}{\leqslant}\frac{3x_{v_1}}{\hat{\lambda}}<x_{v_1}.$ In both cases, $x_{w_1}<x_{v_1}.$ Define $G_3=\hat{G}-w_1v_2+v_1v_2$ if $e(T,S)=2s,$ and $G_3=\hat{G}-w_1v_3+v_1v_2$ if $e(T,S)=2s+1.$ Clearly, $G_3\in \mathcal{G}(m,\theta_{1,2,3})$ and $\lambda(G_3)>\hat{\lambda}$ (based on Lemma~\ref{lem1.1}), a contradiction. Therefore, $e(T,S)\geqslant 2s+2\geqslant 6$ and $t\geqslant 3.$ Hence $m=6k+t+e(T,S)\geqslant 6k+t+2s+2\geqslant 6k+2s+5\geqslant 15.$

In order to obtain a contradiction, we show that $\hat{\lambda}< \frac{1+\sqrt{4m-3}}{2}$ by induction on $m$. If $k=1,$ then by Lemma~\ref{lem4.0}, we have
$$
  {\hat{\lambda}}^2\leqslant \lambda^2(\hat{G}-u_1)+5\leqslant \lambda^2(\hat{G}-u_1-u_2)+8\leqslant \lambda^2(\hat{G}-u_1-u_2-u_3)+9.
$$
Note that $\hat{G}-u_1-u_2-u_3$ is bipartite. Together with Lemma \ref{thm4}, one has $\lambda^2(\hat{G}-u_1-u_2-u_3)\leqslant m-6.$ Hence ${\hat{\lambda}}\leqslant \sqrt{m+3}<\frac{1+\sqrt{4m-3}}{2}$ if $m\geqslant 15.$ Therefore, $\hat{\lambda}< \frac{1+\sqrt{4m-3}}{2}$ holds for $m=15.$

Assume that the result holds for $|E(\hat{G})|\leqslant m-1.$ Now, we consider that $\hat{G}$ contains $m$ edges and $k\,(\geqslant 2)$ $K_4$ blocks. Hence $m\geqslant 21$ and $\hat{\lambda}\geqslant 5.$ Note that $x_{u_1}=x_{u_2}=x_{u_3}$ and $\hat{\lambda}x_{u_1}=2x_{u_1}+x_{\hat{u}}.$ By Lemma \ref{lem:4.0}, we know that  $x_{\hat{u}}\leqslant \frac{1}{\sqrt{2}}$ and so $x^2_{u_1}\leqslant \frac{1}{2(\hat{\lambda}-2)^2}.$ Let $G_4=\hat{G}-\{u_iu_j:1\leqslant i<j\leqslant 3\}.$ Then by the Rayleigh quotient, one has
$$
  \lambda(G_4)\geqslant {\hat{{\bf x}}}^TA(G_4)\hat{{\bf x}}={\hat{{\bf x}}}^TA(\hat{G})\hat{{\bf x}}-6x^2_{u_1}=\hat{\lambda}-6x^2_{u_1}.
$$
By induction, one has $\lambda(G_4)< \frac{1+\sqrt{4m-15}}{2}$ and so $\lambda^2(G_4)-\lambda(G_4)< m-4.$ It follows that
\allowdisplaybreaks
\begin{align*}
  m-4&> \lambda^2(G_4)-\lambda(G_4)\geqslant (\hat{\lambda}-6x^2_{u_1})(\hat{\lambda}-6x^2_{u_1}-1)={\hat{\lambda}}^2-(12x^2_{u_1}+1)\hat{\lambda}+6x^2_{u_1}(6x^2_{u_1}+1)\\
  &>{\hat{\lambda}}^2-\hat{\lambda}-6x^2_{u_1}(2\hat{\lambda}-1)\geqslant {\hat{\lambda}}^2-\hat{\lambda}-\frac{3(2\hat{\lambda}-1)}{(\hat{\lambda}-2)^2}\geqslant {\hat{\lambda}}^2-\hat{\lambda}-3.
\end{align*}
The last inequality follows by $\frac{2\hat{\lambda}-1}{(\hat{\lambda}-2)^2}\leqslant 1$ if $\hat{\lambda}\geqslant 5.$ Hence ${\hat{\lambda}}<\frac{1+\sqrt{4m-3}}{2},$ as desired.

This completes the proof of Theorem \ref{thm1.04}(i).
\end{proof}

\subsection{\normalsize Proof of Theorem \ref{thm1.04}(ii)}
{In this subsection, we present the proof of Theorem \ref{thm1.04}(ii), and assume $\hat{G}$ is a graph in $\mathcal{G}(m,\theta_{1,2,4})$ with the maximum spectral radius.} Note that $\hat{G}$ is $\theta_{1,2,4}$-free, $\hat{G}[N(\hat{u})]$ contains no path of length $4.$ For each component $H$ of $\hat{G}[N(\hat{u})],$ denote by $W_H=\{w\in W:N_H(w)\neq \emptyset\}.$  Thus, $W_{H_i}\cap W_{H_j}=\emptyset$ for any two distinct component $H_i$ and $H_j$ of $\hat{G}[N(\hat{u})],$ unless one of $H_i$ and $H_j$ is an isolated vertex and the other is a star (in particular, vertices in $W_{H_i}\cap W_{H_j}$ must be adjacent to the center vertex of the star). Note that $\hat{\lambda}\geqslant \lambda(S_{\frac{m+3}{2},2})=\frac{1+\sqrt{4m-3}}{2}>5$ if $m\geqslant 22.$ Therefore, ${\hat{\lambda}}^2-\hat{\lambda}\geqslant m-1.$ Together with \eqref{eq:1.1} and \eqref{eq:1.2}, one has
$$
  (m-1)x_{\hat{u}}\leqslant ({\hat{\lambda}}^2-\hat{\lambda})x_{\hat{u}}\leqslant |N(\hat{u})|x_{\hat{u}}+\sum_{u\in N_+(\hat{u})}(d_{N(\hat{u})}(u)-1)x_u-\sum_{u\in N_0(\hat{u})}{{x_u}}+e(W,N(\hat{u}))x_{\hat{u}}.
$$
It follows that
$$
  \left(m-1-|N(\hat{u})|-e(W,N(\hat{u}))+\sum_{u\in N_0(\hat{u})}\frac{x_u}{x_{\hat{u}}}\right)x_{\hat{u}}\leqslant \sum_{u\in N_+(\hat{u})}(d_{N(\hat{u})}(u)-1)x_u.
$$
Denote $\eta(H):=\sum_{u\in V(H)}(d_{H}(u)-1)x_u$ for each non-trivial connected component of $\hat{G}[N(\hat{u})].$ This gives us
\[\label{eq:5.5}
  \left(e(N(\hat{u}))+e(W)-1+\sum_{u\in N_0(\hat{u})}\frac{x_u}{x_{\hat{u}}}\right)x_{\hat{u}}\leqslant \sum_H \eta(H),
\]
where $H$ takes over all non-trivial connected components of $\hat{G}[N(\hat{u})].$

Further on, we need the the following two lemmas.
\begin{lem}\label{cl5.4}
$\hat{G}[N(\hat{u})]$ is $\{C_4,C_5,\ldots\}$-free.
\end{lem}
\begin{proof}
Notice that $\hat{G}$ is $\theta_{1,2,4}$-free. Then $\hat{G}[N(\hat{u})]$ does not contain any path of length $4$ and any cycle of length more than $4$. In particular, if $H$ is a component of $\hat{G}[N(\hat{u})]$ that contains $C_4$, then $H\in\{C_4,C_4+e,K_4\}$; if $H$ is a component of $\hat{G}[N(\hat{u})]$ with circumference $3$ (i.e., the length of the longest cycle in graph $H$ is $3$), then $H\cong K_{1,r}+e$ for some $r\geqslant 2.$ Let $\mathcal{H}$ be the family of components of $\hat{G}[N(\hat{u})]$ each of which contains $C_4$ as a subgraph and let $\mathcal{H}'$ be the family of other non-trivial components of $\hat{G}[N(\hat{u})]$. Then $\mathcal{H}'$ consists of trees and/or unicyclic graphs. Hence, for each $H\in \mathcal{H}',$ we have
$$
  \eta(H)=\sum_{u\in V(H)}(d_{H}(u)-1)x_u\leqslant (2e(H)-|V(H)|)x_{\hat{u}}\leqslant e(H)x_{\hat{u}}.
$$

Since $\hat{G}$ is $\theta_{1,2,4}$-free, one has $W_{H_i}\cap W_{H_j}=\emptyset$ and $e(W_{H_i},W_{H_j})=0$ if ${H_i},{H_j}\in \mathcal{H}$ and $H_i\neq H_j.$ Hence
$$
  e(W)\geqslant \sum_{H\in \mathcal{H}}e(W_H,W).
$$
It follows from \eqref{eq:5.5} that
$$
  \left(\sum_{H\in \mathcal{H}}(e(H)+e(W_H,W))-1\right)x_{\hat{u}}\leqslant \sum_{H\in \mathcal{H}}\eta(H).
$$

Suppose to the contrary that $\mathcal{H}\neq \emptyset.$ In order to get a contradiction, it suffices to show $\eta(H)< (e(H)+e(W_H,W)-1)x_{\hat{u}}$ holds for each $H\in \mathcal{H}.$ Let $H^*\in \mathcal{H}$ with $V(H^*)=\{u_1,u_2,u_3,u_4\}.$ If $W_{H^*}=\emptyset,$ then assume that $x_{u_1}=\max\{x_{u_i}:1\leqslant i\leqslant 4\}.$ Hence $\hat{\lambda}x_{u_1}=\sum_{u\in N(u_1)}x_u\leqslant x_{\hat{u}}+3x_{u_1},$ i.e., $x_{u_1}\leqslant \frac{x_{\hat{u}}}{\hat{\lambda}-3}<\frac{x_{\hat{u}}}{2}$ since $\hat{\lambda}>5.$ It follows that
$$
  \eta(H^*)=\sum_{u\in V(H^*)}(d_{H^*}(u)-1)x_u\leqslant (2e(H^*)-4)x_{u_1}\leqslant (e(H^*)-2)x_{\hat{u}}<(e(H^*)-1)x_{\hat{u}},
$$
as desired. So, in what follows, we assume that $W_{H^*}\neq \emptyset.$ Notice that $d_{N(\hat{u})}(w)=d_{H^*}(w)=1$ if $w\in W_{H^*}.$ Together with Lemma~\ref{lem5.1}(ii), one has $e(W_{H^*},W)\geqslant 1.$ We proceed by considering the following two cases.

{\bf Case 1.} All vertices in $W_{H^*}$ have a unique common neighbor, say $u_1,$ in $H^*.$ That is, $N_W(u_i)=\emptyset$ for $i\in\{2,3,4\}.$ Assume that $x_{u_2}=\max\{x_{u_i}:2\leqslant i\leqslant 4\}.$ Therefore, $\hat{\lambda}x_{u_2}{\leqslant }x_{\hat{u}}+x_{u_1}+x_{u_3}+x_{u_4}\leqslant 2(x_{\hat{u}}+x_{u_2}).$ It follows that $x_{u_2}\leqslant \frac{2x_{\hat{u}}}{\hat{\lambda}-2}<\frac{2x_{\hat{u}}}{3}.$ Hence, 
\allowdisplaybreaks
\begin{align*}
  \eta(H^*)&=\sum_{u\in V(H^*)}(d_{H^*}(u)-1)x_u\leqslant (d_{H^*}(u_1)-1)x_{u_1}+(2e(H^*)-d_{H^*}(u_1)-3)x_{u_2}\\
  &< (d_{H^*}(u_1)-1)x_{\hat{u}}+\left(\frac{4}{3}e(H^*)-\frac{2}{3}d_{H^*}(u_1)-2\right)x_{\hat{u}}\\
  &\leqslant \left(\frac{4}{3}e(H^*)-2\right)x_{\hat{u}}\leqslant e(H^*)x_{\hat{u}}\leqslant (e(H^*)+e(W_{H^*},W)-1)x_{\hat{u}},
\end{align*}
as desired.

{\bf Case 2.} There exist $k\ (\geqslant 2)$ vertices, say $w_1,\ldots,w_k,$ of $W_{H^*}$ such that they have mutual distinct neighbors in $V(H^*).$ Note that $N_{H^*}(w)=N_{H^*}(w')$ if $ww'\in E(\hat{G}[W_{H^*}]).$ Hence  $\{w_1,\ldots,w_k\}$ is an independent set of $\hat{G}.$ It follows from Lemma~\ref{lem5.1}(ii)-(iii) that $d(w_i)\geqslant 2$ for $1\leqslant i\leqslant k$ and $d(w_i)=2$ holds for at most one $w_i.$ Therefore, $\sum_{1\leqslant i\leqslant k}d(w_i)\geqslant 3k-1$ and $e(W_{H^*},W)\geqslant 2k-1.$ Thus,
$$
  e(H^*)\leqslant e(K_4)=6\leqslant e(W_{H^*},W)-2k+7\leqslant e(W_{H^*},W)+3.
$$
So,
$$
  \eta(H^*)=\sum_{u\in V(H^*)}(d_{H^*}(u)-1)x_u\leqslant (2e(H^*)-4)x_{\hat{u}}\leqslant (e(H^*)+e(W_{H^*},W)-1)x_{\hat{u}},
$$
and $\eta(H^*)=(e(H^*)+e(W_{H^*},W)-1)x_{\hat{u}}$ holds only if $k=2,$ $H^*\cong K_4$ and $x_{u_i}=x_{\hat{u}}$ for all $i\in\{1,2,3,4\}.$ On the other hand, if $k=2,$ then there exists a vertex $u_j\in V(H^*)$ satisfying $N_W(u_j)=\emptyset.$ Hence $\hat{\lambda}x_{u_j}\leqslant 4x_{\hat{u}}.$ Thus, $x_{u_j}\leqslant \frac{4x_{\hat{u}}}{\hat{\lambda}}<x_{\hat{u}}$ since $\hat{\lambda}>5.$ Therefore, $\eta(H^*)<(e(H^*)+e(W_{H^*},W)-1)x_{\hat{u}},$ as desired.

This completes the proof.
\end{proof}
In view of Lemma \ref{cl5.4}, one obtains that each component of $\hat{G}[N(\hat{u})]$ is either a tree or a unicyclic graph $K_{1,r}+e$ for some $r\geqslant 2.$ Assume that there are $c$ non-trivial tree components in $\hat{G}[N(\hat{u})].$ Then
$$
  \sum_{H}\eta(H)=\sum_{H}\sum_{u\in V(H)}(d_{H}(u)-1)x_u\leqslant \sum_{H}(2e(H)-|V(H)|)x_{\hat{u}}=(e(N(\hat{u}))-c)x_{\hat{u}},
$$
where $H$ takes over all non-trivial components of $\hat{G}[N(\hat{u})].$ Combining with \eqref{eq:5.5} gives us
\[\label{eq:4.7}
  e(W)\leqslant 1-c-\sum_{u\in N_0(\hat{u})}\frac{x_u}{x_{\hat{u}}}.
\]
Hence $e(W)\leqslant 1$ and $c\leqslant 1.$ In addition, $e(W)=1$ holds only if $c=0$ and $N_0(\hat{u})=\emptyset.$
\begin{lem}\label{cl5.5}
$e(W)=0$ and each component of $\hat{G}[N(\hat{u})]$ is a tree.
\end{lem}
\begin{proof}
Suppose that $e(W)=1$ and let $ww'$ be the unique edge in $\hat{G}[W].$ Then $c=0$ and $N_0(\hat{u})=\emptyset.$ That is, each component of $\hat{G}[N(\hat{u})]$ is isomorphic to a unicyclic graph. Since $G$ is $\theta_{1,2,4}$-free, we know that $w$ and $w'$ have the unique common neighbor $u$ in ${N(\hat{u})}.$ Therefore, $u$ is a cut vertex of $\hat{G},$ a contradiction to Lemma~\ref{lem5.1}(ii). {Hence $e(W)=0.$}

In what follows, we are to prove that each component of $\hat{G}[N(\hat{u})]$ is a tree, i.e., $\hat{G}[N(\hat{u})]$ is $C_3$-free. Suppose that $H^*\cong K_{1,r}+e$ is a component of $\hat{G}[N(\hat{u})]$ for some $r\geqslant 2.$ We consider the following claims.
\begin{claim}\label{cl:5.1}
$H^*\not\cong K_3.$
\end{claim}
\begin{proof}[\bf Proof of Claim \ref{cl:5.1}]
Suppose that $H^*\cong K_3$ with $V(H^*)=\{u_1,u_2,u_3\}.$ If $W_{H^*}= \emptyset,$ then $x_{u_1}=x_{u_2}=x_{u_3}=\frac{x_{\hat{u}}}{\hat{\lambda}-2}<\frac{x_{\hat{u}}}{3}.$ Hence,
$$
  \eta(H^*)=\sum_{1\leqslant i\leqslant 3}(d_{H^*}(u_i)-1)x_{u_i}=3x_{u_1}<x_{\hat{u}}=(e(H^*)-2)x_{\hat{u}}.
$$
Recall that $e(W)=0$ and $\eta(H)\leqslant e(H)x_{\hat{u}}$ for each non-trivial component $H$ of $\hat{G}[N(\hat{u})].$ Together with \eqref{eq:5.5}, one has
$$
  \left(e(N(\hat{u}))-1+\sum_{u\in N_0(\hat{u})}\frac{x_u}{x_{\hat{u}}}\right)x_{\hat{u}}\leqslant \sum_H \eta(H)<(e(N(\hat{u}))-2)x_{\hat{u}},
$$
a contradiction. Hence $W_{H^*}\neq \emptyset.$

It is routine to check that $2\leqslant d(w)\leqslant |V(H^*)|=3$ for each $w\in W_{H^*}.$ If there exists a vertex $w\in W_{H^*}$ such that $d(w)=3,$ then $W_{H^*}=\{w\}.$ If $d(w)=2$ for all $w\in W_{H^*},$ then Lemma~\ref{lem5.1}(iii) implies that all vertices in $W_{H^*}$ share the same neighborhood. Without loss of generality, assume that $N(w)=\{u_1,u_2\}$ for each $w\in W_{H^*}.$ Put $G_1:=\hat{G}-\{wu_1:w\in N_{W_{H^*}}(u_1)\}+\{w\hat{u}:w\in N_{W_{H^*}}(u_1)\}.$ Then in both cases, one has $G_1\in \mathcal{G}(m,\theta_{1,2,4})$ and $\lambda(G_1)>\hat{\lambda}$ (based on Lemma~\ref{lem1.1}), a contradiction.

This completes the proof of Claim \ref{cl:5.1}.
\end{proof}
\begin{claim}\label{cl:5.2}
$W_{H^*}=\emptyset$ and $H^*$ is the unique non-trivial component of $\hat{G}[N(\hat{u})].$
\end{claim}
\begin{proof}[\bf Proof of Claim \ref{cl:5.2}]
In view of Claim \ref{cl:5.1}, we derive that $H^*\cong K_{1,r}+e$ for some $r\geqslant 3.$ Since $\hat{G}$ is $\theta_{1,2,4}$-free, one obtains that $d_{N(\hat{u})}(w)=1$ for each $w\in W_{H^*}.$ On the other hand, Lemma~\ref{lem5.1}(ii) implies that $d(w)\geqslant 2$ for each $w\in W.$ Together with $e(W)=0$, we get  $W_{H^*}=\emptyset.$

Now, we show that $H^*$ is the unique non-trivial component of $\hat{G}[N(\hat{u})].$ Suppose that $\hat{G}[N(\hat{u})]$ contains another non-trivial component $H.$ Note that there are $r-2$ vertices in $V(H^*)$ with degree two in $\hat{G}.$ By Lemma~\ref{lem5.1}(iii), there exist no vertices of degree two out of $H^*.$ That is, $H^*$ is the unique component which contains triangles and so $H$ is a tree. Based on \eqref{eq:4.7}, one has $c=1$ and $N_0(\hat{u})= \emptyset.$ In addition, $W_H\neq \emptyset$ and $d(w)\geqslant 3$ for each $w\in W_H\cup V(H).$ Note that $e(W)=0$ and $W_H\cap W_{H^*}=\emptyset.$ Then $N(w)\subseteq V(H)$ for all $w\in W_H$ and so $2\leqslant \diam(H)\leqslant 3.$

If $\diam(H)=2,$ then $H$ is a star, say $K_{1,r'}$. Assume that $V(H)=\{u_1,\ldots,u_{r'},u\}$ with $d_H(u)=r'\geqslant 2.$ Since no vertices of degree two out of $H^*,$ there exists a vertex $w_1\in W$ such that $u_1w_1\in E(\hat{G})$ and $d(w_1)\geqslant 3.$ Hence $N(w_1)\setminus \{u_1,u\}\neq \emptyset.$ Recall that $N(w_1)\subseteq V(H).$ Without loss of generality, assume that $u_2\in N(w_1).$ If $r'=2,$ i.e., $H\cong K_{1,2},$ then $N(w_1)=V(H).$ For the reason of $\theta_{1,2,4}$-free, we obtain that $W_H=\{w_1\}.$ Let $G_2=\hat{G}-w_1u+w_1\hat{u}.$ It is routine to check that $G_2\in \mathcal{G}(m,\theta_{1,2,4})$ and $\lambda(G_2)>\hat{\lambda}$ (based on Lemma~\ref{lem1.1}), a contradiction. Hence, $r'\geqslant 3.$ Clearly, $\hat{G}[\{\hat{u},u_3,u,u_2,w_1,u_1\}]$ contains $\theta_{1,2,4}$ as a subgraph, a contradiction.

If $\diam(H)=3,$ then $H$ is a double star $D_{a,b}$. Assume that $V(H)=\{u_1,\ldots,u_a,u,v,v_1,\ldots,v_b\},$ where $N_H(u)=\{u_1,\ldots,u_a,v\}$ and $N_H(v)=\{v_1,\ldots,v_b,u\}.$ Since no vertices of degree two out of $H^*,$ there exists a vertex $w_1\in N_W(u_1)$ with $d(w_1)\geqslant 3.$ Hence $N(w_1)\setminus \{u_1,u\}\neq \emptyset.$ Notice that $N(w_1)\subseteq V(H).$ If $v\in N(w_1),$ then $\hat{G}[\{\hat{u},v_1,v,w_1,u_1,u\}]$ contains $\theta_{1,2,4}$ as a subgraph, a contradiction. If $u_i\in N(w_1)$ for some $i\in \{2,\ldots,a\},$ then $\theta_{1,2,4}$ is a subgraph of $\hat{G}[\{\hat{u},v,u,u_i,w_1,u_1\}],$ a contradiction. If $v_i\in N(w_1)$ for some $i\in \{1,\ldots,b\},$ then $\hat{G}[\{\hat{u},u,u_1,w_1,v_i,v\}]$ contains $\theta_{1,2,4}$ as a subgraph, a contradiction. 


This completes the proof of Claim \ref{cl:5.2}.
\end{proof}
In view of Claims \ref{cl:5.1} and \ref{cl:5.2}, we know that $\hat{G}\cong H_{|N_0(\hat{u})|,|W|}\circ (K_1\vee H^*),$ where $H^*\cong K_{1,r}+e$ for some $r\geqslant 3.$ Assume that $V(H^*)=\{u_0,u_1,\ldots,u_r\}$ with $d_{H^*}(u_0)=r$ and $d_{H^*}(u_1)=d_{H^*}(u_2)=2.$ Hence $x_{u_1}=x_{u_2}$ and $x_{u_3}=\cdots=x_{u_r}.$ Note that $\hat{\lambda}x_{u_1}=x_{u_0}+x_{u_2}+x_{\hat{u}}\leqslant x_{u_1}+2x_{\hat{u}}.$ It follows that $x_{u_1}\leqslant \frac{2x_{\hat{u}}}{\hat{\lambda}-1}<\frac{x_{\hat{u}}}{2}$ since $\hat{\lambda}>5.$ Therefore,
$$
  \eta(H^*)=(r-1)x_{u_0}+x_{u_1}+x_{u_2}\leqslant (r-1)x_{\hat{u}}+2x_{u_1}<rx_{\hat{u}}=(e(H^*)-1)x_{\hat{u}}.
$$
Together with $e(W)=0,$ one may see that \eqref{eq:5.5} becomes
$$
  \left(e(H^*)-1+\sum_{u\in N_0(\hat{u})}\frac{x_u}{x_{\hat{u}}}\right)x_{\hat{u}}<(e(H^*)-1)x_{\hat{u}},
$$
a contradiction. This completes the proof.
\end{proof}

Now we come back to prove Theorem \ref{thm1.04}(ii).
\begin{proof}[\bf Proof of Theorem \ref{thm1.04}(ii)]
By Lemmas \ref{cl5.4} and \ref{cl5.5}, one obtains that $\hat{G}[N(\hat{u})]$ consists of $c$ non-trivial tree components and some isolated vertices. Based on \eqref{eq:4.7}, we have $c\leqslant 1$ and equality holds only if $N_0(\hat{u})=\emptyset.$ If $c=0,$ then $\hat{G}$ is bipartite. Together with Lemma \ref{thm4}, one has $\hat{\lambda}\leqslant \sqrt{m}<\frac{1+\sqrt{4m-3}}{2}$ if $m\geqslant 22,$ a contradiction. Hence $c=1$ and $N_0(\hat{u})=\emptyset.$ That is, $\hat{G}[N(\hat{u})]\cong H,$ where $H$ is a non-trivial tree. Recall that $\diam(H)\leqslant 3.$

If $\diam(H)\leqslant 2,$ then $H$ is a star $K_{1,r}$ for some  $r\geqslant 1.$ Assume that $V(H)=\{u_0,u_1,\ldots,u_r\}$ with $d_H(u_0)=r$ and $r\geqslant 1.$ Note that $\hat{\lambda}(x_{u_0}-x_{\hat{u}})=x_{\hat{u}}+\sum_{w\in N_W(u_0)}x_{w}-x_{u_0},$ i.e., $(\hat{\lambda}-1)(x_{u_0}-x_{\hat{u}})=\sum_{w\in N_W(u_0)}x_{w}\geqslant 0$ and equality holds if and only if $N_W(u_0)=\emptyset.$ By the choice of $\hat{u},$ we obtain that $N_W(u_0)=\emptyset.$ {If $W_H\neq \emptyset,$ then let $w_0$ be in $W_H.$ By Lemmas \ref{lem5.1}(ii) and \ref{cl5.5}, one has $d_{H}(w_0)\geqslant 2.$ Assume, without loss of generality, that $\{u_1,u_2\}\in N_H(w_0).$ If $r\geqslant 3,$ then $\hat{G}[\{\hat{u},u_1,w_0,u_2,u_0,u_3\}]$ contains $\theta_{1,2,4}$ as a subgraph, a contradiction. Hence $r=2.$ Recall that $m\geqslant 22.$ Then $|W_H|\geqslant 9.$ It is routine to check that $x_{u_1}>x_{\hat{u}},$ a contradiction.} It follows that $W_H=\emptyset$ and $\hat{G}\cong S_{\frac{m+3}{2},2}.$

If $\diam(H)=3,$ then $H$ is a double star. For the reason of $\theta_{1,2,4}$-free, one has $d_{N(\hat{u})}(w)=1$ for each $w\in W_H.$ Together with Lemmas~\ref{lem5.1}(ii) and \ref{cl5.5}, we get {$W_H=\emptyset,$ i.e., $\hat{G}\cong K_1\vee H.$ Then $\hat{G}$ contains two non-adjacent degree two vertices with distinct neighborhoods, a contradiction to Lemma~\ref{lem5.1}(iii)}.

This completes the proof of Theorem \ref{thm1.04}(ii).
\end{proof}

\section{\normalsize Proofs of Theorems \ref{thm1.3} and \ref{thm1.03}}
{In this section, we give the proofs of Theorems \ref{thm1.3} and \ref{thm1.03}, which characterize graphs among $\mathfrak{G}(m,C_5)$ and $\mathfrak{G}(m,C_6)$ having the maximum spectral radius for $m\geqslant 22$, respectively.} 

Recall that $S_{n,2}^t$ is a graph obtained from $S_{n-t,2}$ by attaching $t$ pendant vertices to a maximum degree vertex of $S_{n-t,2}.$ Before giving the proof, we present a basic lemma, which investigates the spectral radius of $S_{n,2}^t.$
\begin{lem}\label{lem4.00}
Let $m\,(\geqslant 22)$ and $t$ be two positive integers {such that $m\geqslant t+1$.} Denote by $s_t(x):=x^4-mx^2-(m-t-1)x-\frac{t^2-mt+t}{2}$ a real function in $x.$ 
\begin{wst}
\item[{\rm (i)}] The spectral radius of $S_{\frac{m+t+3}{2},2}^t$ is the largest root of $s_t(x)=0;$
\item[{\rm (ii)}] $\min\left\{\lambda(S_{\frac{m+4}{2},2}^1),\lambda(S_{\frac{m+5}{2},2}^2)\right\}>\frac{1+\sqrt{4m-7}}{2};$
\item[{\rm (iii)}] $\lambda(S_{\frac{m+t+3}{2},2}^t)<\lambda(S_{\frac{m+4}{2},2}^1)$ if $m$ is even and $t\,(\geqslant 3)$ is odd, and $\lambda(S_{\frac{m+t+3}{2},2}^t)<\lambda(S_{\frac{m+5}{2},2}^2)$ if $m$ is odd and $t\,(\geqslant 4)$ is even.
\end{wst}
\end{lem}
\begin{proof}
(i) Assume that $V(S_{\frac{m+t+3}{2},2}^t)=\{v,u,u_1,\ldots,u_{\frac{m-t-1}{2}},v_1,\ldots,v_t\},$ where $v$ is the maximum degree vertex, $d(u_i)=2$ for $1\leqslant i\leqslant \frac{m-t-1}{2}$ and $d(v_i)=1$ for $1\leqslant i\leqslant t.$ It is routine to check that $\pi:=\{v\}\cup \{u\}\cup \{u_1,\ldots,u_{\frac{m-t-1}{2}}\}\cup \{v_1,\ldots,v_t\}$ is an equitable partition of $A(S_{\frac{m+t+3}{2},2}^t).$ Hence the quotient matrix of $A(S_{\frac{m+t+3}{2},2}^t)$ corresponding to the partition $\pi$ can be written as
\begin{equation*}
    A(S_{\frac{m+t+3}{2},2}^t)_\pi=\left(
                       \begin{array}{cccc}
                         0 & 1 & \frac{m-t-1}{2} & t \\
                         1 & 0 & \frac{m-t-1}{2} & 0 \\
                         1 & 1 & 0 & 0 \\
                         1 & 0 & 0 & 0 \\
                       \end{array}
                     \right).
\end{equation*}
By a direct calculation, the characteristic polynomial of $A(S_{\frac{m+t+3}{2},2}^t)_\pi$ is $s_t(x)=x^4-mx^2-(m-t-1)x-\frac{t^2-mt+t}{2}.$ Together with Lemma~\ref{lem3.02}, (i) holds.

(ii) By Mathematica 9.0 \cite{W2012}, one obtains that for $m\geqslant 22,$
$$
  s_1\left(\frac{1+\sqrt{4m-7}}{2}\right)<0\ \ \text{and}\ \ s_2\left(\frac{1+\sqrt{4m-7}}{2}\right)<0.
$$
That is to say, $\min\left\{\lambda(S_{\frac{m+4}{2},2}^1),\lambda(S_{\frac{m+5}{2},2}^2)\right\}>\frac{1+\sqrt{4m-7}}{2},$ as desired.

(iii) Note that $m\geqslant t+1.$ If $m=t+1,$ then $S_{\frac{m+t+3}{2},2}^t\cong K_{1,m}$ and $\lambda(K_{1,m})=\sqrt{m}<\frac{1+\sqrt{4m-7}}{2}$ for $m\geqslant 22.$ Together with (ii), our result follows immediately.

Note that {in $S_{\frac{m+t+3}{2},2}^t$}, $m$ does not have the same parity as that of $t$. So in what follows, we consider $m\geqslant t+3.$ If $m$ is even and $t\geqslant 3$ is odd, then
$$
  s_t(x)-s_1(x)=\frac{1}{2}(t-1)(m+2x-t-2).
$$
Hence $s_t(x)-s_1(x)>0$ if $x\geqslant \lambda(S_{\frac{m+4}{2},2}^1).$ Therefore, $\lambda(S_{\frac{m+t+3}{2},2}^t)<\lambda(S_{\frac{m+4}{2},2}^1).$

If $m$ is odd and $t\geqslant 4$ is even, then
$$
  s_t(x)-s_2(x)={\frac{1}{2}(t-2)(m+2x-t-3)}.
$$
Therefore, $s_t(x)-s_2(x)>0$ if $x\geqslant \lambda(S_{\frac{m+5}{2},2}^2).$ It follows that $\lambda(S_{\frac{m+t+3}{2},2}^t)<\lambda(S_{\frac{m+5}{2},2}^2).$
\end{proof}

Let $\hat{G}$ be in $\mathfrak{G}(m,C_5)$ (resp. $\mathfrak{G}(m,C_6)$) such that $\hat{\lambda}:=\lambda(\hat{G})$ is as large as possible. Based on Lemma~\ref{lem2.1}, we know that $\hat{G}$ is connected. Assume that ${\bf \hat{x}}$ is the Perron vector of $\hat{G}$ and $\hat{u}$ is  the extremal vertex of $\hat{G}.$ Note that $S_{\frac{m+4}{2},2}^1\in \mathfrak{G}(m,C_5)$ if $m$ is even, and $S_{\frac{m+5}{2},2}^2\in \mathfrak{G}(m,C_5)$ if $m$ is odd. Based on Lemma \ref{lem4.00}(ii), one has {if $m\geqslant 22$, then} 
$$
  \hat{\lambda}\geqslant \min\left\{\lambda(S_{\frac{m+4}{2},2}^1),\lambda(S_{\frac{m+5}{2},2}^2)\right\}>\frac{1+\sqrt{4m-7}}{2}\geqslant 5.
$$

\subsection{\normalsize Proof of Theorem \ref{thm1.3}}
{In this subsection, we give the proof of Theorem \ref{thm1.3}, and assume $\hat{G}$ is a graph in $\mathfrak{G}(m,C_5)$ such that $\hat{\lambda}:=\lambda(\hat{G})$ is as large as possible.} Since $\hat{G}$ is $C_5$-free, we know that $\hat{G}[N(\hat{u})]$ contains no path of length $3.$ Then the following result holds immediately.

\begin{lem}\label{lem4.1}
Each connected component of $\hat{G}[N(\hat{u})]$ is either a triangle or a star $K_{1,r}$ for some integer $r\geqslant 0,$ where $K_{1,0}$ is a singleton component.
\end{lem}
Recall that $N_0(\hat{u})=\{u\in N(\hat{u}): d_{N(\hat{u})}(u)=0\}.$ For convenience, let $W=V(\hat{G})\setminus N[\hat{u}],$ $W_0=\cup_{u\in N_0(\hat{u})}N_W(u)$ and $N_+(\hat{u})=N(\hat{u})\setminus N_0(\hat{u}).$ Note that $\hat{\lambda}>\frac{1+\sqrt{4m-7}}{2}.$ Hence ${\hat{\lambda}}^2-\hat{\lambda}>m-2.$ In view of \eqref{eq:1.1} and \eqref{eq:1.2}, one has
\begin{align}\label{eq:4.1}
  (m-2)x_{\hat{u}}&<({\hat{\lambda}}^2-\hat{\lambda})x_{\hat{u}}\leqslant |N(\hat{u})|x_{\hat{u}}+\sum_{u\in N_+(\hat{u})}(d_{N(\hat{u})}(u)-1)x_u-\sum_{u\in N_0(\hat{u})}x_u+\sum_{w\in W}d_{N(\hat{u})}(w)x_w\\\notag
  &\leqslant \left(|N(\hat{u})|+2e(N_+(\hat{u}))-|N_+(\hat{u})|-\sum_{u\in N_0(\hat{u})}\frac{x_u}{x_{\hat{u}}}+e(W,N(\hat{u}))\right) x_{\hat{u}}.
\end{align}
It follows from $e(N(\hat{u}))=e(N_+(\hat{u}))$ that
\[\label{eq:4.3}
  e(W)<e(N_+(\hat{u}))-|N_+(\hat{u})|-\sum_{u\in N_0(\hat{u})}\frac{x_u}{x_{\hat{u}}}+2.
\]
Based on Lemma \ref{lem4.1}, one has $e(N(\hat{u}))\leqslant |N_+(\hat{u})|$ and the equality holds if and only if all non-trivial connected components of $\hat{G}[N(\hat{u})]$ are triangles. Together with \eqref{eq:4.3}, we obtain $e(W)\leqslant 1.$

The following lemma gives a clearer local structure of the extremal graph $\hat{G}.$
\begin{lem}\label{lem4.2}
$e(W)=0$ and $W=W_0.$
\end{lem}
\begin{proof}
Suppose that $e(W)=1.$ Then let $w_1w_2$ be the unique edge of $\hat{G}[W].$ By Lemma~\ref{lem2.1}, one has $d(w)\geqslant 2$ for each vertex $w\in W.$ Hence $\min\{d_{N(\hat{u})}(w_1),d_{N(\hat{u})}(w_2)\}\geqslant 1.$ Recall that $\hat{G}$ is $C_5$-free. It follows that $w_1,w_2$ share the same unique common neighbor $u$ in $N(\hat{u}).$ Then $u$ is a cut vertex of $\hat{G},$ which contradicts Lemma~\ref{lem2.1}. Therefore, $e(W)=0.$

{The second part of this lemma can be proved by a similar discussion as Lemma \ref{cl5.2} and the fact that $\hat{G}\not\cong S_{\frac{m+3}{2},2}$.}
\end{proof} 


%
Let $G$ and $H$ be two connected graphs. Define $H\diamond G$ to be the graph obtained by joining a vertex with maximum degree in $G$ and each  vertex of $H$ with an edge.
\begin{lem}\label{lem4.3}
$\hat{G}\cong S_{\frac{m+4}{2},2}^1$ if $m$ is even and $\hat{G}\cong S_{\frac{m+5}{2},2}^2$ if $m$ is odd, or $\hat{G}\cong H_{t,s}\circ R_k$ for some bipartite graph $H_{t,s}$ and some integer $k\geqslant 1.$
\end{lem}
\begin{proof}
By \eqref{eq:4.3} and Lemma \ref{lem4.2}, one obtains
$$
  e(N_+(\hat{u}))>|N_+(\hat{u})|+\sum_{u\in N_0(\hat{u})}\frac{x_u}{x_{\hat{u}}}-2.
$$
Therefore, $e(N_+(\hat{u}))\geqslant |N_+(\hat{u})|-1.$ Together with Lemma \ref{lem4.1}, we have $e(N_+(\hat{u}))=|N_+(\hat{u})|-1$ or $e(N_+(\hat{u}))= |N_+(\hat{u})|.$ It follows that $\hat{G}[N_+(\hat{u})]$ contains at most one tree component. We proceed by considering the following two cases.

{\bf Case 1.} $\hat{G}[N_+(\hat{u})]$ does not contain {any} tree component. Note that $W=W_0.$ Hence $\hat{G}\cong H_{t,s}\circ R_k$ for some bipartite graph $H_{t,s}$ and some integer $k\geqslant 0.$ If $k=0,$ then $\hat{G}$ is a bipartite graph. Together with Lemma~\ref{thm4}, one has $\hat{\lambda}\leqslant \sqrt{m}<\frac{1+\sqrt{4m-7}}{2}$ if $m\geqslant 22,$ a contradiction. Hence $k\geqslant 1,$ as desired.

{\bf Case 2.} $\hat{G}[N_+(\hat{u})]$ contains exactly one tree component. Assume that $\hat{G}[N_+(\hat{u})]$ consists of $k$ triangles and a non-trivial tree component $K_{1,r}$ for some $r\geqslant 1.$ Denote by $V(K_{1,r})=\{u,u_1,\ldots,u_r\}$ with $d_{N(\hat{u})}(u)=r.$ Since $W=W_0,$ we obtain $x_{u_1}=\cdots=x_{u_r}.$ {In addition,}
$$
  \hat{\lambda}x_{u_1}=x_{\hat{u}}+x_u\ \ \text{and}\ \ \hat{\lambda}x_{u}=x_{\hat{u}}+rx_{u_1}.
$$
Therefore, $x_u=\frac{\hat{\lambda}+r}{{\hat{\lambda}}^2-r}x_{\hat{u}}$ and $x_{u_1}=\frac{\hat{\lambda}+1}{{\hat{\lambda}}^2-r}x_{\hat{u}}.$

Suppose that $k\geqslant 1.$ Let $\{\hat{u},v_1,v_2,v_3\}$ be {the vertex set of} a $K_4$ block of $\hat{G}.$ {Clearly,} $x_{v_1}=x_{v_2}=x_{v_3}.$ It follows that $\hat{\lambda}x_{v_1}=2x_{v_1}+x_{\hat{u}},$ that is, $x_{v_1}=\frac{x_{\hat{u}}}{\hat{\lambda}-2}.$ If $r\in\{1,2\},$ then by \eqref{eq:4.1} one has
\begin{align*}
(m-2-|N(\hat{u})|)x_{\hat{u}}&<\sum_{u\in N_+(\hat{u})}(d_{N(\hat{u})}(u)-1)x_u-\sum_{u\in N_0(\hat{u})}x_u+\sum_{w\in W}d_{N(\hat{u})}(w)x_w\\
&\leqslant \left(\frac{3k}{\hat{\lambda}-2}+\frac{(r-1)(\hat{\lambda}+r)}{{\hat{\lambda}}^2-r}\right)x_{\hat{u}}-\sum_{u\in N_0(\hat{u})}x_u+e(N(\hat{u}),W)x_{\hat{u}}\\
&\leqslant \left(\frac{3k+r-1}{\hat{\lambda}-2}+e(N(\hat{u}),W)\right)x_{\hat{u}}.
\end{align*}
Together with $e(W)=0$ and $\hat{\lambda}>5,$ we obtain that
$$
  3k+r-2=e(N_+(\hat{u}))-2<\frac{3k+r-1}{\hat{\lambda}-2}<\frac{3k+r-1}{3}<3k-1,
$$
a contradiction.

Next, we consider $r\geqslant 3.$ {We first give the following claim.
\begin{claim}\label{claim001}
  $\hat{G}\not\cong K_{1,\frac{m-7}{2}}\diamond K_4.$
\end{claim}
\begin{proof}[\bf Proof of Claim \ref{claim001}]
  Suppose that $\hat{G}\cong K_{1,\frac{m-7}{2}}\diamond K_4.$ Then $\pi':\, \{v_1,v_2,v_3\}\cup \{\hat{u}\}\cup \{u_1,\ldots,u_{\frac{m-7}{2}}\}\cup \{u\}$ is an equitable partition of $A(\hat{G})$. Therefore, the quotient matrix of $A(\hat{G})$ corresponds to the partition $\pi'$ is
\begin{equation*}
  A(\hat{G})_{\pi'}=\left(
                     \begin{array}{cccc}
                       2 & 1 & 0 & 0 \\
                       3 & 0 & \frac{m-7}{2} & 1 \\
                       0 & 1 & 0 & 1 \\
                       0 & 1 & \frac{m-7}{2} & 0 \\
                     \end{array}
                   \right).
\end{equation*}
By a direct calculation, the characteristic polynomial of $A(\hat{G})_{\pi}$ is $x^4-2x^3-(m-3)x^2+(m-5)x+\frac{7m-49}{2}$. Let $l(x)=x^4-2x^3-(m-3)x^2+(m-5)x+\frac{7m-49}{2}$ be a real function in $x$ for $x\in[\frac{1+\sqrt{4m-7}}{2},+\infty),$ and let $l^{(i)}(x)$ denote the $i$-th derivative of $l(x).$ It is easy to see that
$$
  l^{(1)}(x)=4x^3-6x^2-2(m-3)x+(m-5),\  l^{(2)}(x)=12x^2-12x-2(m-3)\ \text{and}\ l^{(3)}(x)=24x-12>0.
$$
Hence $l^{(2)}(x)$ is increasing in the interval $[\frac{1+\sqrt{4m-7}}{2},+\infty).$ Then
$$
  l^{(2)}(x)\geqslant l^{(2)}\left(\frac{1+\sqrt{4m-7}}{2}\right)=2(5m-9)>0.
$$
It follows that $l^{(1)}(x)$ is increasing in the interval $[\frac{1+\sqrt{4m-7}}{2},+\infty).$ Therefore,
$$
  l^{(1)}(x)\geqslant l^{(1)}\left(\frac{1+\sqrt{4m-7}}{2}\right)=(m-2)\sqrt{4m-7}-3>0.
$$
Thus, $l(x)$ is an increasing function in the interval $[\frac{1+\sqrt{4m-7}}{2},+\infty)$ and
$$
  l(x)\geqslant l\left(\frac{1+\sqrt{4m-7}}{2}\right)=\frac{1}{2}(7 m-3\sqrt{4 m-7}-52)>0.
$$
Hence, the largest eigenvalue of $ A(\hat{G})_{\pi}$ is less than $\frac{1+\sqrt{4m-7}}{2}$. Together with Lemma \ref{lem3.02}, we know $\hat{\lambda}<\frac{1+\sqrt{4m-7}}{2}$, a contradiction.
\end{proof}}

Now, let $G'=\hat{G}-\{v_iv_j:1\leqslant i<j\leqslant 3\}+\{v_iu:1\leqslant i\leqslant 3\}.$ Then ${G'}\in \mathfrak{G}(m,C_5)$ since $\hat{G}\not\cong K_{1,\frac{m-7}{2}}\diamond K_4.$
For $r\geqslant 3,$ one has
\begin{align*}
  \frac{x_u- x_{v_1}}{x_{\hat{u}}}=\frac{\hat{\lambda}+r}{{\hat{\lambda}}^2-r}-\frac{1}{\hat{\lambda}-2}
  =\frac{(\hat{\lambda}+r)(\hat{\lambda}-2)-({\hat{\lambda}}^2-r)}{({\hat{\lambda}}^2-r)(\hat{\lambda}-2)}
  =\frac{(r-2)\hat{\lambda}-r}{({\hat{\lambda}}^2-r)(\hat{\lambda}-2)}>0,
\end{align*}
the last inequality follows by $\hat{\lambda}>\max\{5,\sqrt{2r}\}$ {(since $K_{2,r}$ is a proper subgraph of $\hat{G}$ and $\lambda(K_{2,r})=\sqrt{2r}$).} Thus, by the Rayleigh quotient one has
$$
  \lambda(G')-\hat{\lambda}\geqslant {\bf \hat{x}}^T(A(G')-A(\hat{G})){\bf \hat{x}}=2x_u(x_{v_1}+x_{v_2}+x_{v_3})-6x_{v_1}^2=6x_{v_1}(x_u-x_{v_1})>0,
$$
i.e., $\lambda(G')>\hat{\lambda},$ which contradicts the choice of $\hat{G}.$


Therefore, $k=0$ and $\hat{G}\cong H_{t,s}\circ S_{r+2,2}$ for some bipartite graph $H_{t,s},$ where $t=|N_0(\hat{u})|$ and $s=|W_0|.$ In order to complete the proof, it suffices to show the following claim.
\begin{claim}\label{claim2}
  $s=0.$
\end{claim}
\begin{proof}[\bf Proof of Claim \ref{claim2}]
Suppose that $s\neq 0,$ i.e., $W_0\neq \emptyset$ and therefore $N_0(\hat{u})\neq \emptyset.$ If $r\in\{1,2\}$ and $m\geqslant 22,$ then by Lemmas~\ref{thm4} and \ref{lem4.0} one has
$$
  {\hat{\lambda}}\leqslant \sqrt{\lambda^2(G-u)+2r+1}\leqslant \sqrt{m+r}< \frac{1+\sqrt{4m-7}}{2},
$$
a contradiction. Hence $r\geqslant 3.$

Recall that $x_u=\frac{\hat{\lambda}+r}{{\hat{\lambda}}^2-r}x_{\hat{u}}.$ Together with \eqref{eq:4.1}, one has
$$
  (m-d(\hat{u})-2)x_{\hat{u}}<\left((r-1)\frac{\hat{\lambda}+r}{{\hat{\lambda}}^2-r}-\sum_{u\in N_0(\hat{u})}\frac{x_u}{x_{\hat{u}}}+e(N(\hat{u}),W)\right)x_{\hat{u}}<\left((r-1)\frac{\hat{\lambda}+r}{{\hat{\lambda}}^2-r}+e(N(\hat{u}),W)\right)x_{\hat{u}}.
$$
By Lemma \ref{lem4.2}, one obtains $r-2=e(N_+(\hat{u}))-2<(r-1)\frac{\hat{\lambda}+r}{{\hat{\lambda}}^2-r},$ that is, $\frac{\hat{\lambda}+r}{{\hat{\lambda}}^2-r}>\frac{r-2}{r-1}.$

On the other hand, note that $\hat{\lambda}>\frac{1+\sqrt{4m-7}}{2}.$ Let $h(x)=(x^2-r)(r-2)-(x+r)(r-1)=(x^2-x)(r-2)-x-2r^2+3r$ be a real function in $x$ for $x\in [\frac{1+\sqrt{4m-7}}{2},+\infty).$ By a direct calculation, we know that the derivation function of $h(x)$ is  $2(r-2)x-(r-1)>0$ and so $h(x)$ is increasing in the interval $[\frac{1+\sqrt{4m-7}}{2},+\infty).$ {Furthermore, in view of Lemma \ref{lem2.1}, one has $d(w)\geqslant 2$ for each $w\in W_0.$ Hence $|N_0(\hat{u})|\geqslant 2$ and $m\geqslant 2r+5.$ Note that $x^2-x=m-2$ if $x=\frac{1+\sqrt{4m-7}}{2}.$ Thus, 
\begin{align*}
  h(x)&\geqslant h\left(\frac{1+\sqrt{4m-7}}{2}\right)=(m-2)(r-2)-\frac{1+\sqrt{4m-7}}{2}-2r^2+3r\\
  &=m-2-\frac{1+\sqrt{4m-7}}{2}+(m-2)(r-3)-2r^2+3r\\
  &>(2r+3)(r-3)-2r^2+3r+9=0.
\end{align*}}
It follows that $\frac{\hat{\lambda}+r}{{\hat{\lambda}}^2-r}<\frac{r-2}{r-1},$ a contradiction.
\end{proof}


In view of Claim \ref{claim2}, one has $s=0$ and so $\hat{G}\cong S_{\frac{m+t+3}{2},2}^t.$ Together with Lemma \ref{lem4.00}(iii), we know that $\hat{G}\cong S_{\frac{m+4}{2},2}^1$ if $m$ is even, and $\hat{G}\cong S_{\frac{m+5}{2},2}^2$ if $m$ is odd. 

This completes the proof.
\end{proof}
Now, we are ready to prove Theorem \ref{thm1.3}.
\begin{proof}[\bf Proof of Theorem \ref{thm1.3}]
{In view of Lemma \ref{lem4.00}, in order to prove {Theorem \ref{thm1.3}}, we only need to prove $\hat{G}\cong S_{\frac{m+4}{2},2}^1$ if $m$ is even and $\hat{G}\cong S_{\frac{m+5}{2},2}^2$ if $m$ is odd. Based on Lemma \ref{lem4.3}, it suffices to show that $\hat{G}\not\cong H_{t,s}\circ R_k$ for any bipartite graph  $H_{t,s}$ and any integer $k\geqslant 1$.} Suppose to the contrary $\hat{G}\cong H_{t,s}\circ R_k,$  where $H_{t,s}=(T,S)$ is a bipartite graph and $k$ is a positive integer. Let $\{u,u_1,u_2,u_3\}$ be {the vertex set of} a $K_4$ block of $\hat{G},$ where $d(u)>3.$ Assume that $T=\{v_1,v_2,\ldots,v_t\}$ with $x_{v_1}\geqslant \cdots\geqslant x_{v_t},$ and $S=\{w_1,w_2,\ldots,w_s\}$ with $x_{w_1}\geqslant \cdots\geqslant x_{w_s}.$ By a similar discussion as Claim \ref{cla3.1}, we obtain that $x_v<x_u$ for all $v\neq u,$ i.e., $u=\hat{u}.$

{If $s=0,$ then $\hat{G}\cong tK_1\diamond R_k$. Similar to Claim \ref{claim001}, one may show $\hat{\lambda}<\frac{1+\sqrt{4m-7}}{2},$ a contradiction.} Therefore, $s\geqslant 1.$

By Lemma~\ref{lem2.1}, one has $d(w_i)\geqslant 2$ for each $w_i\in S.$ Then $t\geqslant 2$ and $e(T,S)\geqslant 2s.$ Together with Lemma~\ref{lem1.1}, we obtain that  $N_S(v_1)=N_S(v_2)=S,$ and $v_iw_j\in E(\hat{G})$ only if $v_pw_q\in E(\hat{G})$ for all $p\leqslant i$ and $q\leqslant j.$ {By a similar discussion as the main proof of Theorem \ref{thm1.04}(i), we can show the following facts:
\begin{wst}
  \item[{\rm (i)}] $s\geqslant 2$;
  \item[{\rm (ii)}] $e(T,S)\geqslant 2s+2\geqslant 6$ and $t\geqslant 3$, which implies $m=6k+t+e(T,S)\geqslant 6k+t+2s+2\geqslant 6k+2s+5$;
  \item[{\rm (iii)}] $k\geqslant 2$.
\end{wst}}



Now, {in order to obtain a contradiction,} we use induction on $m$ to show that $\hat{\lambda}\leqslant \frac{1+\sqrt{4m-7}}{2}.$ If $m=22,$ then $k=2$ and {$s=2.$} It is straightforward to check that $\hat{G}$ is isomorphic to one of the graphs depicted in Figure \ref{fig1}. By a simple calculation, we obtain $\hat{\lambda}<5,$ as desired.
\begin{figure}[!ht]
\centering
  \begin{tikzpicture}[scale = 0.7]
  \tikzstyle{vertex}=[circle,fill=black,minimum size=0.5em,inner sep=0pt]
  \node[vertex] (G_1) at (0,0){};
  \node[vertex] (G_2) at (-0.2,1){};
  \node[vertex] (G_3) at (-0.8,1.5){};
  \node[vertex] (G_4) at (-1.4,1){};
  \node[vertex] (G_5) at (0.2,1){};
  \node[vertex] (G_6) at (0.8,1.5){};
  \node[vertex] (G_7) at (1.4,1){};
  \node[vertex] (G_8) at (-1.5,-1){};
  \node[vertex] (G_9) at (-0.5,-1){};
  \node[vertex] (G_10) at (0.5,-1){};
  \node[vertex] (G_11) at (1.5,-1){};
  \node[vertex] (G_12) at (0,-2){};
  \node[vertex] (G_13) at (-1,-2){};
  \draw[thick] (G_1) -- (G_2)--(G_3)--(G_4)--(G_1);
  \draw[thick] (G_1) -- (G_5)--(G_6)--(G_7)--(G_1);
  \draw[thick] (G_3) --(G_1) -- (G_6);
  \draw[thick] (G_2) -- (G_4);
  \draw[thick] (G_5) -- (G_7);
  \draw[thick] (G_1) -- (G_8)-- (G_13)-- (G_9)-- (G_12)-- (G_10)-- (G_1);
  \draw[thick] (G_1) -- (G_9);
  \draw[thick] (G_1) -- (G_11);
  \draw[thick] (G_10) -- (G_13);
  \draw[thick] (G_8) -- (G_12);
  \end{tikzpicture}
\quad\quad\quad \quad\quad\quad
  \begin{tikzpicture}[scale = 0.7]
  \tikzstyle{vertex}=[circle,fill=black,minimum size=0.5em,inner sep=0pt]
  \node[vertex] (G_1) at (0,0){};
  \node[vertex] (G_2) at (-0.2,1){};
  \node[vertex] (G_3) at (-0.8,1.5){};
  \node[vertex] (G_4) at (-1.4,1){};
  \node[vertex] (G_5) at (0.2,1){};
  \node[vertex] (G_6) at (0.8,1.5){};
  \node[vertex] (G_7) at (1.4,1){};
  \node[vertex] (G_8) at (-1.5,-1){};
  \node[vertex] (G_9) at (-0.5,-1){};
  \node[vertex] (G_10) at (0.5,-1){};
  \node[vertex] (G_11) at (1.5,-1){};
  \node[vertex] (G_12) at (0,-2){};
  \node[vertex] (G_13) at (-1,-2){};
  \draw[thick] (G_1) -- (G_2)--(G_3)--(G_4)--(G_1);
  \draw[thick] (G_1) -- (G_5)--(G_6)--(G_7)--(G_1);
  \draw[thick] (G_3) --(G_1) -- (G_6);
  \draw[thick] (G_2) -- (G_4);
  \draw[thick] (G_5) -- (G_7);
  \draw[thick] (G_1) -- (G_8)-- (G_13)-- (G_9)-- (G_12);
  \draw[thick] (G_10)-- (G_1);
  \draw[thick] (G_1) -- (G_9);
  \draw[thick] (G_1) -- (G_11)-- (G_13);
  \draw[thick] (G_10) -- (G_13);
  \draw[thick] (G_8) -- (G_12);
  \end{tikzpicture}
  \caption{Graphs $H_{t,s}\circ R_k$ with $m=22.$}\label{fig1}
\end{figure}
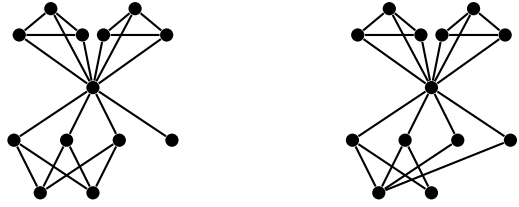

Assume the result holds if $|E(\hat{G})|\leqslant m-1.$ {For the case $|E(\hat{G})|= m$, we can prove $\hat{\lambda}\leqslant \frac{1+\sqrt{4m-7}}{2}$ by applying a similar discussion as that of Theorem \ref{thm1.04}(i).}

This completes the proof.
\end{proof}
\subsection{\normalsize Proof of Theorem \ref{thm1.03}}
{In this subsection, we give the proof of Theorem \ref{thm1.03}, and we assume $\hat{G}$ is a connected graph in $\mathfrak{G}(m,C_6)$ such that $\hat{\lambda}:={\lambda}(\hat{G})$ is as large as possible.}

Let $W=V(\hat{G})\setminus N[\hat{u}]$ and {$W_H= \cup_{v\in V(H)}N_W(v)$} for each component $H$ of $\hat{G}[N(\hat{u})].$ Since $\hat{G}$ is $C_6$-free, one has $W_{H_i}\cap W_{H_j}=\emptyset$ for any two distinct component $H_i$ and $H_j$ of $\hat{G}[N(\hat{u})],$ unless one of $H_i$ and $H_j$ is an isolated vertex and the other is a star $K_{1,r}$ for some $r\geqslant 0$ (in particular, vertices in $W_{H_i}\cap W_{H_j}$ are adjacent to the center vertex of $K_{1,r}$).

Recall that  $\hat{\lambda}>\frac{1+\sqrt{4m-7}}{2}\geqslant 5$ if $m\geqslant 22.$ Therefore, ${\hat{\lambda}}^2-\hat{\lambda}>m-2.$ By \eqref{eq:1.1} and \eqref{eq:1.2}, one has
$$
  (m-2-d(\hat{u}))x_{\hat{u}}<(\hat{\lambda}^2-\hat{\lambda}-d(\hat{u}))x_{\hat{u}}\leqslant \sum_{u\in N_+(\hat{u})}(d_{N(\hat{u})}(u)-1)x_u-\sum_{u\in N_0(\hat{u})}{x_u}+e(W,N(\hat{u}))x_{\hat{u}}.
$$
It follows that
$$
  \left(m-2-d(\hat{u})-e(W,N(\hat{u}))+\sum_{u\in N_0(\hat{u})}\frac{x_u}{x_{\hat{u}}}\right)x_{\hat{u}}<\sum_{u\in N_+(\hat{u})}(d_{N(\hat{u})}(u)-1)x_u.
$$
Put $\theta(H):=\sum_{u\in V(H)}(d_{H}(u)-1)x_u$ for each non-trivial connected component $H$ of $\hat{G}[N(\hat{u})].$ Hence
\[\label{eq:4.5}
  \left(e(N(\hat{u}))+e(W)-2+\sum_{u\in N_0(\hat{u})}\frac{x_u}{x_{\hat{u}}}\right)x_{\hat{u}}<\sum_H \theta(H),
\]
where $H$ takes over all non-trivial connected components of $\hat{G}[N(\hat{u})].$

Firstly, we give the following lemma, which considers the neighborhoods of vertices with degree two in $\hat{G}$.
\begin{lem}\label{cla1}
  Let $u_1$ and $u_2$ be two non-adjacent degree two vertices in $\hat{G}$ with $\sum_{u\in N(u_1)}x_u\geqslant \sum_{u\in N(u_2)}x_u.$ If $\hat{G}-\{u_2u:u\in N(u_2)\}+\{u_2u:u\in N(u_1)\}$ is not isomorphic to $S_{\frac{m+3}{2},2},$ then $N(u_1)=N(u_2).$
\end{lem}
\begin{proof}
Suppose that $N(u_1)\neq N(u_2).$ Let $G=\hat{G}-\{u_2u:u\in N(u_2)\}+\{u_2u:u\in N(u_1)\}.$ Then there exists a vertex, say $w$, in $N_{\hat{G}}(u_1)$ such that $N_{\hat{G}}(w)\subsetneqq N_G(w).$ It is easy to see that
$$
  {\bf \hat{x}}^TA(G){\bf \hat{x}}-{\bf \hat{x}}^TA(\hat{G}){\bf \hat{x}}=2x_{u_2}\left(\sum_{u\in N_{\hat{G}}(u_1)}x_u-\sum_{u\in N_{\hat{G}}(u_2)}x_u\right)\geqslant 0.
$$
Together with Lemma \ref{lem1.1}, one has $\lambda(G)>\hat{\lambda}.$

Next, we are to show that $G$ is $C_6$-free. In fact, if $H\cong C_6$ is a subgraph of $G,$ then $\{u_2u:u\in N_{\hat{G}}(u_1)\}\subseteq E(H).$ It is routine to check that $u_1\not\in V(H).$ Hence we can obtain a $C_6$ in $\hat{G}$ by replacing $u_2$ with $u_1$ in $H,$ a contradiction. It follows that $G\in \mathfrak{G}(m,C_6),$ which contradicts the choice of $\hat{G}.$
\end{proof}
{Note that $P_4,C_4,K_4,K_1\vee (C_4+e),P_5$ and $C_5$ are not induced subgraphs of $S_{\frac{m+3}{2},2}.$ On the other hand, each vertex of $S_{\frac{m+3}{2},2}$ {is of degree $2$ or $\frac{m-1}{2}.$} In particular, $S_{\frac{m+3}{2},2}$ contains no pendant vertices and cut vertices.  Hence, we can determine a graph that is not isomorphic to $S_{\frac{m+3}{2},2}$ according to its induced subgraphs and its vertex degrees.}  Next, we focus on characterizing the structure $\hat{G}[N(\hat{u})].$
\begin{lem}\label{lem4.4}
$\hat{G}[N(\hat{u})]$ is $\{C_4,C_5,\ldots\}$-free.
\end{lem}
\begin{proof}
Notice that $\hat{G}$ is $C_6$-free. Then $\hat{G}[N(\hat{u})]$ is $P_5$-free. It follows that $\hat{G}[N(\hat{u})]$ is $\{C_5,C_6,\ldots\}$-free. In particular, if $H$ is a component of $\hat{G}[N(\hat{u})]$ that contains $C_4$, then $H\in\{C_4,C_4+e,K_4\}$; if $H$ is a component of $\hat{G}[N(\hat{u})]$ with circumference $3$, then $H\cong K_{1,r}+e$ for some $r\geqslant 2.$ Let $\mathcal{H}$ be the family of components of $\hat{G}[N(\hat{u})]$ each of which contains  $C_4$ as a subgraph, and let $\mathcal{H}'$ be the family of other non-trivial components of $\hat{G}[N(\hat{u})]$. Then $\mathcal{H}'$ consists of trees and/or unicyclic graphs. Hence, for each $H\in \mathcal{H}',$ one has
$$
  \theta(H)=\sum_{u\in V(H)}(d_{H}(u)-1)x_u\leqslant (2e(H)-|V(H)|)x_{\hat{u}}\leqslant e(H)x_{\hat{u}}.
$$

Since $\hat{G}$ is $C_6$-free, one has $W_{H_i}\cap W_{H_j}=\emptyset$ and $e(W_{H_i},W_{H_j})=0$ if ${H_i},{H_j}\in \mathcal{H}$ and $H_i\neq H_j.$ Hence
$$
  e(W)\geqslant \sum_{H\in \mathcal{H}}e(W_H,W).
$$
Together with \eqref{eq:4.5}, one obtains
$$
  \left(\sum_{H\in \mathcal{H}}(e(H)+e(W_H,W))-2\right)x_{\hat{u}}<\sum_{H\in \mathcal{H}}\theta(H).
$$

Suppose to the contrary that $\mathcal{H}\neq \emptyset.$ In order to get a contradiction, it suffices to show that $\theta(H)\leqslant (e(H)+e(W_H,W)-2)x_{\hat{u}}$ holds for each $H\in \mathcal{H}.$

Let $H^*$ be in $\mathcal{H}$ with $V(H^*)=\{u_1,u_2,u_3,u_4\}.$ If $W_{H^*}=\emptyset,$ then let $x_{u_1}=\max\{x_{u_i}:1\leqslant i\leqslant 4\}.$ Hence  $\hat{\lambda}x_{u_1}=\sum_{u\in N(u_1)}x_u\leqslant x_{\hat{u}}+3x_{u_1},$ i.e., $x_{u_1}\leqslant \frac{x_{\hat{u}}}{\hat{\lambda}-3}<\frac{x_{\hat{u}}}{2}$ since $\hat{\lambda}> 5.$ Therefore,
$$
  \theta(H^*)=\sum_{u\in V(H^*)}(d_{H^*}(u)-1)x_u\leqslant (2e(H^*)-4)x_{u_1}< (e(H^*)-2)x_{\hat{u}},
$$
as desired. So, in what follows, we assume that $W_{H^*}\neq \emptyset.$ We proceed by distinguishing the following two possible cases.

{\bf Case 1.} All vertices in $W_{H^*}$ have a common neighbor, say $u_1,$ in $H^*.$ That is, $N_W(u_i)=\emptyset$ for $i\in\{2,3,4\}.$ Assume that $x_{u_2}=\max\{x_{u_i}:2\leqslant i\leqslant 4\}.$ 
Therefore, $\hat{\lambda}x_{u_2}{\leqslant }x_{\hat{u}}+x_{u_1}+x_{u_3}+x_{u_4}\leqslant 2(x_{\hat{u}}+x_{u_2}).$ It follows that $x_{u_2}\leqslant \frac{2x_{\hat{u}}}{\hat{\lambda}-2}<\frac{2x_{\hat{u}}}{3}.$ Hence
\begin{align*}
  \theta(H^*)&=\sum_{u\in V(H^*)}(d_{H^*}(u)-1)x_u\leqslant (d_{H^*}(u_1)-1)x_{u_1}+(2e(H^*)-d_{H^*}(u_1)-3)x_{u_2}\\
  &\leqslant (d_{H^*}(u_1)-1)x_{u_1}+\left(\frac{4}{3}e(H^*)-\frac{2}{3}d_{H^*}(u_1)-2\right)x_{\hat{u}}\leqslant \left(\frac{4}{3}e(H^*)-2\right)x_{\hat{u}}\leqslant e(H^*)x_{\hat{u}}.
\end{align*}
If $e(W_{H^*},W)\geqslant 2,$ then $\theta(H^*)\leqslant (e(H^*)+e(W_{H^*},W)-2)x_{\hat{u}},$ as desired. So, it suffices to consider the case  $e(W_{H^*},W)\leqslant 1.$

For the reason of $C_6$-free, one has $d_{N(\hat{u})}(w)=1$ for all $w\in W_{H^*}$ and $N_{H^*}(w)=N_{H^*}(w')$ if $ww'$ is an edge of $\hat{G}[W_{H^*}].$ Together with $W_{H^*}\neq \emptyset$ and Lemma~\ref{lem2.1}, one has $e(W_{H^*},W)=1,$ $|W_{H^*}|=1$ and $d(w)=2$ for the unique vertex $w\in W_{H^*}.$ Let $w'$ be the unique neighbor of $w$ in $W.$ Then
 \begin{equation}\label{eq:4.6}
    \left\{
    \begin{aligned}
        \hat{\lambda}x_{w}&=x_{u_1}+x_{w'},\\
        \hat{\lambda}x_{u_1}&\leqslant 3x_{u_2}+x_{w}+x_{\hat{u}},\\
        \hat{\lambda}x_{u_2}&\leqslant x_{u_1}+2x_{u_2}+x_{\hat{u}}.
    \end{aligned}
    \right.
\end{equation}
The first equation in \eqref{eq:4.6} implies $x_{w}\leqslant \frac{x_{u_1}+x_{\hat{u}}}{\hat{\lambda}}.$ The third inequality in \eqref{eq:4.6} means $x_{u_2}\leqslant \frac{x_{u_1}+x_{\hat{u}}}{\hat{\lambda}-2}.$ Together with the second inequality in \eqref{eq:4.6}, one has
$$
  x_{u_1}\leqslant \frac{{\hat{\lambda}}^2+2\hat{\lambda}-2}{{\hat{\lambda}}^3-2{\hat{\lambda}}^2-4\hat{\lambda}+2}x_{\hat{u}}\leqslant \frac{3}{4}x_{\hat{u}}.
$$
Hence $x_{u_2}\leqslant \frac{7x_{\hat{u}}}{4(\hat{\lambda}-2)}<\frac{7x_{\hat{u}}}{12}.$ Therefore,
\begin{align*}
  \theta(H^*)&\leqslant (d_{H^*}(u_1)-1)x_{u_1}+(2e(H^*)-d_{H^*}(u_1)-3)x_{u_2}\leqslant \left(\frac{3(d_{H^*}(u_1)-1)}{4}+\frac{7(2e(H^*)-d_{H^*}(u_1)-3)}{12}\right)x_{\hat{u}}\\
  &=\left(\frac{7e(H^*)+d_{H^*}(u_1)}{6}-\frac{5}{2}\right)x_{\hat{u}}\leqslant (e(H^*)-1)x_{\hat{u}}=(e(H^*)+e(W_{H^*},W)-2)x_{\hat{u}},
\end{align*}
as desired.

{\bf Case 2.} There exist $k\, (\geqslant 2)$ vertices, say $w_1,\ldots,w_k$, { in $W_{H^*}$ such that they} have mutual distinct neighbors in $V(H^*).$ Recall that $d_{N(\hat{u})}(w)=1$ for all $w\in W_{H^*}$ and $N_{H^*}(w)=N_{H^*}(w')$ if $ww'\in E(\hat{G}[W_{H^*}]).$ Hence  $\{w_1,\ldots,w_k\}$ is an independent set of $\hat{G}.$ Based on Lemma~\ref{lem2.1}, one has $d(w_i)\geqslant 2$ for $1\leqslant i\leqslant k.$ In addition, it follows from Lemma~\ref{cla1} that $d(w_i)=2$ holds for at most one $w_i.$ Therefore, $\sum_{1\leqslant i\leqslant k}d(w_i)\geqslant 3k-1$ and  $e(W_{H^*},W)\geqslant 2k-1.$ Thus,
$$
  e(H^*)\leqslant e(K_4)=6\leqslant e(W_{H^*},W)-2k+7\leqslant e(W_{H^*},W)+3
$$
and $e(H^*)=e(W_{H^*},W)+3$ holds if and only if $k=2$, $H^*\cong K_4$ and $e(W_{H^*},W)=2k-1=3.$
So,
$$
  \theta(H^*)=\sum_{u\in V(H^*)}(d_{H^*}(u)-1)x_u\leqslant (2e(H^*)-4)x_{\hat{u}}\leqslant (e(H^*)+e(W_{H^*},W)-2)x_{\hat{u}},
$$
unless $k=2$, $H^*\cong K_4$ and $e(W_{H^*},W)=3.$

Next, we consider the case $k=2$, $H^*\cong K_4$ and $e(W_{H^*},W)=3.$ Assume, without loss of generality, that $w_i\in N_{W_{H^*}}(u_i)$ for $i\in\{1,2\}$ and  $2\leqslant d(w_1)\leqslant d(w_2).$ Thus, $d(u_3)=d(u_4)=3$ and $x_{u_3}=x_{u_4}.$ {Note that $d_{N(\hat{u})}(w)=1$ for all $w\in W_{H^*}$ and $e(W_{H^*},W)=3.$ Then $d(w_1)=2$ and $d(w_2)\leqslant 3.$ In view of Lemma \ref{cla1}, we obtain $d(w_2)=3.$ It follows that $x_{w_1}\leqslant \frac{2x_{\hat{u}}}{\hat{\lambda}}$ and $x_{w_2}\leqslant \frac{3x_{\hat{u}}}{\hat{\lambda}}.$ Therefore, $x_{w_1}+x_{w_2}\leqslant \frac{5x_{\hat{u}}}{\hat{\lambda}}<x_{\hat{u}}$ since $\hat{\lambda}>5.$ 
If $W_{H^*}\setminus \{w_1,w_2\}\neq \emptyset,$ then let $w$ be such a vertex. Recall that $e(N_W(u_1),N_W(u_2))=0$ and $N_W(u_1)\cap N_W(u_2)=\emptyset.$ Together with Lemmas \ref{lem2.1} and \ref{cla1}, one has $N(w)=\{u_1,w_1\}$ and $N_W(u_1)=\{w,w_1\}.$ It follows that $u_1$ is a cut vertex, which contradicts Lemma \ref{lem2.1}. Thus, $W_{H^*}= \{w_1,w_2\}.$} Hence, 
\begin{equation}\label{eq:4.8}
    \left\{
    \begin{aligned}
        \hat{\lambda}x_{u_1}&=x_{u_2}+2x_{u_3}+x_{\hat{u}}+x_{w_1};\\
        \hat{\lambda}x_{u_2}&=x_{u_1}+2x_{u_3}+x_{\hat{u}}+x_{w_2};\\
        \hat{\lambda}x_{u_3}&=x_{u_1}+x_{u_2}+x_{u_4}+x_{\hat{u}}.
    \end{aligned}
    \right.
\end{equation}
By adding the first two equations in \eqref{eq:4.8}, one has $x_{u_1}+x_{u_2}< \frac{7x_{\hat{u}}}{\hat{\lambda}-1}.$ From the third equation in \eqref{eq:4.8}, one has $x_{u_3}=x_{u_4}\leqslant \frac{3x_{\hat{u}}}{\hat{\lambda}-1}.$  Therefore, $\sum_{1\leqslant i\leqslant 4}x_{u_i}< \frac{13x_{\hat{u}}}{\hat{\lambda}-1}.$ Recall that $\hat{\lambda}>5.$ Hence
$$
  \theta(H^*)=2\sum_{1\leqslant i\leqslant 4}x_{u_i}< \frac{26x_{\hat{u}}}{\hat{\lambda}-1}< 7x_{\hat{u}}=(e(H^*)+e(W_{H^*},W)-2)x_{\hat{u}},
$$
as desired.

This completes the proof.
\end{proof}
In view of Lemma \ref{lem4.4}, we know that each component of $\hat{G}[N(\hat{u})]$ is either a tree or a unicyclic graph $K_{1,r}+e$ for some $r\geqslant 2.$ Let $c$ be the number of non-trivial {tree} components of $\hat{G}[N(\hat{u})].$ Then
$$
  \sum_{H}\theta(H)=\sum_{H}\sum_{u\in V(H)}(d_{H}(u)-1)x_u\leqslant \sum_{H}(2e(H)-|V(H)|)x_{\hat{u}}=(e(N(\hat{u}))-c)x_{\hat{u}},
$$
where $H$ takes over all non-trivial connected components of $\hat{G}[N(\hat{u})].$ For convenience, denote $N_0:=N_0(\hat{u})$ and {$W_0:=\cup_{v\in N_0}N_W(v).$} Combining with \eqref{eq:4.5}, one has
\[\label{eq:4.07}
  e(W)<2-c-\sum_{u\in N_0}\frac{x_u}{x_{\hat{u}}}.
\]
Hence $e(W)\leqslant 1$ and $c\leqslant 1.$

\begin{lem}\label{lem4.5}
$e(W)=0.$
\end{lem}
\begin{proof}
Suppose that $e(W)=1.$ Let $w_1w_2$ be the unique edge of $\hat{G}[W].$ 
Applying Lemma~\ref{lem2.1} yields $d(w)\geqslant 2$ for each $w\in W.$ Then $d_{N(\hat{u})}(w_i)\geqslant 1$ for $i\in \{1,2\}.$ On the other hand, in view of \eqref{eq:4.07}, one obtains $c=0.$ That is to say, each non-trivial component of $\hat{G}[N(\hat{u})]$ is a unicyclic graph $K_{1,r}+e$ for some $r\geqslant 2.$ Hence each vertex $w$ in $W_H$ satisfies $N_{N(\hat{u})}(w)\subseteq V(H),$  where $H$ is a  non-trivial component of $\hat{G}[N(\hat{u})].$

Since $\hat{G}$ is $C_6$-free, one of the following holds:
\begin{wst}
\item[{\rm (i)}] $w_1,w_2\in W_H$ for some non-trivial component $H$ of $\hat{G}[N(\hat{u})]$ and they must have a unique common neighbor $u$ in $V(H);$
\item[{\rm (ii)}] $N_{N(\hat{u})}(w_1),N_{N(\hat{u})}(w_2)\subseteq N_0.$ 
\end{wst}
If item (i) holds, then $u$ is a cut vertex, which contradicts Lemma~\ref{lem2.1}. Hence item (ii) holds. Then $|N(w_1)\cap N(w_2)|\leqslant 2.$ Otherwise, one may assume $\{v_1,v_2,v_3\}\subseteq N(w_1)\cap N(w_2),$ which implies that $\hat{u}v_1w_1v_2w_2v_3\hat{u}$ is a $C_6$ of $\hat{G},$ a contradiction. We proceed by distinguishing the following three cases.

{\bf Case 1.} $|N(w_1)\cap N(w_2)|=2.$ In this case, $d_{N(\hat{u})}(w_1)=d_{N(\hat{u})}(w_2)=2$ since $\hat{G}$ is $C_6$-free.  Assume that $N(w_1)\cap N(w_2)=\{v_1,v_2\}.$ We claim that $N_W(v_1)=N_W(v_2).$ Otherwise, there exists a vertex $w\in (N_W(v_1)\setminus N_W(v_2))\cup (N_W(v_2)\setminus N_W(v_1)).$ Assume, without loss of generality, that $w\in N_W(v_1)\setminus N_W(v_2).$ Recall that $d(w)\geqslant 2.$ Then we can find a vertex, say $v_3,$ in $N_{N_0}(w).$ It is obvious that $\hat{u}v_2w_2v_1wv_3\hat{u}$ is a $C_6$ of $\hat{G},$ a contradiction. Similarly, we can show that each vertex $w$ in $N_W(v_1)$ satisfies $N_{N(\hat{u})}(w)=\{v_1,v_2\}.$  Let $G_1=\hat{G}-\{wv_1:w\in N_W(v_1)\}+\{w\hat{u}:w\in N_W(v_1)\}.$ Clearly, $G_1\in \mathfrak{G}(m,C_6)$ and $\lambda(G_1)>\hat{\lambda}$ (in view of Lemma~\ref{lem1.1}), a contradiction.

{\bf Case 2.} $|N(w_1)\cap N(w_2)|=1.$ Assume that $N(w_1)\cap N(w_2)=\{v_1\}.$ By Lemma~\ref{lem2.1}, we know that $v_1$ is not a cut vertex. Hence $\max\{d_{N_0}(w_1),d_{N_0}(w_2)\}\geqslant 2.$ Assume, without loss of generality, that $d_{N_0}(w_1)\geqslant 2.$ Thus, there exists a vertex $v_2\in N_0$ such that $w_1$ and $v_2$ are adjacent. It follows that $d(w_2)=2.$ Otherwise, there exits a vertex $v_3\in N_0$ that is adjacent to $w_2.$ Then $\hat{u}v_3w_2v_1w_1v_2\hat{u}$ is a $C_6$ in $\hat{G},$ a contradiction. Applying Lemma~\ref{cla1} yields that each vertex out of $N(v_1)\cap N(w_1)$ has degree at least $3.$ Therefore, there is a vertex $w_3\in N_W(v_2)$ with $d(w_3)\geqslant 3.$ It follows that $N_{N(\hat{u})}(w_3)\setminus \{v_1,v_2\}\neq \emptyset.$ Then let $v_4$ be a vertex in $N_{N(\hat{u})}(w_3)\setminus \{v_1,v_2\}$. Hence $\hat{u}v_1w_1v_2w_3v_4\hat{u}$ is a $C_6$ in $\hat{G},$ a contradiction. 


{\bf Case 3.} $|N(w_1)\cap N(w_2)|=0.$ Firstly, we consider the case that $N_{N_0}(w_1)\cup N_{N_0}(w_2)$ contains a vertex, say $v_1,$ with degree two in $\hat{G}$. Without loss of generality, assume $v_1\in N_{N_0}(w_1).$ Let $v_2$ be a vertex in $N_{N_0}(w_2).$ Then Lemma~\ref{cla1} implies that all vertices out of {$N(v_1)\cup (N(\hat{u})\cap N(w_1))$} have degree at least three. Thus, $d(v_2)\geqslant 3,$ i.e., there exists a vertex $w_3\in N_W(v_2)$ with $d(w_3)\geqslant 3.$ Hence there are two distinct vertices $v_3,v_4\in N_{N_0}(w_3)$ and each of which has degree at least $3.$ Note that $v_3$ is not adjacent to $w_1$ (resp. $w_2$). Otherwise, $\hat{u}v_1w_1v_3w_3v_4\hat{u}$ (resp. $\hat{u}v_2w_2v_3w_3v_4\hat{u}$) is a $C_6$ in $\hat{G},$ a contradiction. So, $N_W(v_3)\setminus \{w_1,w_2,w_3\}\neq \emptyset$ and let $w_4$ be such a vertex. Then $d(w_4)\geqslant 3.$ Note that $w_4v_2,w_4v_4\not\in E(\hat{G}).$ Otherwise, either $\hat{u}v_2w_4v_3w_3v_4\hat{u}$ or $\hat{u}v_4w_4v_3w_3v_2\hat{u}$ is a $C_6$ in $\hat{G},$ a contradiction. Hence there exists a vertex $v_5$ in {$N_0$} such that $w_4v_5\in E(\hat{G}).$ Whereas, $\hat{u}v_5w_4v_3w_3v_2\hat{u}$ is a $C_6$ of $\hat{G}$, which is also a contradiction. 

Now, we assume that all vertices in $N_{N_0}(w_1)\cup N_{N_0}(w_2)$ have degree at least three. For $i\in \{1,2\}$, let $v'_i\in N_{N_0}(w_i)$ and $w'_{i+2}$ be a vertex in $N_W(v'_i)\setminus \{w_1,w_2\}.$ If $w'_3=w'_4$ and $\max\{d(w_1),d(w_2)\}\geqslant 3,$ then we can find a $C_6$ in $\hat{G},$ a contradiction. If $w'_3=w'_4$ and ${d(w_1)=d(w_2)}=2,$ then by Lemma \ref{cla1}, one obtains that all vertices in $V(\hat{G})\setminus \{w_1,w_2\}$ have degree at least $3.$ Hence $d(w'_3)\geqslant 3,$ which implies that there exists a vertex $v'_3\in N_0$ such that $v'_3w'_3\in E(\hat{G}).$ Note that $d(v'_3)\geqslant 3.$ Then there is a vertex $w'_5\in N_W(v'_3)$ with $d(w'_5)\geqslant 3.$ Hence  $w'_5$ is not adjacent to $v'_2$. Otherwise, $\hat{u}v'_2w'_5v'_3w'_3v'_1\hat{u}$ is a $C_6$ in $\hat{G},$ a contradiction. Thus, there exists a vertex $v'_4\in N_0$ such that $w'_5v'_4\in E(\hat{G}).$ Then $\hat{u}v'_4w'_5v'_3w'_3v'_2\hat{u}$ is a $C_6$ in $\hat{G},$ a contradiction. It follows that $w'_3\neq w'_4,$ i.e., $N_W(v'_1)\cap N_W(v'_2)=\emptyset.$

For $i\in\{3,4\},$ since $d(w'_i)\geqslant 2,$ there exists a vertex $v''_i\in N_0$ such that $w'_i$ and $v''_i$ are adjacent in $\hat{G}.$ Note that $v''_3\neq v''_4,$ {which also implies $N_{N_0}(w'_3)\cap N_{N_0}(w'_4)=\emptyset.$} Otherwise, $\hat{u}v'_2w'_4v''_3w'_3v'_1\hat{u}$ is a $C_6$ in $\hat{G},$ a contradiction. It follows from Lemma \ref{cla1} that $\max\{d(w'_3),d(w'_4)\}\geqslant 3.$ Without loss of generality, assume that $d(w'_3)\geqslant 3.$ Hence there is a vertex $v''_5\in N_0$ such that $w'_3v''_5\in E(\hat{G}).$ It is easy to see that $v''_3$ is not adjacent to $w_1,w_2,w'_4.$ If $d(v''_3)\geqslant 3,$ then there is a vertex $w''_5\in N_W(v''_3).$  Notice that $w''_5$ is not adjacent to $v'_1,v'_2,v''_4,v''_5.$ Hence we can find a vertex $v''_6\in N_0$ such that $v''_6w''_5\in E(\hat{G}).$ It is obvious that $\hat{u}v''_6w''_5v''_3w'_3v'_1\hat{u}$ is a $C_6$ in $\hat{G},$ a contradiction. Hence $d(v''_3)=2.$ Applying Lemma~\ref{cla1} again  yields that $d(w'_4)\geqslant 3$ and $d(v''_4)\geqslant 3.$ Furthermore, $w'_4$ is not adjacent to $v'_1,v''_3,v''_5.$ Then there is a vertex $v''_7\in N_0$ such that $v''_7w'_4\in E(\hat{G}).$ By a similar discussion as before, we can also get a contradiction.

This completes the proof.
\end{proof}
The following two results characterize the structure of $\hat{G}$ by considering $\hat{G}[N(\hat{u})]$ contains {triangles} and $\hat{G}[N(\hat{u})]$ is a forest, respectively.
\begin{lem}\label{lem4.6}
 If $\hat{G}[N(\hat{u})]$ contains {triangles}, then {$\hat{G}\cong K_1\vee (S^1_{\frac{m-1}{2}}\cup K_1)$ if $m$ is odd and $\hat{G}\cong K_1\vee S^1_{\frac{m}{2}}$} if $m$ is even.
\end{lem}
\begin{proof}
Assume that $\hat{G}[N(\hat{u})]$ contains a component $H^*\cong K_{1,r}+e$ for some $r\geqslant 2.$ We first consider the following three claims. {The first one can be proved by a similar discussion as Claim \ref{cl:5.1}, so we omit its proof here.}
\begin{claim}\label{cl4.1}
  $H^*\not\cong K_3.$
\end{claim}


\begin{claim}\label{cl4.2}
$W_{H^*}=\emptyset$ and $H^*$ is the unique non-trivial component of $\hat{G}[N(\hat{u})].$
\end{claim}
\begin{proof}[\bf Proof of Claim \ref{cl4.2}]
In view of Claim \ref{cl4.1}, we may assume that $H^*\cong K_{1,r}+e$ for some $r\geqslant 3.$ Since $\hat{G}$ is $C_6$-free, one obtains that $d_{N(\hat{u})}(w)=1$ for each $w\in W_{H^*}.$ On the other hand, Lemma~\ref{lem2.1} implies that $d(w)\geqslant 2$ for each $w\in W.$ Together with Lemma \ref{lem4.5}, one has $W_{H^*}=\emptyset.$

Next, we show that $H^*$ is the unique non-trivial component of $\hat{G}[N(\hat{u})].$ Note that there are $r-2$ vertices in $V(H^*)$ with degree two in $\hat{G}.$ By Lemma \ref{cla1}, we obtain that there exist no vertices with degree two out of $H^*.$ Therefore, $H^*$ is the unique component of $\hat{G}[N(\hat{u})]$ which contains triangles. Suppose that $\hat{G}[N(\hat{u})]$ contains a non-trivial tree component $H.$ Based on \eqref{eq:4.07}, one has $c=1,$ i.e., $H$ is the unique non-trivial tree component of $\hat{G}[N(\hat{u})].$ In addition, $W_H\neq \emptyset$ and $d(w)\geqslant 3$ for each $w\in W_H\cup V(H).$

{If $\diam(H)=2$ or $3$, then by a similar discussion as Claim \ref{cl:5.2}, we may get a contradiction.} 


If $\diam(H)=1,$ then $H\cong K_2$ with $V(H)=\{u_1,u_2\}.$ Note that $\min\{d(u_1),d(u_2)\}\geqslant 3.$ There exists a vertex $w_i$ with $d(w_i)\geqslant 3$ such that $u_iw_i\in E(\hat{G})$ for $i\in\{1,2\}.$ If $w_1=w_2,$ then we can find a vertex $u_3$ in $N_0\cap N(w_1)$ with $d(u_3)\geqslant 3.$ It follows that there is a vertex $w_{3}\in N_W(u_3)$ with $d(w_{3})\geqslant 3.$ Note that $w_3$ is not adjacent to $u_1$ (resp. $u_2$). Otherwise, $\hat{u}u_3w_3u_1w_1u_2\hat{u}$ (resp. $\hat{u}u_3w_3u_2w_1u_1\hat{u}$) is a $C_6$ in $\hat{G}$, a contradiction. Hence we can find a vertex $u_4\in N_0$ such that $w_3u_4\in E(\hat{G}).$ However, $\hat{u}u_1w_1u_3w_3u_4\hat{u}$ is a $C_6$ of $\hat{G},$ a contradiction. Thus, $w_1\neq w_2,$ i.e., $N_W(u_1)\cap N_W(u_2)=\emptyset.$

It follows that  there are two distinct vertices $u'_3,u'_4\in N_0$ with $\min\{d(u'_3),d(u'_4)\}\geqslant 3$ such that $w_1u'_3,w_2u'_4\in E(\hat{G}).$ Note that $u'_3$ is not adjacent to $w_2.$ Otherwise, $\hat{u}u_1w_1u'_3w_2u_2\hat{u}$ is a $C_6$ in $\hat{G},$ a contradiction. Then there is a vertex $w'_3\in N_W(u'_3)$ with $d(w'_3)\geqslant 3.$ Also, for the reason of $C_6$-free, one obtains that $w'_3$ is not adjacent to $u_2,u'_4.$  Then we can find a vertex $u'_5\in N_{N_0}(w'_3).$ Whereas, $\hat{u}u'_5w'_3u'_3w_1u_1\hat{u}$ forms a $C_6$ in $\hat{G},$ a contradiction.

This completes the proof of Claim \ref{cl4.2}.
\end{proof}
\begin{claim}\label{cl4.3}
$W_0=\emptyset.$
\end{claim}
\begin{proof}[\bf Proof of Claim \ref{cl4.3}]
Suppose that $W_0\neq\emptyset.$ Then $N_0\neq \emptyset.$ It follows from Lemmas~\ref{lem2.1} and \ref{cla1} that $d(w)\geqslant 2$ for each $w\in W_0$ and $d(u)\neq 2$ for each $u\in N_0\cup W_0.$ If $d(u)\geqslant 3$ for some $u\in N_0,$ then assume that $\{w_1,w_2\}\subseteq N_W(u)$ with $\min\{d(w_1),d(w_2)\}\geqslant 3.$ Hence there are two distinct vertices $u_1$ and $u_2$ in $N_0$ such that $u_1w_1,u_2w_2\in E(\hat{G}).$ It is obvious that $\hat{u}u_2w_2uw_1u_1\hat{u}$ forms a $C_6$ in $\hat{G},$ a contradiction. Hence $d(u)=1$ for each $u\in N_0.$ That is, $W_0=\emptyset.$

This completes the proof of Claim \ref{cl4.3}.
\end{proof}

Now we come back to show Lemma \ref{lem4.6}.

In view of Claims \ref{cl4.1}-\ref{cl4.3}, we have $\hat{G}\cong K_1\vee (H^*\cup (m-2r-2)K_1),$ where $H^*\cong K_{1,r}+e$ for some $r\geqslant 3.$ Without loss of generality, assume that $V(H^*)=\{u,u_1,u_2,\ldots,u_r\}$ with $d_{H^*}(u)=r$ and $d_{H^*}(u_1)=d_{H^*}(u_2)=2.$ Note that $\pi'':=(V(\hat{G})\setminus (\{\hat{u}\}\cup V(H^*)))\cup \{\hat{u}\}\cup \{u_1,u_2\}\cup \{u\}\cup \{u_3,\ldots,u_r\}$ is an equitable partition of $A(\hat{G}).$ Hence the quotient matrix of $A(\hat{G})$ with respect to the partition $\pi''$ is
\begin{equation*}
    A(\hat{G})_{\pi''}=\left(
       \begin{array}{ccccc}
         0 & 1 & 0 & 0 & 0 \\
         m-2r-2 & 0 & 2 & 1 & r-2 \\
         0 & 1 & 1 & 1 & 0 \\
         0 & 1 & 2 & 0 & r-2 \\
         0 & 1 & 0 & 1 & 0 \\
       \end{array}
     \right).
\end{equation*}
Denote by $f_r(x)$ the characteristic polynomial of $A(\hat{G})_{\pi}.$ Then
$$
  f_r(x)=x^5- x^4- (m-1)x^3-(2r-m+5)x^2-(2r^2-mr+4)x+2r^2-mr-2r+2m-4.
$$
It is easy to see that $f_r(x)-f_{r+1}(x)=2(x^2+x)+(4r-m)(x-1).$

Firstly, we assume that $m\leqslant 4r+11.$ {Based on Lemma \ref{lem3.02}, one obtains that the largest root of $f_r(x)=0$ equals to $\hat{\lambda}.$ {Clearly, $\hat{\lambda}>\lambda(K_2)=1.$} If $m\leqslant 4r,$ then $f_r(x)-f_{r+1}(x)>0$ for all $x\in(1,+\infty)$ and all positive integer $r.$ If $4r+1\leqslant m\leqslant 4r+11,$} then the discriminant of $f_r(x)-f_{r+1}(x)$ is $(4r-m+2)^2+8(4r-m)<0.$ Hence $f_r(x)-f_{r+1}(x)>0$ for all $x\in(1,+\infty)$ and all positive integer $r.$ Therefore, the largest eigenvalue of $\hat{G}$ is less than that of $K_1\vee ((K_{1,r+1}+e)) \cup (m-2r-4)K_1).$ By the choice of $\hat{G},$ we know that  {$\hat{G}\cong K_1\vee (S^1_{\frac{m-1}{2}}\cup K_1)$ if $m$ is odd and $\hat{G}\cong K_1\vee S^1_{\frac{m}{2}}$ if $m$ is even.}

Next, we consider the case $m\geqslant 4r+12$ and $3\leqslant r\leqslant 7.$ 
Recall that $\lambda(S_{\frac{m+4}{2},2}^1)$ and $\lambda(S_{\frac{m+5}{2},2}^2)$ are, respectively, the largest zeros of $s_1(x)$ and $s_2(x)$, where $s_1(x)$ and $s_2(x)$ are two real functions defined in Lemma \ref{lem4.00}. For $m\geqslant 4r+12$ and $3\leqslant r\leqslant 7$, using Mathematica 9.0 \cite{W2012} gives us
$$
  f_r(x)-x\cdot s_1(x)>0\ \ \text{and}\ \ f_r(x)-x\cdot s_2(x)>0.
$$
Combining with Lemma \ref{lem3.02}, one has $\lambda(S_{\frac{m+4}{2},2}^1)>\hat{\lambda}$ for even $m$,  and $\lambda(S_{\frac{m+5}{2},2}^2)>\hat{\lambda}$ for odd $m$, which contradicts the choice of $\hat{G}.$

Now, we assume $m\geqslant 4r+12$ and $r\geqslant 8.$ Note that $x_{u_1}=x_{u_2}$ and $x_{u_3}=\cdots=x_{u_r}.$ Hence $\hat{\lambda}x_{u_1}=x_{\hat{u}}+x_{u_2}+x_{u}$ and $\hat{\lambda}x_{u_3}=x_{\hat{u}}+x_{u}.$ Hence $x_{u_1}= \frac{x_u+x_{\hat{u}}}{\hat{\lambda}-1}$  and   $x_{u_3}=\frac{x_{\hat{u}}+x_{u}}{\hat{\lambda}}.$ Thus,
$$
  \hat{\lambda}x_{u}=x_{\hat{u}}+2x_{u_1}+(r-2)x_{u_3}=x_{\hat{u}}+ \left(\frac{2}{\hat{\lambda}-1}+\frac{r-2}{\hat{\lambda}}\right)(x_{\hat{u}}+x_{u}),
$$
which is equivalent to
$$
  x_{u}=\frac{{\hat{\lambda}}^2+(r-1)\hat{\lambda}-r+2}{{\hat{\lambda}}^3-{\hat{\lambda}}^2-r\hat{\lambda}+r-2}x_{\hat{u}}.
$$

Since $m\geqslant 4r+12>2r+12,$ one has $N_0\neq \emptyset.$ Let $G_3=\hat{G}-u_1u_2+u_0u,$ where $u_0$ is a vertex in $N_0.$ Recall that $\hat{\lambda}> \frac{1+\sqrt{4m-7}}{2}.$ If $r \geqslant 8,$ then by the Rayleigh quotient and Mathematica 9.0 \cite{W2012}, one has
\begin{align*}
  \lambda(G_3)-\hat{\lambda}&\geqslant {\bf \hat{x}}^T(A(G_3)-A(\hat{G})){\bf \hat{x}}=2(x_{u_0}x_u-x_{u_1}x_{u_2})=2\left(\frac{x_ux_{\hat{u}}}{\hat{\lambda}}-\left(\frac{x_u+x_{\hat{u}}}{\hat{\lambda}-1}\right)^2\right)\\
  &=\frac{2x_ux_{\hat{u}}}{({\hat{\lambda}-1})^2}\left(\left({\hat{\lambda}}+\frac{1}{\hat{\lambda}}-4\right)-\left(\frac{x_u}{x_{\hat{u}}}+\frac{x_{\hat{u}}}{x_{u}}\right)\right)>0.
\end{align*}
Hence $\lambda(G_3)>\hat{\lambda}$ if $r \geqslant 8$ and $m\geqslant 4r+12.$ Furthermore, it is routine to check that $G_3\in \mathfrak{G}(m,C_6),$ then we also get a contradiction.

This completes the proof.
\end{proof}
\begin{lem}\label{lem4.7}
If $\hat{G}[N(\hat{u})]$ is a forest, then $\hat{G}\cong S_{\frac{m+4}{2},2}^1$ if $m$ is even and $\hat{G}\cong S_{\frac{m+5}{2},2}^2$ if $m$ is odd.
\end{lem}
\begin{proof}
Assume that there are exactly $c$ non-trivial tree components in $\hat{G}[N(\hat{u})]$. In view of \eqref{eq:4.5}, one has $c\leqslant 1.$
If $c=0,$ then by Lemma \ref{lem4.5}, we obtain that $\hat{G}$ is bipartite. Combining with Lemma \ref{thm4}, one has $\hat{\lambda}\leqslant \sqrt{m}\leqslant \frac{1+\sqrt{4m-7}}{2}$ if $m\geqslant 22,$ a contradiction. Hence $c=1$ and let $H$ be the unique non-trivial tree component of $\hat{G}[N(\hat{u})].$ Note that $\hat{G}$ is $C_6$-free. Then $\diam(H)\leqslant 3.$

If $\diam(H)=3,$ then $H$ is a double star $D_{a,b}$. Assume that $V(D_{a,b})=\{u_1,\ldots,u_a,u,u',u'_1,\ldots,u'_b\}$ with  {$d_{D_{a,b}}(u)=a+1\geqslant b+1=d_{D_{a,b}}(u')\geqslant 2$}. Hence $d_{N(\hat{u})}(w)=1$ for each vertex $w\in W_H.$ Together with Lemma~\ref{lem2.1}, one has $d_W(w)\geqslant 1$ for each vertex $w\in W_H.$ It follows from Lemma \ref{lem4.5} that $W_H=\emptyset$. Therefore, $\hat{G}$ contains $a+b$ non-adjacent vertices with degree two in $V(H).$

If $b\geqslant 2$ or $N_0\neq \emptyset,$ then let $G_1=\hat{G}-u'u'_b+uu'_b.$ It is routine to check that $x_u\geqslant x_{u'}.$ Then Lemma~\ref{lem1.1} implies that $\lambda(G_1)>\hat{\lambda}.$ In addition, $G_1\in \mathfrak{G}(m,C_6),$ which contradicts the choice of $\hat{G}.$ Hence $b=1$ and $N_0= \emptyset.$ That is to say, $\hat{G}\cong K_1\vee D_{\frac{m-5}{2},1}$ and $m$ is odd. {Hence $m\geqslant 23$.} It is straightforward to check that $\pi_1:={\{\hat{u}\}}\cup \{u\}\cup \{u'\}\cup \{u'_1\}\cup \{u_1,\ldots,u_{\frac{m-5}{2}}\}$ is an equitable partition of $A(\hat{G}).$ Hence the quotient matrix of $A(\hat{G})$ corresponding to the partition $\pi_1$ can be written as follows:
\begin{equation*}
A(\hat{G})_{\pi_1}=\left(
  \begin{array}{ccccc}
    0 & 1 & 1 & 1 & \frac{m-5}{2} \\
    1 & 0 & 1 & 0 & \frac{m-5}{2} \\
    1 & 1 & 0 & 1 & 0 \\
    1 & 0 & 1 & 0 & 0 \\
    1 & 1 & 0 & 0 & 0 \\
  \end{array}
\right).
\end{equation*}
Let $G_2$ be the graph obtained from $\hat{G}-u'u_1'$ by attaching a new pendant vertex $v$ to $\hat{u}.$ {Clearly, $G_2\in \mathfrak{G}(m,C_6).$} It is easy to see that $\pi_2:={\{\hat{u}\}}\cup \{u\}\cup \{u'_1\}\cup\{v\}\cup \{u_1,\ldots,u_{\frac{m-5}{2}},u'\}$ is an equitable partition of $A(G_2).$ Hence the quotient matrix of $A(G_2)$ corresponding to the partition $\pi_2$ is
\begin{equation*}
A(G_2)_{\pi_2}=\left(
  \begin{array}{ccccc}
    0 & 1 & 1 & 1 & \frac{m-3}{2} \\
    1 & 0 & 0 & 0 & \frac{m-3}{2} \\
    1 & 0 & 0 & 0 & 0 \\
    1 & 0 & 0 & 0 & 0 \\
    1 & 1 & 0 & 0 & 0 \\
  \end{array}
\right).
\end{equation*}
By a direct calculation, we obtain
$$
  \det(xI-A(\hat{G})_{\pi_1})=x^5-mx^3-(m-1)x^2+\frac{3m-15}{2}x+m-5\ \text{and}\ \det(xI-A(G_2)_{\pi_2})=x^5-mx^3-(m-3)(x^2-x).
$$
Consequently,
$$
 \det(xI-A(\hat{G})_{\pi_1})-\det(xI-A(G_2)_{\pi_2})=-2x^2+\frac{m-9}{2}x+m-5.
$$
{Note that $\frac{1+\sqrt{4m-7}}{2}<\hat{\lambda}<\frac{1+\sqrt{4m-3}}{2}$ (the right inequality follows by Theorem \ref{thm004}). Let $h(x)=-2x^2+\frac{m-9}{2}x+m-5$ be a real function in $x$ for $x\in [\frac{1+\sqrt{4m-7}}{2},\frac{1+\sqrt{4m-3}}{2}].$ Hence 
$$
  h(x)\geqslant \min\left\{h\left(\frac{1+\sqrt{4m-7}}{2}\right),h\left(\frac{1+\sqrt{4m-3}}{2}\right)\right\}.
$$
On the other hand, by some calculations we find that if $m\geqslant 25$, then
\begin{align*}
    h\left(\frac{1+\sqrt{4m-7}}{2}\right)&=-2(m-2)+{\left(\frac{m-9}{2}-2\right)}\frac{1+\sqrt{4m-7}}{2}+m-5\\
    &\geqslant {-m-1+\frac{5(m-13)}{2}=\frac{3m-67}{2}>0}
\end{align*}
and
\begin{align*}
    h\left(\frac{1+\sqrt{4m-3}}{2}\right)&=-2(m-1)+{\left(\frac{m-9}{2}-2\right)}\frac{1+\sqrt{4m-3}}{2}+m-5\\
    &> {-m-3+\frac{5(m-13)}{2}=\frac{3m-71}{2}>0.}
\end{align*}
Therefore, $h(x)>0$ if $x\in [\frac{1+\sqrt{4m-7}}{2},\frac{1+\sqrt{4m-3}}{2}]$ and $m\geqslant 25.$ {Moreover, if $m=23$, then by a direct calculation, we may get that the largest eigenvalue of $A(\hat{G})_{\pi_1}$ is less than that of $A(G_2)_{\pi_2}$.} Together with Lemma \ref{lem3.02}, one has $\lambda(G_2)>\hat{\lambda}$ for all $m\geqslant 23,$ which contradicts the choice of $\hat{G}$.}

So, in what follows, we assume that $\diam(H)\leqslant 2.$ That is, $H$ is a star $K_{1,r}$ for some $r\geqslant 1.$ Assume that $V(H)=\{u_0,u_1,\ldots,u_r\}$ with $d_{H}(u_0)=r.$

If $N_0=\emptyset,$ then $\hat{\lambda}x_{\hat{u}}=\sum_{i=0}^rx_{u_i}\leqslant (r+1)x_{\hat{u}}.$ Hence $r\geqslant 4.$ Since $\hat{G}\not\cong S_{\frac{m+3}{2},2},$ one has $W_H\neq \emptyset.$ If there exists a vertex $w\in W_H$ such that $|N(w)\cap\{u_1,\ldots,u_r\}|\geqslant 2,$ then we can assume that $u_1,u_2\in N(w)$. In this case, $\hat{u}u_1wu_2u_0u_3\hat{u}$ is a $C_6$ of $\hat{G},$ a contradiction. Hence $|N(w)\cap\{u_1,\ldots,u_r\}|\leqslant 1$ for all $w\in W_H.$ Together with Lemmas~\ref{lem2.1} and \ref{lem4.5}, one has $d(w)=2$ and $u_0\in N(w)$ for all $w\in W_H.$ Hence
$$
  \hat{\lambda}(x_{u_0}-x_{\hat{u}})=x_{\hat{u}}+\sum_{w\in W_H}x_w-x_{u_0}>0.
$$
That is, $x_{u_0}-x_{\hat{u}}>0,$ which contradicts the choice of $\hat{u}.$

Next, we consider the case $N_0\neq \emptyset.$ {Suppose that $W_H\neq \emptyset$. Then Lemma~\ref{lem2.1} implies $d(w)\geqslant 2$ for each $w\in W_H.$} Assume $N_0=\{v_1,v_2,\ldots,v_t\}$  and $W_0=\{w_1,w_2,\ldots,w_s\}.$ We proceed by distinguishing the following three cases.

{\bf Case 1.} $r\geqslant 3.$ In this case, by a similar discussion as the case $N_0=\emptyset,$ we obtain that $u_0\in N(w)$ and $|N(w)\cap \{u_1,\ldots,u_r\}|\leqslant 1$ for all $w\in W_H.$ If there exists a vertex, say $w_1$ in $W_H$ such that $|N(w_1)\cap \{u_1,\ldots,u_r\}|= 1,$ then assume that $u_1\in N(w_1).$ Since $\hat{G}$ is $C_6$-free, one obtains that the neighborhood of each vertex in $W_H$ is $\{u_0,u_1\}.$ Let $G_3=\hat{G}-\{u_0w:w\in W_H\}+\{\hat{u}w:w\in W_H\}.$ Obviously, $G_3\in  \mathfrak{G}(m,C_6)$ and $\lambda(G_3)>\hat{\lambda}$ (based on Lemma~\ref{lem1.1}), a contradiction.

Now, assume that $|N(w)\cap \{u_1,\ldots,u_r\}|=0$ for all $w\in W_H.$ Then $d(u_i)=2$ for $i\in\{1,\ldots,r\}.$ Note that $N_0\neq \emptyset$ and $W_H\neq \emptyset.$ Based on Lemma~\ref{cla1}, we obtain that all vertices in $V(G)\setminus \{u_1,\ldots,u_r\}$ have degree at least $3.$ Let $w_1\in W_H$ with $\{u_0,v_1,v_2\}\subseteq N(w_1),$ where $v_1,v_2\in N_0.$ Note that $d(v_1)\geqslant 3.$ Therefore, there exists a vertex $w_2\in N_W(v_1).$ Clearly, $w_2$ is not adjacent to $u_0$ (resp. $v_2$). Otherwise, $\hat{u}u_0w_2v_1w_1v_2\hat{u}$ (resp. $\hat{u}u_0w_1v_1w_2v_2\hat{u}$) forms a $C_6$ in $\hat{G},$ a contradiction. Hence there is a vertex $v_3\in N_0$ such that $v_3w_2\in E(\hat{G}).$ It follows that $\hat{u}u_0w_1v_1w_2v_3\hat{u}$ is a $C_6$ in $\hat{G},$ a contradiction.


{\bf Case 2.} $r=2.$ Let $w$ be a vertex in $W_H.$ If $w$ is adjacent to $u_1,$ then $d_{N_0}(w)=0$ (since $\hat{G}$ is $C_6$-free). We claim that $u_2\in N_{H}(w).$ Otherwise, {Lemmas~\ref{lem2.1} and \ref{lem4.5} imply} $u_0\in N_H(w)$ and $d(w)=2.$ Based on Lemma~\ref{cla1}, one has $d(u_2)\geqslant 3.$ Thus, there exists a vertex $w'\in N_W(u_2)$ with $d(w')\geqslant 3.$ Hence $N_H(w')=V(H)$ and so $\hat{u}u_1wu_0w'u_2\hat{u}$ is a $C_6$ of $\hat{G},$ a contradiction. It follows that $N_W(u_1)\subseteq N_W(u_2).$ Similarly, we can show that $N_W(u_2)\subseteq N_W(u_1).$ That is, $N_W(u_1)= N_W(u_2).$

If there exists a vertex $w\in W_H$ with $N_H(w)=V(H),$ then $W_H=\{w\}.$ Let $G_4=\hat{G}-u_0w+\hat{u}w.$ Then $G_4\in  \mathfrak{G}(m,C_6)$ and $\lambda(G_4)>\hat{\lambda}$ (by Lemma~\ref{lem1.1}), a contradiction. Hence for each vertex $w\in W_H,$ either $N_H(w)=\{u_1,u_2\}$ or $N_H(w)=\{u_0\}$ holds. Furthermore, if $w$ is a vertex in $N_W(u_1),$ then $d(w)=2.$ 

{\bf Subcase 2.1.} $N_W(u_1)\neq \emptyset.$ Then Lemma \ref{cla1} implies that each vertex in $V(\hat{G})\setminus N_W(u_1)$ has degree at least $3.$ If $N_W(u_0)=\emptyset,$ then let $G_5=\hat{G}-\{u_1w:w\in N_W(u_1)\}+\{\hat{u}w:w\in N_W(u_1)\}.$ It is routine to check that $G_5\in  \mathfrak{G}(m,C_6)$ and $\lambda(G_5)>\hat{\lambda}$ (based on Lemma~\ref{lem1.1}), a contradiction.

Next, we consider that $N_W(u_0)\neq \emptyset.$ Assume that $w_{1}\in N_{W}(u_0).$ Hence $d(w_1)\geqslant 3$ and there are two distinct vertices, say $v_{1}$ and $v_{2},$ in $N(w_{1})\cap N_0.$ On the other hand, $\min\{d(v_{1}),d(v_{2})\}\geqslant 3.$ Then there is a vertex $w_2\in N_{W_0}(v_{1}).$ Furthermore, $w_2$ is not adjacent to $v_{2}$ (resp. $u_0$). Otherwise, $\hat{u}u_0w_1v_{1}w_2v_{2}\hat{u}$ (resp. $\hat{u}u_0w_2v_{1}w_1v_{2}\hat{u}$) is a $C_6$ in $\hat{G},$ a contradiction. Thus, we can find a vertex $v_{3}\in N_{N_0}(w_2).$ It follows that $\hat{u}u_0w_1v_{1}w_2v_{3}\hat{u}$  is a $C_6$ in $\hat{G},$ a contradiction.

{\bf Subcase 2.2.} $N_W(u_1)= \emptyset.$ Since $W_H\neq \emptyset,$ one has $N_W(u_0)\neq \emptyset.$ Note that $d(u_1)=d(u_2)=2.$ Then Lemma~\ref{cla1} implies that each vertex in $V(\hat{G})\setminus \{u_1,u_2\}$ has degree at least $3.$ By a similar discussion as that of Subcase~2.1, we can also get a contradiction.

{\bf Case 3.} $r=1.$ Note that $d(w)\geqslant 2$ for each $w\in W.$ We proceed by considering the following three subcases.

{\bf Subcase 3.1.} $|N_W(u_0)\cap N_W(u_1)|\geqslant 2.$ For the reason of $C_6$-free, we obtain that either $N(w)=\{u_0,u_1\}$ or $N(w)\cap\{u_0,u_1\}=\emptyset$ holds for each $w\in W.$ Let $G_6=\hat{G}-\{u_1w:w\in N_W(u_1)\}+\{\hat{u}w:w\in N_W(u_1)\}.$ Note that $N_0\neq \emptyset.$ Hence $G_6\not\cong S_{\frac{m+3}{2},2}$ and so $G_6\in \mathfrak{G}(m,C_6).$ By Lemma \ref{lem1.1}, one has $\lambda(G_6)>\hat{\lambda}$, which  contradicts the choice of $\hat{G}$.

{\bf Subcase 3.2.} $|N_W(u_0)\cap N_W(u_1)|=1.$ In this case, for the reason of $C_6$-free, one has $|N_W(u_0)|=|N_W(u_1)|=1.$ Let $w\in N_W(u_0)\cap N_W(u_1).$ If $N_0\cap N(w)=\emptyset$ or each vertex in $N_0\cap N(w)$ has degree two, then let $G_7=\hat{G}-u_1w+\hat{u}w.$ Clearly, $G_7\in \mathfrak{G}(m,C_6)$ and $\lambda(G_7)>\hat{\lambda}$ (in view of Lemma \ref{lem1.1}),  a contradiction. Thus, there exists a vertex $v_1\in N_0\cap N(w)$ with $d(v_1)\geqslant 3.$ Let $w_1$ be in $N_W(v_1)\setminus \{w\}.$ Note that $u_0,u_1\not\in N(w_1).$ Hence there is a vertex $v_2\in N_{N_0}(w_1)\setminus \{v_1\}.$ Thus, $\hat{u}u_1wv_1w_1v_2\hat{u}$ is a $C_6$ in $\hat{G},$ a contradiction. 

{\bf Subcase 3.3.} $|N_W(u_0)\cap N_W(u_1)|=0.$  Without loss of generality, assume that $d(u_0)\geqslant d(u_1).$ If $d(u_1)=2,$ then together with Lemma \ref{lem4.5} we obtain that $\hat{G}-u_1$ is a bipartite graph. Combining with Lemmas \ref{thm4} and \ref{lem4.0}, one has if $m\geqslant 22,$ then
$$
  \hat{\lambda}\leqslant \sqrt{\lambda^2(\hat{G}-u_1)+3}\leqslant \sqrt{m+1}<\frac{1+\sqrt{4m-7}}{2},
$$
a contradiction. Hence $d(u_0)\geqslant d(u_1)\geqslant 3$.

For the reason of $C_6$-free, each vertex in $N_W(u_0)$ has no common neighbor with vertices in $N_W(u_1).$ Let $w_{i+1}$ be a vertex in $N_W(u_i)$ and $v_{i+1}$ be a vertex in $N_{N_0}(w_{i+1})$ for $i\in\{0,1\}.$ Based on Lemma \ref{cla1}, we have $\max\{d(v_1),d(v_2)\}\geqslant 3$. Without loss of generality, assume that $d(v_1)\geqslant 3.$ That is, there is a vertex $w_3\in N_{W}(v_1).$ Note that $w_3$ is not adjacent to $u_1$ (resp. $v_2$). Otherwise, $\hat{u}v_1w_3u_1w_2v_2\hat{u}$ (resp. $\hat{u}v_1w_3v_2w_2 u_1\hat{u}$) is a $C_6$ in $\hat{G},$ a contradiction. If $d(w_3)\geqslant 3$ or $d(w_3)=2$ and $w_3u_0\not\in E(\hat{G}),$ then we can find a vertex $v_3\in N_0$ that is adjacent to $w_3.$ Therefore, $\hat{u}v_3w_3v_1w_1u_0\hat{u}$  forms a $C_6$ in $\hat{G},$ a contradiction. If $d(w_3)=2$ and $w_3u_0\in E(\hat{G}),$ then Lemma \ref{cla1} implies $d(v_2)\geqslant 3.$ Hence there is a vertex $w_4\in N_W(v_2)$ with $d(w_4)\geqslant 3.$ Similarly, we can show that $w_4$ is not adjacent to $u_0$ and $v_1$. Thus, there exists a vertex $v_4\in N_0$ that is adjacent to $w_4.$ Whereas, $\hat{u}v_4w_4v_2w_2u_1\hat{u}$  is a $C_6$ in $\hat{G},$ a contradiction. 


Combining with Cases 1-3, we know that $W_H=\emptyset.$ Based on Claim \ref{claim2}, one has $W_0=\emptyset.$ Therefore, $\hat{G}\cong S_{\frac{m+t+3}{2},2}^t.$ Together with Lemma~\ref{lem4.00}(iii), we know that $\hat{G}\cong S_{\frac{m+4}{2},2}^1$ if $m$ is even, and $\hat{G}\cong S_{\frac{m+5}{2},2}^2$ if $m$ is odd.

This completes the proof.
\end{proof}
Now, we are ready to prove Theorem \ref{thm1.03}.
\begin{proof}[\bf Proof of Theorem \ref{thm1.03}]
Let $s_1(x)$ and $s_2(x)$ be two functions defined in Lemma \ref{lem4.00}, i.e.,
$$
  s_1(x)=x^4-mx^2+2x-mx+\frac{m}{2}-1\ \text{and}\ s_2(x)=x^4-mx^2-(m-3)(x-1).
$$
Then $\lambda(S_{\frac{m+4}{2},2}^1)$ and $\lambda(S_{\frac{m+5}{2},2}^2)$ are the largest roots of $s_1(x)=0$ and $s_2(x)=0,$ respectively. Furthermore, by using quotient matrix, we derive that the spectral radii of $K_1\vee S_{\frac{m}{2}}^1$ and $K_1\vee (S_{\frac{m-1}{2}}^1\cup K_1)$ equal the largest zeros of $g_1(x)$ and $g_2(x)$ respectively, where
$$
  g_1(x)=x^4-x^3+x^2-mx^2-3x+m-6\ \text{and}\ g_2(x)=x^5-x^4+ x^3-m x^3-2x^2+\frac{3m-17}{2}x-\frac{m-7}{2}.
$$
Together with Lemmas \ref{lem4.6} and \ref{lem4.7}, we know that in order to prove Theorem \ref{thm1.03}, it suffices to compare the spectral radii of $S_{\frac{m+4}{2},2}^1$ and $K_1\vee S_{\frac{m}{2}}^1$ if $m$ is even, and the spectral radii of $S_{\frac{m+5}{2},2}^2$ and $K_1\vee (S_{\frac{m-1}{2}}^1\cup K_1)$ if $m$ is odd.

Firstly, we assume that $m\geqslant 22$ is even. Using Mathematica 9.0 \cite{W2012} gives us
$$
  s_1\left(\frac{1+\sqrt{4m-5}}{2}\right)<0\ \text{and}\ s_1(x)>0\ \text{if}\ x\geqslant \frac{1+\sqrt{4m-4}}{2}
$$
and
$$
  g_1\left(\frac{1+\sqrt{4m-7}}{2}\right)<0\ \text{and}\ g_1(x)>0\ \text{if}\ x\geqslant \frac{1+\sqrt{4m-3}}{2}.
$$
It follows that
$$
  \frac{1+\sqrt{4m-5}}{2}<\lambda(S_{\frac{m+4}{2},2}^1)<\frac{1+\sqrt{4m-4}}{2}
$$
{and}
$$
  \frac{1+\sqrt{4m-7}}{2}<\lambda(K_1\vee S_{\frac{m}{2}}^1)<\frac{1+\sqrt{4m-3}}{2}.
$$

It is straightforward to check that $s_1(x)-g_1(x)=x^3-x^2+ 5x-m x-\frac{m}{2}+5.$ Using Mathematica 9.0 \cite{W2012} gives us
$$
  s_1(x)-g_1(x)>0\ \text{for}\ x\in\left(\frac{1+\sqrt{4m-5}}{2},\frac{1+\sqrt{4m-4}}{2}\right)\ \text{and}\ 22\leqslant m\leqslant 72
$$
and
$$
  g_1(x)-s_1(x)>0\ \text{for}\ x\in\left(\frac{1+\sqrt{4m-7}}{2},\frac{1+\sqrt{4m-3}}{2}\right)\ \text{and}\ m\geqslant 90.
$$
Therefore, for even $m$, $\lambda(S_{\frac{m+4}{2},2}^1)<\lambda(K_1\vee S_{\frac{m}{2}}^1)$ if $22\leqslant m\leqslant 72$; and $\lambda(S_{\frac{m+4}{2},2}^1)>\lambda(K_1\vee S_{\frac{m}{2}}^1)$ if $m\geqslant 90.$ Furthermore, applying Mathematica 9.0 \cite{W2012} to calculate the spectral radii of the corresponding quotient matrices, we get $\lambda(S_{\frac{m+4}{2},2}^1)>\lambda(K_1\vee S_{\frac{m}{2}}^1)$ if $74\leqslant m\leqslant 88$. That is to say, for even $m,$ 
$$
 \hat{G}\cong \left\{
                \begin{array}{ll}
                  K_1\vee S_{\frac{m}{2}}^1, & \hbox{if $22\leqslant m\leqslant 72;$} \\
                  S_{\frac{m+4}{2},2}^1, & \hbox{if $m\geqslant 74.$}
                \end{array}
              \right.
$$

Now, we consider that $m\geqslant 23$ is odd.
Using Mathematica 9.0 \cite{W2012} gives us
$$
xs_2\left(\frac{1+\sqrt{4m-7}}{2}\right)<0\ \text{and}\ xs_2(x)>0\ \text{if}\ x\geqslant \frac{1+\sqrt{4m-6}}{2};
$$
and
$$
g_2\left(\frac{1+\sqrt{4m-7}}{2}\right)<0\ \text{and}\ g_2(x)>0\ \text{if}\ x\geqslant \frac{1+\sqrt{4m-5}}{2}.
$$
It follows that
$$
  \frac{1+\sqrt{4m-7}}{2}<\lambda(S_{\frac{m+5}{2},2}^2)<\frac{1+\sqrt{4m-6}}{2}
$$
{and}
$$
  \frac{1+\sqrt{4m-7}}{2}<\lambda(K_1\vee (S_{\frac{m-1}{2}}^1\cup K_1))<\frac{1+\sqrt{4m-5}}{2}.
$$

It is straightforward to check that  $xs_2(x)-g_2(x)=x^4-x^3+5x^2-m x^2-\frac{m-11}{2}x-\frac{m-7}{2}.$ Using Mathematica~9.0 \cite{W2012} again yields
$$
  xs_2(x)-g_2(x)>0\ \text{for}\ x\in\left(\frac{1+\sqrt{4m-7}}{2},\frac{1+\sqrt{4m-6}}{2}\right)\ \text{and}\ 23\leqslant m\leqslant 71;
$$
and
$$
  g_2(x)-xs_2(x)>0\ \text{for}\ x\in\left(\frac{1+\sqrt{4m-7}}{2},\frac{1+\sqrt{4m-5}}{2}\right)\ \text{and}\ m\geqslant 89.
$$
Therefore, for odd $m$, $\lambda(S_{\frac{m+5}{2},2}^2)<\lambda(K_1\vee (S_{\frac{m-1}{2}}^1\cup K_1))$ if $23\leqslant m\leqslant 71;$ and $\lambda(S_{\frac{m+5}{2},2}^2)>\lambda(K_1\vee (S_{\frac{m-1}{2}}^1\cup K_1))$ if  $m\geqslant 89$. In addition, using Mathematica 9.0 \cite{W2012} to calculate the spectral radii of the corresponding quotient matrices, we get that $\lambda(S_{\frac{m+5}{2},2}^2)>\lambda(K_1\vee (S_{\frac{m-1}{2}}^1\cup K_1))$ if $73\leqslant m\leqslant 87.$ That is to say, for odd $m,$ one has $
 \hat{G}\cong K_1\vee (S_{\frac{m-1}{2}}^1\cup K_1)$ if $23\leqslant m\leqslant 71$ and $\hat{G}\cong S_{\frac{m+5}{2},2}^2$ if $m\geqslant 73.$

This completes the proof. 
\end{proof}
\section{\normalsize Concluding remarks}
In this paper, we consider the Brualdi-Hoffman-Tur\'an type problem of graphs, which is interesting and challenging. Its discussion lies heavily on the structure of graphs.
Firstly the unique maximal graph among $\mathcal{G}(m,\theta_{1,2,3})$ with $m\geqslant 8$ and  $\mathcal{G}(m,\theta_{1,2,4})$ with $m\geqslant 22$ are characterized, respectively. 
Then, all the maximal graphs among $\mathcal{G}(m,C_5)$ (resp. $\mathcal{G}(m,C_6)$) excluding the book graph are characterized for $m\geqslant 26$, which extends the results of \cite{Lou,ZLS}.  

{Note that Zhai, Lin and Shu \cite{ZLS} also characterized graphs of size $m$ having the largest spectral radius among $\mathcal{G}(m,\{\theta_{1,2,2},\theta_{1,2,3}\})$, i.e.,
\begin{thm}[\cite{ZLS}]
Let $G$ be a graph of size $m.$ If $G\in \mathcal{G}(m,\{\theta_{1,2,2},\theta_{1,2,3}\})$ with $m\geqslant 9,$ then $\lambda(G)\leqslant \sqrt{m}.$ Equality holds if and only if $G$ is a complete bipartite graph or one of $S_7^3,\,S_8^2$ and $S_9^1$.
\end{thm}

Motivated by the last two main results in our current paper, it is natural to determine the graph having the maximum spectral radius among $\mathfrak{G}(m,\{\theta_{1,2,2},\theta_{1,2,3}\})$, where
$$
\mathfrak{G}(m,\{\theta_{1,2,2},\theta_{1,2,3}\})=\mathcal{G}(m,\{\theta_{1,2,2},\theta_{1,2,3}\})\setminus \{K_{a,b}:\ m=ab\}.
$$
In fact, we can solve this problem by applying methods similar to those of Theorems \ref{thm1.3} and \ref{thm1.03}, and we omit its procedure here (see Appendix for its detailed proof).
\begin{thm}\label{thm1.2}
 Let $G$ be in $\mathfrak{G}(m,\{\theta_{1,2,2},\theta_{1,2,3}\})$ with $m\geqslant 26.$ Then $\lambda(G)\leqslant \rho_3(m),$ where $\rho_3(m)$ is the largest zero of $x^3-x^2-(m-1)x+m-3.$ Equality holds if and only if $G\cong S_m^1.$
\end{thm}

Furthermore, based on Theorem \ref{thm1.04}, it is interesting to investigate the following problem.
\begin{pb}\label{pb2}
  Which graph attains the maximum spectral radius among $\mathcal{G}(m,\theta_{1,2,3})$ (resp. $\mathcal{G}(m,\theta_{1,2,4})$) excluding the book graph?
\end{pb}
Very recently, Fang and You \cite{FY} settled the former case in Problem \ref{pb2}, i.e., they determined all graphs having the maximum spectral radius among $\mathcal{G}(m,\theta_{1,2,3})\setminus \{S_{\frac{m+3}{2},2}\}$. For the latter case, if $m$ is even, it was solved by Liu and Wang \cite{Liu2022},  i.e., they characterized all graphs having the maximum spectral radius among $\mathcal{G}(m,\theta_{1,2,4})\setminus \{S_{\frac{m+3}{2},2}\}$ for even $m$.}

Notice that, in 2017, Nikiforov \cite{0007} introduced the \textit{$A_{\alpha}$-matrix} of a graph $G$, which is a convex combination of $D(G)$ and $A(G)$, i.e.,
$$
  A_{\alpha}(G)=\alpha D(G)+(1-\alpha)A(G),\ \ \ 0\leqslant \alpha\leqslant 1,
$$
where $D(G)={\rm diag}(d_1,d_2,\ldots, d_n)$ is the diagonal matrix whose diagonal entries are the vertex degrees in $G$.

It is interesting to consider the Brualdi-Solheid-Tur\'{a}n type problem and Brualdi-Hoffman-Tur\'an type problem on the above novel matrix of graphs. Very recently, Tian, Chen, and Cui~\cite{Tian-21} identified the graphs with no 4-cycle and 5-cycle having the largest $A_\alpha$-spectral radius. Naturally, it is interesting to generalize all the results of this paper to $A_\alpha$-matrix. We will do it in the near future.

\section*{\normalsize Appendix}
Now, we give the detailed proof of Theorem \ref{thm1.2}. Before giving the proof, we present the following four lemmas.

For $2\leqslant s\leqslant t$, let $K_{s,t}^-$ be the graph obtained from the complete bipartite graph $K_{s,t}$ by deleting an edge, and let $K_{s,t}^+$ denote the graph obtained from $K_{s,t}$ by adding a new edge $xy$, where $x$ is a new vertex and $y$ is a vertex in the part of order $s.$ \begin{lem}[\cite{Cheng}]\label{thm3}
Let $G$ be an isolated-free bipartite graph with $m\geqslant 3$ edges such that $\lambda(G)$ is as large as possible. Assume that $G$ is not a complete bipartite graph.
\begin{wst}
\item[{\rm (i)}] If $m$ is odd, then $G\cong K_{2,q}^-$ with $q=\frac{m+1}{2};$
\item[{\rm (ii)}] If $m$ is even, $m-1$ is a prime and $m+1$ is not a prime, then $G\cong K_{p_1,q_1}^-,$ where $p_1\geqslant 3$ is the least integer that divides $m+1$ and $q_1=\frac{m+1}{p_1};$
\item[{\rm (iii)}] If $m$ is even, $m+1$ is a prime and $m-1$ is not a prime, then $G\cong K_{p_2,q_2}^+,$ where $p_2\geqslant 3$ is the least integer that divides $m-1$ and $q_2=\frac{m-1}{p_2};$
\item[{\rm (iv)}] If $m$ is even and neither $m+1$ nor $m-1$ is a prime, then $G\in \{K_{p_1,q_1}^-,K_{p_2,q_2}^+\},$ where $p_i$ and $q_i$ are defined as  (ii)-(iii) for $i=1,2;$
\item[{\rm (v)}] If $m\geqslant 4,$ then both $m+1$ and $m-1$ are primes if and only if $\lambda(G)\leqslant \sqrt{\frac{m+\sqrt{m^2-4(m-1-\sqrt{m-1})}}{2}}.$
\end{wst}
\end{lem}

Similar to Lemma \ref{lem2.1}, we give the following result.
\begin{lem}\label{lem02.1}
Let $G$ be a graph having the maximum spectral radius among $\mathfrak{G}(m,\{\theta_{1,2,2},\theta_{1,2,3}\}),$ where $m\geqslant 3$. Then $G$ is connected. Furthermore, if $u$ is the extremal vertex of $G,$ then there exists no cut vertex in $V(G)\setminus\{u\},$ and so $d(v)\geqslant 2$ for all $v\in V(G)\setminus N[u].$
\end{lem}
\begin{proof}
At first, suppose to the contrary that $G$ is disconnected. Let $G_1$ and $G_2$ be two connected components of $G$ and ${\lambda}(G)=\lambda(G_1).$ Choose $v_1v_2\in E(G_2)$ and assume $v'$ is a vertex with minimum degree in $G_1$. {Let $G'$ be a graph obtained from $G-v_1v_2+v_1v'$ by deleting all of its isolated vertices.} Notice that both $\theta_{1,2,2}$ and $\theta_{1,2,3}$ are $2$-connected. Hence $G'$ is $\{\theta_{1,2,2},\theta_{1,2,3}\}$-free. Clearly, $G_1$ is a proper subgraph of $G'$ and $G'$ is not a complete bipartite graph. Therefore, $\lambda(G')>\lambda(G_1)={\lambda}(G),$ which contradicts the choice of $G$. Hence, $G$ is connected.

Suppose that there is at least one cut vertex of $G$ in $V(G)\setminus\{u\}.$ Let $B$ be an end-block of $G$ with $u\not\in V(B)$ and let $v\in V(B)$ be a cut vertex of $G.$  Denote by $D_{i,j}$ the double star which consists of two stars $K_{1,i+1}$ and $K_{1,j+1}$ joined together so that they share an edge. In what follows we are to show $G\not\cong D_{m-2,1}.$

Recall that $S_m^1$ is obtained from $K_{1,m-1}$ by adding an edge between its two pendant vertices. Therefore, {$S_m^1\in \mathfrak{G}(m,\{\theta_{1,2,2},\theta_{1,2,3}\}).$} Hence $\lambda(G)\geqslant \lambda(S_m^1).$ Assume that ${V(K_{1,m-1})}=\{u_1,u_2,\ldots,u_m\}$ and $u_1$ is the center vertex. Denote by ${\bf y}=(y_{u_1},\ldots,y_{u_m})^T$ the Perron vector of $K_{1,m-1}.$ By a direct calculation, we get $y_{u_1}=\frac{1}{\sqrt{2}}$ and $y_{u_2}=\cdots=y_{u_m}=\frac{1}{\sqrt{2(m-1)}}.$ On the other hand, by Rayleigh quotient, we know that
\[\label{eq:2.00}
  \lambda(S_m^1)=\max_{{\bf x}\in \mathbb{R}^m}\frac{{\bf x}^TA(S_m^1){\bf x}}{{\bf x}^T{\bf x}}\geqslant {\bf y}^TA(S_m^1){\bf y}={\bf y}^TA(K_{1,m-1}){\bf y}+2y_{u_2}^2=\sqrt{m-1}+\frac{1}{m-1}.
\]
By a direct calculation we get
$$
  \lambda(D_{m-2,1})=\sqrt{\frac{m+\sqrt{m^2-4m+8}}{2}}<\sqrt{m-1}+\frac{1}{m-1}\leqslant \lambda(S_m^1)\leqslant{\lambda}(G).
$$
{Hence, $G\not\cong D_{m-2,1}.$}

Let $G''=G-\{wv:w\in V(B)\cap N(v)\}+\{wu:w\in V(B)\cap N(v)\}.$ Since $G\not\cong D_{m-2,1},$ one has $G''\not\cong K_{1,m}.$ Together with the fact that $u$ is a cut vertex of $G'',$ we obtain $G''$ is not a complete bipartite graph. It is routine to check that $G''$ is $\{\theta_{1,2,2},\theta_{1,2,3}\}$-free. That is, $G''\in \mathfrak{G}(m,\{\theta_{1,2,2},\theta_{1,2,3}\}).$ In view of Lemma~\ref{lem1.1}, we have $\lambda(G'')>{\lambda}(G),$ a contradiction.

This completes the proof.
\end{proof}

The following lemma determines the spectral radius of $S_m^1.$
\begin{lem}\label{claim1}
For $m\geqslant 4,$ the spectral radius of $S_m^1$ is equal to both the largest root of $x^4-mx^2-2x+m-3=0$ and that of $x^3-x^2-(m-1)x+m-3=0.$
\end{lem}
\begin{proof}
Assume that $V(S_m^1)=\{r_1,r_2,\ldots,r_m\}$ with $d_{S_m^1}(r_1)=m-1$ and $d_{S_m^1}(r_2)=d_{S_m^1}(r_3)=2.$ Hence, both $\pi_3:=\{r_1\}\cup \{r_2\}\cup \{r_3\}\cup \{r_4,\ldots,r_m\}$ and $\pi_4:=\{r_1\}\cup \{r_2,r_3\}\cup \{r_4,\ldots,r_m\}$ are equitable partitions of $A(S_m^1).$ The corresponding quotient matrices of $A(S_m^1)$ are given, respectively, as
\begin{equation*}
    \text{$A(S_m^1)_{\pi_3}=\left(
       \begin{array}{cccc}
         0 & 1 & 1 & m-3 \\
         1 & 0 & 1 & 0 \\
         1 & 1 & 0 & 0 \\
         1 & 0 & 0 & 0 \\
       \end{array}
     \right)$\ \ \ and\ \ \   $A(S_m^1)_{\pi4}=\left(
       \begin{array}{ccc}
         0 & 2 & m-3 \\
         1 & 1 & 0 \\
         1 & 0 & 0 \\
       \end{array}
     \right)$}.
\end{equation*}
By some calculations we obtain the characteristic polynomial of $A(S_m^1)_{\pi_3}$ is $x^4-mx^2-2x+m-3,$ whereas the characteristic polynomial of $A(S_m^1)_{\pi_4}$ is $x^3-x^2-(m-1)x+m-3$. By Lemma \ref{lem3.02}, our result follows immediately.
\end{proof}

Let $\hat{G}$ be in $\mathfrak{G}(m,\{\theta_{1,2,2},\theta_{1,2,3}\})$ with $m\geqslant 26$ such that $\lambda(\hat{G})$ is as large as possible. Denote $\hat{\lambda}:=\lambda(\hat{G}).$ In view of Lemma~\ref{lem02.1}, one obtains that $\hat{G}$ is connected. Assume $\bf{\hat{x}}$ is the Perron vector of $\hat{G}$ and $\hat{u}$ is the extremal vertex of $\hat{G}.$

Since $\hat{G}$ is $\theta_{1,2,2}$-free, we know that $\hat{G}[N(\hat{u})]$ consists of isolated vertices and isolated edges. Hence $N_i(\hat{u})=\{u\in N(\hat{u}): d_{N(\hat{u})}(u)=i\}$ for $i\in\{0,1\}$ and  $N(\hat{u}) = N_0(\hat{u})\cup N_1(\hat{u}).$ Let $W=V(\hat{G})\setminus N[\hat{u}]$ and $N^2_i(\hat{u})=\{w\in N^2(\hat{u}): d_{N_i(\hat{u})}(w)\geqslant 1\}$ for $i\in \{0,1\}.$ Since $\hat{G}$ is $\{\theta_{1,2,2},\theta_{1,2,3}\}$-free, one obtains that  $N^2_0(\hat{u})\cap N^2_1(\hat{u})=\emptyset$ and  $d_{N(\hat{u})}(w)=1$ for all $w\in N^2_1(\hat{u}).$ Note that $S_m^1\in \mathcal{G}(m,\{\theta_{1,2,2},\theta_{1,2,3}\})$ and $S_m^1$ is not a complete bipartite graph. Together with \eqref{eq:2.00}, one has for $m\geqslant 26,$
$$
  \hat{\lambda}\geqslant \lambda(S_m^1)\geqslant \sqrt{m-1}+\frac{1}{m-1}>\sqrt{m-1}\geqslant 5.
$$
Combining with \eqref{eq:1.2}, we have
\begin{align}\label{eq:3.2}
  (m-1-d(\hat{u}))x_{\hat{u}}&<{\hat{\lambda}}^2x_{\hat{u}}{=} \sum_{u\in N_1(\hat{u})}x_u+\sum_{w\in N^2_1(\hat{u})}x_w+\sum_{w\in N^2_0(\hat{u})}d_{N(\hat{u})}(w)x_w\\\notag
  &\leqslant (2e(N_1(\hat{u}))+e(N^2(\hat{u}),N(\hat{u})))x_{\hat{u}}=(2e(N(\hat{u}))+e(W,N(\hat{u})))x_{\hat{u}}.
\end{align}
It follows that $e(W)-1<e(N(\hat{u})),$ i.e., $e(W)\leqslant e(N(\hat{u})).$ 

\begin{lem}\label{lem3.1}
If $m\geqslant 26,$ then $e(W)=0.$
\end{lem}
\begin{proof}
{Recall that $e(W)\leqslant e(N(\hat{u}))=e(N_1(\hat{u})).$} If $N_1(\hat{u})=\emptyset,$ then $e(W)=0,$ as desired. So, in what follows, we assume $N_1(\hat{u})\neq \emptyset$ and $E(\hat{G}[N_1(\hat{u})])=\{u_{2i-1}u_{2i}:1\leqslant i\leqslant t\}.$ Define $\varphi(u_i):=\sum_{w\in N_W(u_i)}d_W(w)$ for $1\leqslant i\leqslant 2t.$ Let
$$
  A_0=\{i:\varphi(u_{2i-1})+\varphi(u_{2i})=0\},\ A_1=\{i:\varphi(u_{2i-1})+\varphi(u_{2i})=1\}\ \text{and}\ A_2=\{i:\varphi(u_{2i-1})+\varphi(u_{2i})\geqslant 2\}.
$$
Since $\hat{G}$ is $\{\theta_{1,2,2},\theta_{1,2,3}\}$-free, one has $N_W(u_{i})\cap N_W(u_{j})=\emptyset$ for $1\leqslant i<j\leqslant 2t.$ Then
\[\label{eq:3.04}
  2e(W)=\sum_{w\in W}d_W(w)\geqslant \sum_{w\in N_1^2(\hat{u})}d_W(w)\geqslant |A_1|+2|A_2|.
\]
This gives us
\[\label{eq:3.4}
  e(W)\geqslant \frac{1}{2}|A_1|+|A_2|.
\]

In order to complete the proof, we show the following two claims.
\begin{claim}\label{cl3.1}
$x_{u_{2i-1}}+x_{u_{2i}}<x_{\hat{u}}$ holds for all $i\in A_0\cup A_1.$
\end{claim}
\begin{proof}[\bf Proof of Claim \ref{cl3.1}]
Recall that $d_{N(\hat{u})}(w)=1$ for each $w\in N_1^2(\hat{u}).$ Together with Lemma~\ref{lem02.1}, we have $d_{W}(w)\geqslant 1$ for each $w\in N_1^2(\hat{u}).$

If $i\in A_0,$ then {$N_W(u_{2i-1})\cup N_W(u_{2i})=\emptyset.$ Hence} $x_{u_{2i-1}}=x_{u_{2i}}=\frac{x_{\hat{u}}}{\hat{\lambda}-1}.$ Note that $\hat{\lambda}>5.$ Hence
$$
  x_{u_{2i-1}}+x_{u_{2i}}\leqslant \frac{2x_{\hat{u}}}{\hat{\lambda}-1}<\frac{2x_{\hat{u}}}{4}<x_{\hat{u}},
$$
as desired.

Now, we consider that $i\in A_1.$ {Recall that $N_W(u_{2i-1})\cap N_W(u_{2i})=\emptyset.$} Then $|N_W(u_{2i-1})\cup N_W(u_{2i})|=1.$ Without loss of generality, we assume that {$|N_W(u_{2i-1})|=1.$ Let} $w\in N_W(u_{2i-1})$ and $w'$ be the unique neighbor of $w$ in $W.$ Hence $\hat{\lambda}x_{u_{2i}}=x_{\hat{u}}+x_{u_{2i-1}}$ and $\hat{\lambda}x_w=x_{u_{2i-1}}+x_{w'}\leqslant \hat{\lambda}x_{u_{2i}}.$ It follows that
$$
  {\hat{\lambda}}^2 x_{u_{2i-1}}=\hat{\lambda}(x_{\hat{u}}+x_{u_{2i}}+x_w)\leqslant (\hat{\lambda}+2)x_{\hat{u}}+2x_{u_{2i-1}}.
$$
Therefore,
$$
  x_{u_{2i-1}}\leqslant \frac{\hat{\lambda}+2}{{\hat{\lambda}}^2-2}x_{\hat{u}}\ \text{and}\ x_{u_{2i}}=\frac{x_{\hat{u}}+x_{u_{2i-1}}}{\hat{\lambda}}\leqslant \frac{\hat{\lambda}+1}{{\hat{\lambda}}^2-2}x_{\hat{u}}.
$$
Together with $\hat{\lambda}>5,$ we obtain
$$
  x_{u_{2i-1}}+x_{u_{2i}}\leqslant \frac{2\hat{\lambda}+3}{{\hat{\lambda}}^2-2}x_{\hat{u}}<x_{\hat{u}},
$$
as desired.
\end{proof}
In view of Claim \ref{cl3.1}, one has
$$
  \sum_{u\in N_1(\hat{u})}x_u=\sum_{i\in A_0\cup A_1\cup A_2}(x_{u_{2i-1}}+x_{u_{2i}})\leqslant (|A_0|+|A_1|+2|A_2|)x_{\hat{u}}
  =(e(N_1(\hat{u}))+|A_2|)x_{\hat{u}}.
$$
Together with \eqref{eq:3.2},  we have
$$
  (m-d(\hat{u})-1)x_{\hat{u}}<(e(N_1(\hat{u}))+|A_2|+e(W,N(\hat{u})))x_{\hat{u}},
$$
which is equivalent to  $m-d(\hat{u})-1<e(N_1(\hat{u}))+|A_2|+e(W,N(\hat{u})).$ Hence $e(W)\leqslant |A_2|.$ Combining with \eqref{eq:3.4}, we obtain $A_1=\emptyset$ and $e(W)=|A_2|.$ Therefore, all inequalities in \eqref{eq:3.04} are equalities. That is, $d_W(w)=0$ for each vertex $w\in W\setminus N_1^2(\hat{u})$ and $\varphi(u_{2i-1})+\varphi(u_{2i})=2$ for each $i\in A_2.$ 

\begin{claim}\label{cl3.2}
$x_{u_{2i-1}}+x_{u_{2i}}< x_{\hat{u}}$ for each $i\in A_2.$
\end{claim}
\begin{proof}[\bf Proof of Claim \ref{cl3.2}]
Recall that $\varphi(u_{2i-1})+\varphi(u_{2i})=2$ for each $i\in A_2.$ So, we proceed by distinguishing the following two possible cases.

{\bf Case 1.} $\varphi(u_{2i-1})=1$ and $\varphi(u_{2i})=1.$ In this case, {based on the definition of $\varphi(u_{i}),$ we can assume that $N_W(u_{j})=\{w_{j}\}$ and $d_W(w_j)=1,$ where $j\in \{2i-1,2i\}.$} Hence 
\begin{equation*}
    \left\{
    \begin{aligned}
        \hat{\lambda}x_{u_{2i-1}}&=x_{\hat{u}}+x_{u_{2i}}+x_{w_{2i-1}};\\
        \hat{\lambda}x_{u_{2i}}&=x_{\hat{u}}+x_{u_{2i-1}}+x_{w_{2i}};\\
        \hat{\lambda}x_{w_{2i-1}}&{\leqslant 2x_{\hat{u}}};\\
        \hat{\lambda}x_{w_{2i}}&{\leqslant 2x_{\hat{u}}}.
    \end{aligned}
    \right.
\end{equation*}
By a direct calculation, we get
$$
  x_{u_{2i-1}}+x_{u_{2i}}\leqslant {\frac{2\hat{\lambda}+4}{{\hat{\lambda}}^2-\hat{\lambda}}}x_{\hat{u}}< x_{\hat{u}},
$$
the last inequality follows by $\hat{\lambda}>5.$

{\bf Case 2.} $\varphi(u_{2i-1})=2$ and $\varphi(u_{2i})=0,$ or $\varphi(u_{2i-1})=0$ and $\varphi(u_{2i})=2.$ Without loss of generality, we assume that $\varphi(u_{2i-1})=2$ and $\varphi(u_{2i})=0.$

If $u_{2i-1}$ is adjacent to a vertex with degree $2$ in $\hat{G}[W],$ then assume that $N_W(u_{2i-1})=\{w_{{2i-1}}\}.$  Hence
\begin{equation*}
    \left\{
    \begin{aligned}
        \hat{\lambda}x_{u_{2i-1}}&=x_{\hat{u}}+x_{u_{2i}}+x_{w_{2i-1}};\\
        \hat{\lambda}x_{u_{2i}}&=x_{\hat{u}}+x_{u_{2i-1}};\\
        \hat{\lambda}x_{w_{2i-1}}&\leqslant 3x_{\hat{u}}.
    \end{aligned}
    \right.
\end{equation*}
By a direct calculation, one has
$$
  x_{u_{2i-1}}+x_{u_{2i}}\leqslant \frac{2\hat{\lambda}+3}{{\hat{\lambda}}^2-\hat{\lambda}}x_{\hat{u}}< x_{\hat{u}},
$$
the last inequality follows by $\hat{\lambda}>5.$

If $u_{2i-1}$ is adjacent to two vertices with degree $1$ in $\hat{G}[W],$ then by a similar discussion as above, we can derive that
$x_{u_{2i-1}}+x_{u_{2i}}< x_{\hat{u}}.$

This completes the proof of Claim \ref{cl3.2}.
\end{proof}
Combining with \eqref{eq:3.2} and Claims \ref{cl3.1}-\ref{cl3.2}, one has
$$
  (m-d(\hat{u})-1)x_{\hat{u}}<(|A_0|+|A_2|+e(W,N(\hat{u})))x_{\hat{u}}=(e(N(\hat{u}))+e(W,N(\hat{u})))x_{\hat{u}},
$$
that is $e(W)-1<0,$ i.e., $e(W)=0.$

This completes the proof.
\end{proof}

In the following, we give the proof of Theorem \ref{thm1.2}.
\begin{proof}[\bf Proof of Theorem \ref{thm1.2}]
Based on Lemma \ref{claim1}, it suffices to show $\hat{G}\cong S_m^1.$ If $N_1(\hat{u})=\emptyset,$ then $\hat{G}$ is a bipartite graph with color classes $X\cup Y$. Assume, without loss of generality, that $|X|\leqslant |Y|$. Assume that $X=\{u_1,\ldots,u_p\}$ and $Y=\{v_1,\ldots,v_q\}.$ Let $g(x):=x^4-mx^2-2x+m-3$ be a real function in $x$. In view of Lemma \ref{claim1}, we know that  $\lambda(S_m^1)$ is the largest root of $g(x)=0.$ Note that $\hat{G}$ is not a complete bipartite graph. We proceed by considering the following four possible cases to show that $\hat{\lambda}<\lambda(S_m^1).$


{\bf Case 1.} $m$ is odd, or $m$ is even, $m-1$ is a prime and $m+1$ is not a prime. Based on Lemma \ref{thm3}{(i)-(ii)}, one has $\hat{G}\cong K^-_{p,q},$ where $p\geqslant 2$ is the least integer such that $p\mid (m+1)$ {and $q=\frac{m+1}{p}.$} Without of loss of generality, assume that $u_pv_q\not\in E(\hat{G}).$ Then  $\pi_5:=(X\setminus \{u_p\})\cup \{u_p\}\cup (Y\setminus\{v_q\})\cup \{v_q\}$ is an equitable partition of $A(K^-_{p,q}).$ It is routine to check that the quotient matrix of $A(K^-_{p,q})$ corresponding to the partition $\pi_1$ can be written as
\begin{equation*}
    A(K^-_{p,q})_{\pi_5}=\left(
  \begin{array}{cccc}
    0 & 0 & \frac{m+1}{p}-1 & 1 \\
    0 & 0 & \frac{m+1}{p}-1 & 0 \\
    p-1 & 1 & 0 & 0 \\
    p-1 & 0 & 0 & 0 \\
  \end{array}
\right).
\end{equation*}
By a direct calculation, we know that the characteristic polynomial of $A(K^-_{p,q})_{\pi_1}$ is $f_1(x)=x^4-mx^2-p-\frac{m+1}{p}+m+2.$ Hence $f_1(x)-g(x)=2x-p-\frac{m+1}{p}+5.$ In view of Lemma \ref{lem3.02}, if one hopes to prove $\lambda(K^-_{p,q})<\lambda(S_m^1),$ {it suffices to show $f_1(x)-g(x)>0$ if $x\geqslant \lambda(S_m^1).$ Then, we just need to prove} $\lambda(S_m^1)>\frac{p}{2}+\frac{m+1}{2p}-\frac{5}{2}.$ By Mathematica~9.0 \cite{W2012}, we get
$$
  \max\left\{g\left(\frac{p}{2}+\frac{m+1}{2p}-\frac{5}{2}\right):m\geqslant 26\ \text{and}\ 2\leqslant p\leqslant \sqrt{m+1}\right\}<0,
$$
which implies that $\lambda(S_m^1)>\frac{p}{2}+\frac{m+1}{2p}-\frac{5}{2}$ for $m\geqslant 26,$ as desired.

{\bf Case 2.} $m$ is even, $m+1$ is a prime and $m-1$ is not a prime. In view of Lemma \ref{thm3}{(iii)}, one has $\hat{G}\cong K^+_{p,q-1},$ where $p\geqslant 3$ is the least integer such that $p\mid (m-1)$ {and $q=\frac{m-1}{p}+1.$} Without of loss of generality, assume that $d(u_p)=q$ and $d(v_q)=1.$ Therefore, $\pi_6:=(X\setminus \{u_p\})\cup \{u_p\}\cup (Y\setminus\{v_q\})\cup \{v_q\}$ is an equitable partition of $A(K^+_{p,q-1}).$ Hence the quotient matrix of $A(K^+_{p,q-1})$ with respect to the partition $\pi_6$ is
\begin{equation*}
    A(K^+_{p,q-1})_{\pi_6}=\left(
  \begin{array}{cccc}
    0 & 0 & \frac{m-1}{p} & 0 \\
    0 & 0 & \frac{m-1}{p} & 1 \\
    p-1 & 1 & 0 & 0 \\
    0 & 1 & 0 & 0 \\
  \end{array}
\right).
\end{equation*}
It is straightforward to check that the characteristic polynomial of $A(K^+_{p,q-1})_{\pi_6}$ is $f_2(x)=x^4-mx^2-\frac{m-1}{p}+m-1.$ Hence $f_2(x)-g(x)=2x-\frac{m-1}{p}+2.$ In order to prove $\lambda(K^+_{p,q-1})<\lambda(S_m^1),$ {it suffices to show $f_2(x)-g(x)> 0$ for all $x\geqslant \lambda(S_m^1)$ (based on Lemma \ref{lem3.02}). Hence we just need to prove} $\lambda(S_m^1)>\frac{m-1}{2p}-1.$ Using Mathematica 9.0 \cite{W2012} gives us
$$
  \max\left\{g\left(\frac{m-1}{2p}-1\right):m\geqslant 26\ \text{and}\ 3\leqslant p\leqslant \sqrt{m-1}\right\}<0.
$$
Therefore, $\lambda(S_m^1)>\frac{m-1}{2p}-1$ for $m\geqslant 26,$ as desired.

{\bf Case 3.} $m$ is even and neither $m+1$ nor $m-1$ is a prime. In view of Lemma \ref{thm3}{(iv)}, one has $\hat{G}\in \{K^-_{p,q},K^+_{p,q-1}\},$ where $p\geqslant 3$ is the least integer such that $p\mid (m-1)$ or $p\mid (m+1).$ By Cases 1-2, we obtain that $\hat{\lambda}<\lambda(S_m^1).$

{\bf Case 4.} $m$ is even, both $m+1$ and $m-1$ are primes. In this case, $\hat{\lambda}\leqslant \sqrt{\frac{m+\sqrt{m^2-4(m-1-\sqrt{m-1})}}{2}}$ (based on Lemma \ref{thm3}{(v)}). By a direct calculation, one obtains that for $m\geqslant 26,$
$$
  \hat{\lambda}\leqslant \sqrt{\frac{m+\sqrt{m^2-4(m-1-\sqrt{m-1})}}{2}}<\sqrt{m-1}+\frac{1}{m-1}\leqslant \lambda(S_m^1).
$$

Together with Cases 1-4, we know that $\hat{\lambda}<\lambda(S_m^1),$ which contradicts the choice of $\hat{G}$.

In what follows, we consider the case $N_1(\hat{u})\neq \emptyset.$ Assume that $E(\hat{G}[N_1(\hat{u})])=\{u_{2i-1}u_{2i}:1\leqslant i\leqslant t\}.$ Note that $d_{N(\hat{u})}(w)=1$ for each $w\in N_1^2(\hat{u}).$ Together with Lemmas~\ref{lem02.1} and \ref{lem3.1}, one has $N_1^2(\hat{u})=\emptyset$ and so $W=N_0^2(\hat{u}).$ Hence $d(u_{2i-1})=d(u_{2i})=2$ and $x_{u_{2i-1}}=x_{u_{2i}}$ for $1\leqslant i\leqslant t.$ Therefore, $\hat{\lambda}x_{u_{2i-1}}=x_{u_{2i}}+x_{\hat{u}}.$ That is, $x_{u_{2i-1}}=\frac{x_{\hat{u}}}{\hat{\lambda}-1}.$ Then \eqref{eq:3.2} implies that
$$
  (m-d(\hat{u})-1)x_{\hat{u}}<\left(\frac{2e(N_1(\hat{u}))}{\hat{\lambda}-1}+e(W,N(\hat{u}))\right)x_{\hat{u}}.
$$
It follows that $-1=e(W)-1<(\frac{2}{\hat{\lambda}-1}-1)e(N_1(\hat{u}))<0$. If $e(N_1(\hat{u}))\geqslant 2,$ then $\hat{\lambda}<3+\frac{2}{e(N_1(\hat{u}))-1}\leqslant 5,$ which contradicts the fact that $\hat{\lambda}> 5.$ Hence $e(N_1(\hat{u}))=1.$ That is to say, $\hat{G}\cong H_{|N_0(\hat{u})|,|N_0^2(\hat{u})|}\circ C_3,$ where $H_{|N_0(\hat{u})|,|N_0^2(\hat{u})|}=(T,S)$ is a bipartite graph with $T=N_0(\hat{u})$ and $S=N_0^2(\hat{u})=W.$ Put $V(C_3):=\{\hat{u},u_1,u_2\}.$ Then $x_{u_1}=x_{u_2}=\frac{x_{\hat{u}}}{\hat{\lambda}-1}.$ We claim that $x_{w}<x_{u_1}$ for all $w\in S.$ Otherwise, there exists a vertex $w\in S$ such that $x_{w}\geqslant x_{u_1}.$ Let $G_1=\hat{G}-u_1u_2+wu_2.$ Then $G_1$ is a bipartite graph and so $G_1\in \mathfrak{G}(m,\{\theta_{1,2,2},\theta_{1,2,3}\})$. Applying Lemma~\ref{lem1.1} yields $\lambda(G_1)>\hat{\lambda},$ which contradicts the choice of $\hat{G}.$ 

Next, we show that $S=\emptyset.$ Otherwise $|S|\geqslant 1.$ If $|S|=1,$ then assume that $S=\{w\}.$ In view of \eqref{eq:3.2}, one obtains
\begin{align*}\notag
  d(w)x_{\hat{u}}=(m-d(\hat{u})-1)x_{\hat{u}}<x_{u_1}+x_{u_2}+d(w)x_w< \frac{d(w)+2}{\hat{\lambda}-1}x_{\hat{u}}< \frac{d(w)+2}{4}x_{\hat{u}}.
\end{align*}
Hence $d(w)<\frac{2}{3},$ a contradiction. Thus, $|S|\geqslant 2.$

Without loss of generality, we assume that $T=\{v_1,\ldots,v_t\}$ and $S=\{w_1,\ldots,w_s\}$ with $x_{v_i}\geqslant x_{v_j}$ and $x_{w_i}\geqslant x_{w_j}$ if $i<j.$ Based on Lemma~\ref{lem1.1}, one has $v_iw_j\in E(\hat{G})$ only if $v_pw_q\in E(\hat{G})$ for all $p\leqslant i$ and $q\leqslant j.$ Note that $e(W)=0$ and $d(w)\geqslant 2$ for each $w\in W.$ Then $e(T,S)\geqslant 2s\geqslant 4.$ 
It is routine to check that $\hat{\lambda}(x_{\hat{u}}-x_{w})\geqslant x_{u_1}+x_{u_2}=\frac{2x_{\hat{u}}}{\hat{\lambda}-1}$ for each $w\in S.$ That is, $x_{w}\leqslant (1-\frac{2}{\hat{\lambda}(\hat{\lambda}-1)})x_{\hat{u}}$ if $w\in S.$ For convenience, denote by $\tilde{e}:=e(T,S).$ In view of \eqref{eq:3.2}, one obtains
\begin{align}\notag
  \tilde{e}x_{\hat{u}}&=(m-d(\hat{u})-1)x_{\hat{u}}<x_{u_1}+x_{u_2}+\sum_{w\in S}d_{T}(w)x_w\\\notag
  &\leqslant \left(\frac{2}{\hat{\lambda}-1}+\sum_{w\in S}d_{T}(w)\left(1-\frac{2}{\hat{\lambda}(\hat{\lambda}-1)}\right)\right)x_{\hat{u}}\\\label{eq:3.5}
  &\leqslant \left(\frac{2}{\hat{\lambda}-1}+\tilde{e}\left(1-\frac{2}{\hat{\lambda}(\hat{\lambda}-1)}\right)\right)x_{\hat{u}}.
\end{align}
If $\hat{\lambda}\leqslant \tilde{e},$ then $\frac{2}{\hat{\lambda}-1}-\frac{2\tilde{e}}{\hat{\lambda}(\hat{\lambda}-1)}\leqslant 0.$ By \eqref{eq:3.5}, one has $\tilde{e}<\tilde{e},$ a contradiction. Hence $\hat{\lambda}>\tilde{e}.$

It follows form Theorem \ref{thm004}(i) that $\hat{\lambda}<\sqrt{m}.$ Then $|T|+3+\tilde{e}=m> \hat{\lambda}^2>{\tilde{e}^2},$ that is, $t+2\geqslant \tilde{e}(\tilde{e}-1).$ Recall that $x_{v_1}\geqslant x_{v_2}\geqslant \cdots\geqslant x_{v_t}\geqslant \frac{x_{\hat{u}}}{\hat{\lambda}}.$ For each $w\in {S},$ we have
$$
  \hat{\lambda}(x_{\hat{u}}-x_{w})\geqslant x_{u_1}+x_{u_2}+(t-d_{T}(w))x_{v_t}\geqslant \left(\frac{2}{\hat{\lambda}-1}+\frac{t-d_{T}(w)}{\hat{\lambda}}\right)x_{\hat{u}}.
$$
Hence $x_w\leqslant (1-\frac{2}{\hat{\lambda}(\hat{\lambda}-1)}-\frac{t-d_{T}(w)}{{\hat{\lambda}}^2})x_{\hat{u}}$ {if $w\in S$.} By a similar discussion as \eqref{eq:3.5}, one obtains
$$
  \tilde{e}x_{\hat{u}}<\left(\frac{2}{\hat{\lambda}-1}+\sum_{w\in S}d_{T}(w)\left(1-\frac{2}{\hat{\lambda}(\hat{\lambda}-1)}-\frac{t-d_{T}(w)}{{\hat{\lambda}}^2}\right)\right)x_{\hat{u}}.
$$
That is,
\begin{align*}
 0 &<\frac{2}{\hat{\lambda}-1}-\sum_{w\in S}\left(\frac{2}{\hat{\lambda}(\hat{\lambda}-1)}+\frac{t-d_{T}(w)}{{\hat{\lambda}}^2}\right)d_{T}(w)
 =\frac{2\hat{\lambda}-2\tilde{e}}{\hat{\lambda}(\hat{\lambda}-1)}-\frac{t\tilde{e}-\sum_{w\in S}d^2_{T}(w)}{{\hat{\lambda}}^2}\\
 &\leqslant \frac{1}{\hat{\lambda}}\left(\frac{2\hat{\lambda}-2\tilde{e}}{\hat{\lambda}-1}-\frac{t\tilde{e}-\tilde{e}^2}{\hat{\lambda}}\right)
 =\frac{1}{{\hat{\lambda}}^2(\hat{\lambda}-1)}({2{\hat{\lambda}}^2-2\hat{\lambda}\tilde{e}}-{t(\hat{\lambda}-1)\tilde{e}+\tilde{e}^2(\hat{\lambda}-1)}).
\end{align*}
It follows that
\begin{align*}
    0&<2{\hat{\lambda}}^2-2\hat{\lambda}\tilde{e}-t(\hat{\lambda}-1)\tilde{e}+\tilde{e}^2(\hat{\lambda}-1)\\
    &\leqslant 2{\hat{\lambda}}^2-2\hat{\lambda}\tilde{e}-(\tilde{e}(\tilde{e}-1)-2)(\hat{\lambda}-1)\tilde{e}+\tilde{e}^2(\hat{\lambda}-1)\\
    &=2{\hat{\lambda}}^2-\tilde{e}^2(\tilde{e}-2)(\hat{\lambda}-1)-2\tilde{e}\\
    &< 2{\hat{\lambda}}^2-\tilde{e}^2(\tilde{e}-2)^2.
\end{align*}
Thus, ${2}{\hat{\lambda}}^2>{\tilde{e}^2(\tilde{e}-2)^2}.$ Recall that $\hat{\lambda}<\sqrt{m}.$ Therefore, $2(t+3+\tilde{e})>\tilde{e}^2(\tilde{e}-2)^2,$ i.e., $t>{\frac{\tilde{e}^4}{2}-2\tilde{e}^3+2\tilde{e}^2-\tilde{e}-3}\geqslant \tilde{e}-2.$ Hence there are at least two pendant vertices, say $v_{t-1}$ and $v_t,$ in $N_0(\hat{u}).$ Since $d(w)\geqslant 2$ for each vertex $w\in {S},$ we know that $w_s$ is adjacent to at most $\frac{\tilde{e}}{s}$ vertices in ${T}$ and $s\leqslant \frac{\tilde{e}}{2}.$ Recall that $x_{w_1}<\frac{x_{\hat{u}}}{\hat{\lambda}-1}.$ Hence 
$$
  \hat{\lambda}x_{v_1}\leqslant x_{\hat{u}}+sx_{w_1}\leqslant \left(1+\frac{s}{\hat{\lambda}-1}\right)x_{\hat{u}}\ \ \text{and}\ \ \hat{\lambda}x_{w_s}\leqslant \frac{\tilde{e}}{s}x_{v_1}.
$$
Then
$$
  x_{v_1}\leqslant \frac{\hat{\lambda}+s-1}{{\hat{\lambda}}(\hat{\lambda}-1)}x_{\hat{u}}\ \ \text{and}\ \ x_{w_s}\leqslant \frac{\tilde{e}}{s}\frac{\hat{\lambda}+s-1}{{\hat{\lambda}}^2(\hat{\lambda}-1)}x_{\hat{u}}.
$$
It follows that
$$
  x_{v_1}x_{w_s}\leqslant \frac{\tilde{e}(\hat{\lambda}+s-1)^2}{s{\hat{\lambda}}^3(\hat{\lambda}-1)^2}x^2_{\hat{u}}.
$$
Using Mathematica 9.0 \cite{W2012} gives us ${\tilde{e}(\hat{\lambda}+s-1)^2}<{s{\hat{\lambda}}(\hat{\lambda}-1)^2}$ if $\hat{\lambda}\geqslant \frac{\tilde{e}(\tilde{e}-2)}{\sqrt{2}}$ and $\tilde{e}\geqslant 2s\geqslant 4.$ Hence $x_{v_1}x_{w_s}<\frac{x^2_{\hat{u}}}{{\hat{\lambda}}^2}= x_{v_{t-1}}x_{v_t}.$ Let $G_2=\hat{G}-v_1w_s+v_{t-1}v_t.$ Clearly, $G_2$ is a non-bipartite graph in  $\mathcal{G}(m,\{\theta_{1,2,2},\theta_{1,2,3}\}),$ {i.e., $G_2\in \mathfrak{G}(m,\{\theta_{1,2,2},\theta_{1,2,3}\}).$} Furthermore, by Lemma~\ref{lem1.1}, one has $\lambda(G_2)>\hat{\lambda},$ which contradicts the choice of $\hat{G}.$ Hence $S=\emptyset$ and $\hat{G}\cong S_m^1.$

This completes the proof.
\end{proof}
\end{document}